\documentclass[12pt]{article}
\title{Flops and Poisson deformations of symplectic varieties}
\author{Yoshinori Namikawa}
\date{ }
\usepackage{amsmath,amssymb,amsthm, amscd}
\chardef\bslash=`\\

\def\cit{{\mathbb C}}

\def\0{{\mathcal O}}
\def\g{{\mathfrak g}}

\begin{document}
\maketitle

\section{Introduction} 
In [Na] we have dealt with a deformation of a 
{\em projective} symplectic variety. This paper, 
on the contrary, deals with a deformation of 
a {\em local} symplectic variety. More exactly, 
we mean by a local symplectic variety, a normal 
variety $X$ satisfying 
\begin{enumerate} 
\item there is a birational projective morphism 
from $X$ to an affine normal variety $Y$, 
\item there is an everywhere non-degenerate 
d-closed 2-form $\omega$ on the regular part 
$U$ of $X$ such that, for any resolution 
$\pi: \tilde{X} \to X$ with $\pi^{-1}(U) \cong U$, 
$\omega$ extends to a regular 2-form on $\tilde{X}$. 
\end{enumerate} 
In the remainder, we call such a variety a {\em convex} 
symplectic variety. A convex symplectic variety has 
been studied in [K-V], [Ka 1] and [G-K]. 
One of main difficulties we meet is the fact that  
tangent objects 
$\mathbf{T}^1_X$ and $\mathbf{T}^1_Y$ are {\em not} 
finite dimensional, since $Y$ may possibly have 
non-isolated singularities; hence the usual deformation theory 
does not work well. Instead, in [K-V], [G-K], they 
introduced a Poisson scheme and studied a {\em Poisson 
deformation} of it. A Poisson deformation is 
the deformation of the pair of a scheme itself 
and a Poisson structure on it. When $X$ is a convex 
symplectic variety, $X$ admits a natural Poisson 
structure induced from a symplectic 2-form $\omega$; hence 
one can consider its Poisson deformations. Then they 
are controlled by the {\em Poisson cohomology}. The 
Poisson cohomology has been extensively 
studied by Fresse [Fr 1], [Fr 2]. 
In some good cases, it can be described by well-known 
topological data (Corollary 10). The first application of the Poisson 
deformation theory is the following two results:    
\vspace{0.15cm}

{\bf Corollary 25}. {\em Let $Y$ be an affine 
symplectic variety with a good $\mathbf{C}^*$-action 
and assume that the Poisson structure of $Y$ is positively 
weighted.  
Let 
$$ X \stackrel{f}\rightarrow Y \stackrel{f'}\leftarrow 
X'$$ 
be a diagram such that, } 
\begin{enumerate}
\item $f$ (resp. $f'$) is a crepant, birational, projective 
morphism. 
\item $X$ (resp. $X'$) has only terminal singularities. 
\item $X$ (resp. $X'$) is {\bf Q}-factorial. 
\end{enumerate} 
{\em Then both $X$ and $X'$ have locally trivial deformations 
to an affine variety $Y_t$ obtained as a Poisson  
deformation of $Y$. In particular, $X$ and $X'$ have the same 
kind of singularities.} 
\vspace{0.12cm} 

A typical situation of Corollary 25 is a symplectic flop. 
At this moment, we need 
the ``good $\mathbf{C}^*$ contition" to make sure the 
existence of an algebraization of certain formal Poisson 
deformation. For the exact definition of a good 
$\mathbf{C}^*$-action, see Appendix. 
But even if $Y$ does not have such an action, 
one can prove:  
\vspace{0.15cm}

{\bf Corollary 31}. {\em Let $Y$ be an affine symplectic 
variety. 
Let 
$$ X \stackrel{f}\rightarrow Y \stackrel{f'}\leftarrow 
X'$$ 
be a diagram such that, } 
\begin{enumerate}
\item $f$ (resp. $f'$) is a crepant, birational, projective 
morphism. 
\item $X$ (resp. $X'$) has only terminal singularities. 
\item $X$ (resp. $X'$) is {\bf Q}-factorial. 
\end{enumerate}
{\em If $X$ is smooth, then $X'$ is smooth.}
\vspace{0.2cm}

The proofs of Corollaries 25 and 31 are essentially based 
on [Ka 1], where he proved that the smoothness is preserved 
in a symplectic flop under certain assumptions. 
Corollaries 25 and 31 are local versions of 
Corollary 1 of [Na]. More general facts can be found in 
Corollary 30.   

The following is the second application: 
\vspace{0.15cm}
 
{\bf Corollary 28}. {\em Let $Y$ be an affine 
symplectic variety with a good $\mathbf{C}^*$-action. 
Assume that the Poisson structure of $Y$ is positively weighted, and 
$Y$ has only terminal singularities.  
Let 
$ f: X \rightarrow Y $ 
be a crepant, birational, projective 
morphism such that $X$ has only terminal singularities and 
such that $X$ is {\bf Q}-factorial. 
Then the following are equivalent.} 
\vspace{0.12cm} 

(a) {\em $X$ is non-singular.} 

(b) {\em $Y$ is smoothable by a Poisson 
deformation.} 
\vspace{0.12cm}

In the proof of Corollary 28, we observe that 
the pro-representable hulls (= formal Kuranishi spaces) 
of the Poisson deformations of $X$ and $Y$ are isomorphic.  
Here we just use the assumption that $Y$ has only terminal 
singularities. Thus, any formal Poisson 
deformation of $Y$ is obtained from that of $X$ by the 
contraction map; this makes it possible for us to obtain 
(a) from (b). But, what we really want, is just 
that the formal Kuranishi space for $X$ 
{\em dominates} that for $Y$. The author believes that 
this would be true if $Y$ does not have terminal 
singularities. So our final goal would be 
the following conjecture: 
\vspace{0.15cm}

{\bf Conjecture}. {\em 
Let $Y$ be an affine 
symplectic variety with a good $\mathbf{C}^*$-action. 
Assume that the Poisson structure of 
$Y$ is positively weighted. Then the following are 
equivalent.} 
\vspace{0.15cm}

(1) {\em $Y$ has a crepant projective resolution.} 
\vspace{0.15cm} 

(2) {\em $Y$ has a smoothing by a Poisson deformation.}
\vspace{0.2cm}

The contents of this paper are as follows. 
In \S 2 we introduce the Poisson cohomology of a 
Poisson algebra according 
to Fresse [Fr 1], [Fr 2]. In Propositions 5, we shall 
prove that a Poisson deformation of a Poisson algebra is 
controlled by the Poisson cohomology. In particular, 
when the Poisson algebra is smooth, the Poisson cohomology 
is computed by the Lichnerowicz-Poisson complex. 
Since the Lichnerowicz-Poisson complex is defined also 
for a smooth Poisson scheme, one can define the 
Poisson cohomology for a smooth Poisson scheme.  
In \S 3, we restrict ourselves to the Poisson structures 
attached to a convex symplectic variety $X$. When $X$ is 
smooth, the Poisson cohomology can be identified with 
the truncated De Rham cohomology (Proposition 9). 
When $X$ has only terminal singularities, 
its Poisson deformations are the same as those of  
the regular locus $U$ of $X$. Thus the Poisson deformations 
of $X$ are controlled by the truncated De Rham cohomology 
of $U$. Theorem 14 and Corollary 15 
assert that, the Poisson deformation functor of a convex 
symplectic variety with terminal singularities, has a 
pro-representable hull and it is unobstructed. 
These are more or less already known. 
But we reproduce them here so that they fit our aim and our context. 
(see also [G-K], Appendix). 
Kaledin's twistor 
deformation is also easily generalized to our singular case; 
but this generalization is very useful in the proof of 
Corollary 25.  
At the end of this section we shall prove the following 
two key results:
\vspace{0.12cm}

{\bf Theorem 17}. {\em Let $X$ be 
a convex symplectic variety with terminal 
singularities. Let $(X, \{\;, \;\})$ 
be the Poisson structure induced by 
the symplectic form on the regular part. 
Assume that $X^{an}$ is 
{\bf Q}-factorial. Then any Poisson 
deformation of $(X, \{\; , \;\})$ is 
locally trivial as a flat deformation 
(after forgetting Poisson structure).} 
\vspace{0.15cm}

{\bf Theorem 19}. {\em Let $X$ be a convex symplectic 
variety with terminal singularities. Let $L$ 
be a (not necessarily ample) line bundle on $X$. 
Then the twistor deformation $\{X_n\}_{n \geq 1}$ of 
$X$ associated 
with $L$ is locally trivial 
as a flat deformation.} 
\vspace{0.15cm}

\S 4 deals with a convex symplectic 
variety with a good $\mathbf{C}^*$-action. The main results 
of this section are Corollary 25 and Corollary 28 explained 
above. These are actually corollaries to Theorem 19 and 
Theorem 17 respectively. 
In \S 5 we consider the general case where $Y$ does not 
have a good $\mathbf{C}^*$-action. Corollary 30 is a similar 
statement to Corollary 25 in the general case; but for 
the lack of algebraizations, 
it is not clear, at this moment, how the singularities of 
$X'$ are related with those of $X$. Finally, we shall 
prove Corollary 31 explained above.  In \S 6 one can find 
a concrete example of a Poisson deformation (Example 32). 
Example 33 is an example of a singular symplectic flop.   
The final section is Appendix, where some well-known 
results on good $\mathbf{C}^*$-actions are proved. 
The main result of Appendix is Corollary A.10. 
For a non-compact variety, the analytic 
category and the algebraic category are usually quite 
different. However, Corollary A.10 asserts that when we 
have a good $\mathbf{C}^*$-action, they are well-combined.  
  
The author would like to thank A. Fujiki for the discussion 
on Lemma A.8 in the Appendix, and D. Kaledin for pointing 
out that Corollary 30 is not sufficient for us to claim that $X$ 
and $X'$ have the same kind of singularities. 
\vspace{0.15cm}

\section{Poisson deformations} 
(i) {\bf Harrison cohomology}:  
Let $S$ be a commutative 
{\bf C}-algebra and let $A$ be a commutative $S$-algebra. 
Let $S_n$ be the $n$-th symmetric 
group. Then $S_n$ acts from the left hand side on the 
$n$-tuple tensor product $A\otimes_S...\otimes_SA$ as 
$$\pi(a_1 \otimes ... \otimes a_n) := a_{\pi^{-1}(1)} \otimes 
... \otimes a_{\pi^{-1}(n)},$$ 
where $\pi \in S_n$. This action extends naturally to an action 
of the group algebra $\mathbf{C}[S_n]$ on $A\otimes_S...\otimes_SA$. 
For $0 < r < n$, an element 
$\pi \in S_n$ is called a 
{\em pure shuffle} of type $(r, n-r)$ if $\pi(1) < ... < \pi(r)$ 
and $\pi(r+1) <... < \pi(n)$. Define an element $s_{r,n-r} \in 
\mathbf{C}[S_n]$ by  
$$ s_{r,n-r} := \Sigma \mathrm{sgn}(\pi)\pi, $$ 
where the sum runs through all pure shuffles of type $(r, n-r)$. 
Let $N$ be the $S$-submodule of $A\otimes_SA ...\otimes_SA$ 
generated by all elements 
$$\{s_{r,n-r}(a_1 \otimes ... \otimes a_n)\}_{0 < r < n, \, 
a_i \in A}. $$
Define $\mathrm{ch}_n(A/S) := (A\otimes_SA ...\otimes_SA)/N$. 
Let $M$ be an $A$-module. Then the {\em Harrison chain} 
$\{\mathrm{ch}_{\cdot}(A/S; M)\}$ is defined as follows: 
\vspace{0.12cm}  

\begin{enumerate} 
\item $\mathrm{ch}_n(A/S; M) := \mathrm{ch}_n(A/S)\otimes_SM$ 
\item the boundary map $\partial_n : \mathrm{ch}_n(A/S; M) \to 
\mathrm{ch}_{n-1}(A/S; M)$ is defined by 
 $\partial_n(a_1 \otimes ... \otimes a_n \otimes m) :=$  
$$a_1 \otimes ... \otimes a_nm + \Sigma_{1 \leq i \leq n-1}
(-1)^{n-i}a_1 \otimes ... \otimes a_ia_{i+1}\otimes ... \otimes a_n 
\otimes m + (-1)^n a_2 \otimes ... \otimes a_n \otimes a_1m. $$ 
\end{enumerate}  

We define the $n$-th {\em Harrison homology} $\mathrm{Har}_n(A/S; M)$ 
just as the $n$-th homology of $\mathrm{ch}_{\cdot}(A/S; M)$.  

The {\em Harrison cochain} $\{\mathrm{ch}^{\cdot}(A/S;M)\}$ 
is defined as follows 
\vspace{0.12cm} 
     
\begin{enumerate} 
\item $\mathrm{ch}^n(A/S; M) := \mathrm{Hom}_S(\mathrm{ch}_n(A/S), M)$ 
\item the coboundary map $d^n: \mathrm{ch}^n(A/S; M) \to 
\mathrm{ch}^{n+1}(A/S; M)$ is defined by 
 $(d^nf)(a_1 \otimes ... \otimes a_{n+1}):= (-1)^{n+1} 
a_1 f(a_2\otimes ... \otimes a_{n+1})$ $$+ \Sigma_{1 \leq i \leq n}
(-1)^{n+1-i}f(a_1 \otimes ... \otimes a_ia_{i+1}\otimes ... \otimes a_{n+1})  
 + f(a_1 \otimes ... \otimes a_n)a_{n+1}. $$ 
\end{enumerate}

We define the $n$-th {\em Harrison cohomology} $\mathrm{Har}^n(A/S; M)$ 
as the $n$-th cohomolgy of $\mathrm{ch}^{\cdot}(A/S; M)$. 
\vspace{0.15cm} 

{\bf Example 1}. Assume that $S = \mathbf{C}$. 
Then $\mathrm{ch}_2(A/\mathbf{C}; A) = \mathrm{Sym}^2_{\mathbf{C}}(A)
\otimes_{\mathbf{C}}A$ and $\mathrm{ch}_1(A/\mathbf{C}; A) 
= A \otimes_{\mathbf{C}}A$. The boundary map $\partial_2$ is defined 
as 
$$ \partial_2([a_1\otimes a_2]\otimes a) := a_1\otimes a_2a - 
a_1a_2 \otimes a + a_2 \otimes a_1a.$$ 
We see that $\mathrm{im}(\partial_2)$ is a right $A$-submodule 
of $A\otimes_{\mathbf{C}}A$.  Let $I \subset A \otimes_{\mathbf{C}}A$ be 
the ideal generated by all elements of the form $a\otimes b - 
b\otimes a$ with $a, b \in A$. Then we have a homomorphism of 
right $A$-modules $I \to (A\otimes_{\mathbf{C}}A)/\mathrm{im}(\partial_2)$. 
One can check that its kernel coincides with $I^2$. Hence we have 
$$\Omega^1_{A/\mathbf{C}} \cong \mathrm{Har}_1(A/\mathbf{C}; A).$$ 
In fact, the Harrison chain $\mathrm{ch}_{\cdot}(A/\mathbf{C}; A)$ 
is quasi-isomorphic to the cotangent complex $L^{\cdot}_{A/\mathbf{C}}$ for 
a $\mathbf{C}$-algebra $A$ (cf. [Q]). 
\vspace{0.15cm}

Let $A$ and $S$ be the same as above. We put $S[\epsilon] 
:= S\otimes_{\mathbf{C}}\mathbf{C}[\epsilon]$, where 
$\epsilon^2 = 0$. Let us consider the set of all 
$S[\epsilon]$-algebra structures of the $S[\epsilon]$-module 
$A\otimes_S S[\epsilon]$ such that they induce the original 
$S$-algebra $A$ if we take the tensor product of  
$A\otimes_S S[\epsilon]$ and $S$ over $S[\epsilon]$. 
We say that two elements of this set are {\em equivalent} if 
and only if there is an isomorphism of $S[\epsilon]$-algebras 
between them which induces the identity map of $A$ over 
$S$. We denote by $D(A/S, S[\epsilon])$ the set of such equivalence 
classes. Fix an $S[\epsilon]$-algebra structure $(A\otimes_S 
S[\epsilon], *)$. Here $*$ just means the corresponding ring 
structure. Then we define $\mathrm{Aut}(*, S)$ to be the 
set of all $S[\epsilon]$-algebra automorphisms of 
$(A\otimes_S S[\epsilon], *)$ which induces the identity map 
of $A$ over $S$.  
\vspace{0.15cm}

{\bf Proposition 2}. {\em Assume that 
$A$ is a free $S$ module.} 

(1) {\em There is a one-to-one correspondence 
between $\mathrm{Har}^2(A/S; A)$ and \\
$D(A/S, S[\epsilon])$.} 
\vspace{0.12cm} 

(2) {\em There is a one-to-one correspondence between 
$\mathrm{Har}^1(A/S; A)$ and \\ 
$\mathrm{Aut}(*, S)$.}  
\vspace{0.12cm}

{\em Proof}. We shall only give a proof to (1). 
The proof of (2) is left to the readers. 
Denote by $*$ a ring structure on 
$A\otimes_S S[\epsilon] = A \oplus A\epsilon$. For 
$a, b \in A$,  
write $$a * b = ab + \epsilon \phi (a,b)$$ with some 
$\phi : A \times A \to A$. The multiplication 
of an element of $S[\epsilon]$ and an element of 
$A\otimes_S S[\epsilon]$ should coincides with the action 
of $S[\epsilon]$ as the $S[\epsilon]$-module; hence 
$a * \epsilon = a\epsilon$ and $a\epsilon * \epsilon = 0$. 
Then $$a * (b\epsilon) = a * (b * \epsilon) = (a * b)* \epsilon$$ 
$$= \{ab + \phi (a,b)\epsilon\}*\epsilon = ab \epsilon 
+ \phi (a,b)* (\epsilon * \epsilon) = ab\epsilon.$$ 
Similarly, we have $(a\epsilon)*(b\epsilon) = 0$.     
Therefore, $*$ is determined completely by $\phi$. 
By the commutativity of $*$, $\phi \in 
\mathrm{Hom}_S(\mathrm{Sym}^2_S(A), A)$. 
By the associativity: $(a*b)*c = a*(b*c)$, we get 
$$\phi (ab,c) + c\phi(a,b) = \phi(a,bc) + 
a\phi(b,c).$$ This condition is equivalent to that 
$\phi \in \mathrm{Ker}(d^2)$, where $d^2$ is the 
2-nd coboundary map of the Harrison cochain. 
Next let us observe when two ring structures 
$*$ and $*'$ are equivalent. As above, we write 
$a*b = ab + \epsilon \phi(a,b)$ and 
$a*'b = ab + \epsilon \phi'(a,b)$. 
Assume that a map $\psi: A \oplus A\epsilon 
\to A \oplus A\epsilon$ gives an equivalence. 
Then, for $a \in A$, write $\psi(a) = a + f(a)\epsilon$ 
with some $f: A \to A$. One can show that  
$\psi(a\epsilon) = a\epsilon$. 
Since $\psi(a)*'\psi(b) = \psi(a*b)$, 
we see that 
$$\phi'(a,b) -\phi(a,b) = f(ab) -af(b) -bf(a).$$ 
This implies that $\phi' - \phi \in \mathrm{im}(d^1)$.  
\vspace{0.15cm}

{\bf Remark 3}. Assume that $S$ is an Artinian ring and 
$A$ is flat over $S$. Then $A$ is a free $S$-module and 
for any flat extension $A'$ of $A$ over $S[\epsilon]$, 
$A' \cong A\otimes_S S[\epsilon]$ as an $S[\epsilon]$-module. 
\vspace{0.15cm}     
 
(ii) {\bf Poisson cohomology} (cf. [Fr 1,2]): 
Let $A$ and $S$ be the same as (i). 
Assume that $A$ is a free $S$-module. Let us consider the 
graded free $S$-module 
$\mathrm{ch}_{\cdot}(A/S) := \oplus_{0 < m}\mathrm{ch}_m(A/S)$ and take 
its super-symmetric algebra $\mathcal{S}(\mathrm{ch}_{\cdot}(A/S))$. 
By definition, $\mathcal{S}(\mathrm{ch}_{\cdot}(A/S))$ is 
the quotient of the tensor algebra $T(\mathrm{ch}_{\cdot}(A/S)) := 
\oplus_{0 \le n}(\mathrm{ch}_{\cdot}(A/S))^{\otimes n}$ by the 
two-sided ideal  
$M$ generated by the elements of the form: 
$a \otimes b - (-1)^{pq}b \otimes a$, where $a \in \mathrm{ch}_p(A/S)$ 
and $b \in \mathrm{ch}_q(A/S)$. We denote by 
$\bar{\mathcal{S}}(\mathrm{ch}_{\cdot}(A/S))$ the truncation of the degree 
$0$ part. In other words, 
$$\bar{\mathcal{S}}(\mathrm{ch}_{\cdot}(A/S)) := 
\oplus_{0 < n}(\mathrm{ch}_{\cdot}(A/S))^{\otimes n}/M.$$      
        
Now let us consider the graded $A$-module 
$$\bar{\mathcal{S}}(\mathrm{ch}_{\cdot}(A/S))\otimes_{S}A 
:= \mathrm{ch}_{\cdot}(A/S)\otimes_{S}A \oplus 
(\mathcal{S}^2(\mathrm{ch}_{\cdot}(A/S))\otimes_{S}A) \oplus 
 ...$$ 
The Harrison boundary maps $\partial$ on 
$\mathrm{ch}_{\cdot}(A/S)\otimes_{S}A$ naturally extends to 
those on $\mathcal{S}^n(\mathrm{ch}_{\cdot}(A/S))\otimes_{S}A$. 
In fact, for $a_i \in \mathrm{ch}_{p_i}(A/S)\otimes_{S}A$, 
$1 \leq i \leq n$, denote by $a_1 \cdot \cdot \cdot a_n\in 
\mathcal{S}^n_A(\mathrm{ch}_{\cdot}(A/S)\otimes_{S}A)$  
their super-symmetric product. We then define $\partial$ 
inductively as 
$$\partial(a_1 ... a_n) := \partial(a_1)a_2... a_n  
 + (-1)^{p_1}a_1\cdot\partial(a_2 ... a_n).$$  
In this way, each $\mathcal{S}^n_A(\mathrm{ch}_{\cdot}(A/S)\otimes_{S}A) 
= \mathcal{S}^n(\mathrm{ch}_{\cdot}(A/S))\otimes_{S}A$ 
becomes a chain complex. 
By taking the dual, 
$$\mathrm{Hom}_A
(\mathcal{S}^n(\mathrm{ch}_{\cdot}(A/S))\otimes_{S}A, A)$$  
$$= \mathrm{Hom}_S(\mathcal{S}^n(\mathrm{ch}_{\cdot}(A/S)), A)$$ 
becomes a cochain complex: 
  
\begin{equation}
\begin{CD}
@A{d}AA @. @. @. \\
\mathrm{Hom}(\mathrm{ch}_3, A) @. {...} @. {...} @.\\
@A{d}AA  @A{d}AA  @A{d}AA  @. \\
\mathrm{Hom}(\mathrm{ch}_2, A) @. \;\;\mathrm{Hom}(\mathrm{ch}_2\otimes
\mathrm{ch}_1, A)  @. {...} @.\\ 
@A{d}AA  @A{d}AA @A{d}AA @. \\ 
\mathrm{Hom}(\mathrm{ch}_1, A) @. 
\mathrm{Hom}(\wedge^2\mathrm{ch}_1, A) @. 
\mathrm{Hom}(\wedge^3\mathrm{ch}_1, A) @. \;\;{...}  
\end{CD}
\end{equation}
\vspace{0.15cm}
  
Here we abbreviate $\mathrm{ch}_i(A/S)$ by $\mathrm{ch}_i$ 
and $\mathrm{Hom}_S(...)$ by $\mathrm{Hom}(...)$.   
We want to make the diagram above into a double complex 
when $A$ is a Poisson $S$-algebra. 
\vspace{0.12cm}

{\bf Definition}. A Poisson $S$-algebra $A$ is a commutaive 
$S$-algebra with an $S$-linear map 
$$\{\; , \;\}: \wedge^2_{S}A \to A$$ 
such that 
\begin{enumerate}
\item $\{a, \{b, c\}\} + \{b, \{c, a\}\} + \{c, \{a, b\}\} = 0$ 
\item $\{a, bc\} = \{a, b\}c + \{a, c\}b.$
\end{enumerate}

We assume now that $A$ is a Poisson $S$-module such 
that $A$ is a free $S$-module. We put  $\bar{T}_S(A) 
 := \oplus_{0 < n}(A)^{\otimes n}$.  
We shall introduce an $S$-bilinear bracket product $$[\; , \;] :         
\bar{T}_S(A) \times 
\bar{T}_S(A) \to 
\bar{T}_S(A)$$ 
in the following manner. Take two elements from $\bar{T}_S(A)$: 
$f = f_1\otimes ... \otimes f_p$   
and $g = g_1\otimes ... \otimes g_q$. Here 
each $f_i$ and each $g_i$ are elements of $A$. 
Let $\pi \in S_{p+q}$ be a pure shuffle of type 
$(p, q)$. For the convention, we put $f_{i+p} := 
g_i$. Then the shuffle product is defined as 
$$ f\cdot g := \Sigma\mathrm{sgn}(\pi)f_{\pi(1)}\otimes ... 
\otimes f_{\pi(p+q)},$$ where the sum runs through all 
pure shuffle of type $(p,q)$. 
For each term of the sum (which is indexed by 
$\pi$), let $I_{\pi}$ be the set of all $i$ such that 
$\pi(i) \le p$ and $\pi(i+1) \ge p+1$ (which implies 
that $f_{\pi(i+1)} = g_{\pi(i+1)-p}$). Then we define 
$[f,g]$ as 
$$\Sigma\mathrm{sgn}(\pi)(\Sigma_{i \in I_{\pi}}(-1)^{i+1}
f_{\pi(1)}\otimes ... \otimes \{f_{\pi(i)}, f_{\pi(i+1)}\} 
\otimes ... \otimes f_{\pi(p+q)}). $$  

The bracket $[\;, \;]$ induces that on $\mathrm{ch}_{\cdot}(A/S)$ 
by the quotient map $\bar{T}_S(A) \to \mathrm{ch}_{\cdot}(A/S)$. 
By abuse of notation, we denote by $[\;, \;]$ the induced bracket. 
We are now in a position to define coboundary maps 
$$\delta: \mathrm{Hom}_S(\mathcal{S}^{s-1}(\mathrm{ch}_{\cdot}
(A/S)), A) \to \mathrm{Hom}_S(\mathcal{S}^s(\mathrm{ch}_{\cdot}
(A/S)), A)$$ so that $\mathrm{Hom}_S(\bar{\mathcal{S}}
(\mathrm{ch}_{\cdot}(A/S)), A)$ is made into a double complex 
together with $d$ already defined. We 
take an element of the form $x_1\cdot\cdot\cdot x_s$ from $\mathcal{S}^s
(\mathrm{ch}_{\cdot}(A/S))$ with each $x_i$ being a homogenous 
element of $\mathrm{ch}_{\cdot}(A/S)$.  

For $f \in \mathrm{Hom}_S(\bar{\mathcal{S}}^{s-1}
(\mathrm{ch}_{\cdot}(A/S)), A)$, we define \vspace{0.12cm} 

$$\delta(f)(x_1...x_s) := \sum_{1 \le i \le s}(-1)^{\sigma(i)}
\overline{[x_i, f(x_1\cdot\cdot\cdot\breve{x_i}\cdot\cdot\cdot x_s)]}$$ 
$$- \sum_{i < j}(-1)^{\tau(i,j)}
f([x_i, x_j]\cdot\cdot\cdot\breve{x_i}\cdot\cdot\cdot \breve{x_j}
\cdot\cdot\cdot x_s).$$ \vspace{0.12cm}
   
Here $\overline{[\; , \;]}$ is the composite of $[\; , \;]$ and 
the truncation map $\mathrm{ch}_{\cdot}(A/S) \to \mathrm{ch}_1(A/S)(=A)$. 
Moreover, 
$$\sigma(i) := 
\mathrm{deg}(x_i)
\cdot (\mathrm{deg}(x_1) + ... +  
\mathrm{deg}(x_{i-1}))$$ 
and 
$$ \tau(i,j) := \mathrm{deg}(x_i)(\mathrm{deg}x_1 + ... + 
\mathrm{deg}x_{i-1})$$ 
$$+ \mathrm{deg}(x_j)(\mathrm{deg}(x_1) + ... + 
\breve{\mathrm{deg}(x_i)}  
+ ... + \mathrm{deg}(x_{j-1})).$$
   
We now obtain a double complex $(\mathrm{Hom}_S
(\bar{\mathcal{S}}(\mathrm{ch}_{\cdot}(A/S)), A), d, \delta)$. 
The $n$-th {\em Poisson cohomology} 
$\mathrm{HP}^n(A/S)$ for a Poisson $S$-algebra 
$A$ is the $n$-th cohomology of the total complex (by 
$d + \delta$) of 
this double complex.  

\begin{equation}
\begin{CD}
@A{d}AA @. @. @. \\
\mathrm{Hom}(\mathrm{ch}_3, A) @>{\delta}>> {...} @. {...} @.\\
@A{d}AA  @A{d}AA  @A{d}AA  @. \\
\mathrm{Hom}(\mathrm{ch}_2, A) @>{\delta}>> \mathrm{Hom}(\mathrm{ch}_2\otimes
\mathrm{ch}_1, A)  @>{\delta}>> {...} @.\\ 
@A{d}AA  @A{d}AA @A{d}AA @. \\ 
\mathrm{Hom}(\mathrm{ch}_1, A) @>{\delta}>> 
\mathrm{Hom}(\wedge^2\mathrm{ch}_1, A) @>{\delta}>> 
\mathrm{Hom}(\wedge^3\mathrm{ch}_1, A) @. \;\;{...}  
\end{CD}
\end{equation}
\vspace{0.15cm}  
 
{\bf Example 4}. We shall calculate $\delta$ explicitly in a 
few cases. As in the diagram above, we abbreviate $\mathrm{Hom}_S$ 
by $\mathrm{Hom}$, and $\mathrm{ch}_i(A/S)$ by $\mathrm{ch}_i$. 
\vspace{0.12cm}

(i) Assume that $f \in \mathrm{Hom}(\mathrm{ch}_1, A)$. 
$$\delta(f)(a \wedge b) = \{a, f(b)\} + \{f(a), b\} - f(\{a, b\}).$$ 
\vspace{0.12cm}

(ii) Assume that $\varphi \in \mathrm{Hom}(\mathrm{ch}_2, A)$. 
For $(a,b) \in \mathrm{Sym}^2_S(A)(= \mathrm{ch}_2)$, and for 
$c \in A (= \mathrm{ch}_1)$, 
$$\delta(\varphi)((a,b)\cdot c) = [(a,b), \varphi(c)] + [c, \varphi(a,b)] 
- \varphi([(a,b),c])$$ 
$$= \{c, \varphi(a,b)\} - \varphi(\{c,b\},a) - \varphi(\{c,a\},b). $$ 
\vspace{0.12cm}

(iii) Assume that $\psi \in \mathrm{Hom}(\wedge^2\mathrm{ch}_1, A)$. 
$$\delta(\psi)(a\wedge b \wedge c) = \{a, \psi(b,c)\} + \{b, \psi(c,a)\} 
+ \{c, \psi(a,b)\}$$ 
$$ + \psi(a, \{b,c\}) + \psi(b, \{c,a\}) + \psi(c, \{a,b\}).$$ 
\vspace{0.12cm}

Let $A$ be a Poisson $S$-algebra such that $A$ is a free 
$S$-module. We put $S[\epsilon] := 
S\otimes_{\mathbf{C}}\mathbf{C}[\epsilon]$, where 
$\epsilon^2 = 0$. Let us consider the set of all 
Poisson $S[\epsilon]$-algebra structures on the 
$S[\epsilon]$-module $A\otimes_SS[\epsilon]$ such 
that they induce the original Poisson $S$-algebra 
$A$ if we take the tensor product of $A\otimes_SS[\epsilon]$ 
and $S$ over $S[\epsilon]$. We say that two elements of 
this set are {equivalent} if and only if there is an 
isomorphism of Poisson $S[\epsilon]$-algebras between 
them which induces the identity map of $A$ over $S$. 
We denote by $PD(A/S, S[\epsilon])$ the set of such 
equivalence classes. Fix a Poisson $S[\epsilon]$-algebra 
structure $(A\otimes_SS[\epsilon], *, \{\; , \;\})$. 
Then we define $\mathrm{Aut}(*, \{\; , \;\}, S)$ to 
be the set of all automorphisms of Poisson 
$S[\epsilon]$-algebras of 
$(A\otimes_SS[\epsilon], *, \{\; , \;\})$ which 
induces the identity map of $A$ over $S$. 

{\bf Proposition 5}. (1) {\em There is a 
one-to-one correspondence between \\
$HP^2(A/S)$ and $PD(A/S, S[\epsilon])$.} 
\vspace{0.12cm}

(2) {\em There is a one-to-one 
correspondence between $HP^1(A/S)$ and \\
$\mathrm{Aut}(*, \{\; , \;\}, S)$.} 
\vspace{0.12cm}

{\em Proof}. (1): As explained in Proposition 2, 
giving an $S[\epsilon]$-algebra structure $*$ on 
$A \oplus A\epsilon$ is equivlaent to giving $\varphi 
\in \mathrm{Hom}_S(\mathrm{Sym}^2_S(A), A)$ with 
$d(\varphi) = 0$ such that 
$a * b = ab + \epsilon\varphi(a,b)$. 
Assume that $\{\; , \;\}_{\epsilon}$ is a Poisson bracket on 
$(A \oplus A\epsilon, *)$ which is an extension of the 
original Poisson bracket $\{\; , \;\}$ on $A$. We put 
$$\{a, b\}_{\epsilon} = \{a,b\} + \psi(a,b)\epsilon.$$ 
Since $\{a, b\epsilon\}_{\epsilon} = \{a, b*\epsilon\}_{\epsilon} 
= \{a,b\}\epsilon$ and $\{a\epsilon, b\epsilon\} = 0$, 
the Poisson structure $\{\; , \;\}_{\epsilon}$ is 
completely determined by $\psi$. By the skew-commutativity 
of $\{\; , \;\}$, $\psi \in \mathrm{Hom}_S(\wedge^2_SA, A)$. 
The equality 
$$\{a, b*c\}_{\epsilon} = \{a, b\}_{\epsilon}*c + 
\{a, c\}_{\epsilon}*b$$ 
is equivalent to the equality
$$(\star): \;\; \psi(a, bc) - c\psi(a,b) - b\psi(a,c) $$ 
$$= \varphi(\{a,b\},c) + \varphi(\{a,c\},b) 
- \{a, \varphi(b,c)\}.$$ 
The equality 
$$\{a, \{b, c\}_{\epsilon}\}_{\epsilon} + 
\{b, \{c, a\}_{\epsilon}\}_{\epsilon} + 
\{c, \{a, b\}_{\epsilon}\}_{\epsilon} = 0$$ 
is equivalent to the equality 
$$(\star\star): \psi(a, \{b,c\}) + \psi(b, \{c,a\}) + \psi(c, \{a,b\})$$  
$$+ \{a, \psi(b,c)\} + \{b, \psi(c,a)\} + 
\{c, \psi(a,b)\} = 0.$$ 
  
We claim that the equality $(\star)$ means $\delta(\varphi) 
+ d(\psi) = 0$ in the diagram: 
$$ \mathrm{Hom}(\mathrm{Sym}^2(A),A) \stackrel{\delta}\to  
\mathrm{Hom}(\mathrm{Sym}^2(A)\otimes A, A) \stackrel{d}\leftarrow 
\mathrm{Hom}(\wedge^2A, A).$$   
By Example 4, (ii), we have shown that 
$$\delta(\varphi)((a,b)\cdot c) = \{c, \varphi(a,b)\} - 
\varphi(\{c,b\},a) - \varphi(\{c,a\},b).$$ 
On the other hand, for the Harrison boundary map  
$$\partial: \mathrm{Sym}^2(A)\otimes_S A \otimes_S A \to 
\wedge^2A \otimes_S A,$$ we have 
$$\partial ((a,b)\otimes c \otimes 1) = b(a\wedge c) - 
ab \wedge c + a(b \wedge c).$$ Since $d$ is defined as the 
dual map of $\partial$, we see that 
$$d\psi((a,b)\cdot c) = \psi(c,ab) - a\psi(c,b) - b\psi(c,a).$$ 
As a consequence, we get 
$$(\delta \varphi + d \psi)((a,b)\cdot c) = \psi(c,ab) 
- a\psi(c,b) -b\psi(c,a)$$ 
$$+ \{c, \varphi(a,b)\} - \varphi(\{c,b\},a) - 
\varphi(\{c,a\}, b).$$ 
By changing $a$ and $c$ each other, we conclude that 
$\delta(\varphi) 
+ d(\psi) = 0$.  

By the equality $(\star\star)$ and Example 4, (iii), we 
see that $(\star\star)$ means $\delta(\psi) = 0$ for the 
map $\delta: \mathrm{Hom}(\wedge^2A, A) \to 
\mathrm{Hom}(\wedge^3A, A)$.    
Next, let us observe when two Poisson structures   
$(\varphi, \psi)$ and $({\varphi}', {\psi}')$ 
(on $A \oplus A\epsilon$) are equivalent. 
Assume that, for $f \in \mathrm{Hom}_S(A,A)$,  
$$ \chi_f : A \oplus A\epsilon \to A \oplus A\epsilon $$ 
gives such an equivalence between both Poisson structures, 
where $\chi_f(a) = a + f(a)\epsilon$, $\chi_f(a\epsilon) 
= a\epsilon$ for $a \in A$.    
Since $\chi_f$ gives an equivalence of $S[\epsilon]$-algebras, 
$$(\varphi' - \varphi)(a,b) = f(ab) - af(b) -bf(a) = -d(f)(a,b)$$ 
by Proposition 2. 
The map $\chi_f$ must be compatible with two Poisson 
structure: 
$$\{\chi_f(a), \chi_f(b)\}'_{\epsilon} 
= \chi_f(\{a,b\}_{\epsilon}).$$
The left hand side equals 
$$\{a,b\} + [\phi'(a,b) + \{a,f(b)\} + \{f(a),b\}]\epsilon.$$ 
The right hand side equals 
$$\{a,b\} + [f(\{a,b\}) + \psi(a,b)]\epsilon.$$ 
Thus, we have 
$$(\psi' - \psi)(a,b) = -\delta(f)(a,b),$$ and the proof 
of (1) is now complete. 
We omit the proof of (2).  
\vspace{0.15cm}

We next consider the case where $A$ is formally smooth over 
$S$.  
We put $\Theta_{A/S} := \mathrm{Hom}_A
(\Omega^1_{A/S}, A)$. 
We make $\oplus_{i>0}\wedge^i_A\Theta_{A/S}$ 
into a complex by defining the coboundary map 
$$\delta : \wedge^i\Theta_{A/S} \to \wedge^{i+1}\Theta_{A/S}$$ 
as 
$$\delta f(da_1 \wedge ... \wedge da_{i+1}) := 
\sum_{j}(-1)^{j+1}\{a_j, f(da_1 \wedge ... \wedge \breve{da_j} 
\wedge ... \wedge da_{i+1})\} $$ 
$$ - \sum_{j < k}(-1)^{j+k+1}f(d\{a_j, a_k\})\wedge da_1 \wedge 
... \wedge \breve{da_j}\wedge ... \wedge \breve{da_k} \wedge 
... \wedge da_{i+1}),$$ 
for $f \in \wedge^i\Theta_{A/S} = \mathrm{Hom}_A
(\Omega^i_{A/S}, A)$. 
This complex is called the {\em Lichnerowicz-Poisson} complex.   
One can connect this complex with our Poisson cochain 
complex $\mathcal{C}^{\cdot}(A/S)$. In fact, 
there is a map $\mathrm{ch}_1\otimes_S A  
\to \Omega^1_{A/S}$ (cf. Example 1). This map induces, for each 
$i$, 
$\wedge^i\mathrm{ch}_1(A/S)\otimes_S A  
\to \Omega^i_{A/S}$. By taking the dual, we get 
$$\wedge^i\Theta_{A/S} \to {\mathrm{Hom}}_A
(\wedge^i\mathrm{ch}_1(A/S)\otimes_S A, A) 
= \mathrm{Hom}(\wedge^i\mathrm{ch}_i(A/S), A).$$ 
By these maps, we have a map of complexes 
$$\wedge^{\cdot}\Theta_{A/S} \to \mathcal{C}^{\cdot}(A/S).$$ 

{\bf Proposition 6}. {\em For a Poisson $S$-algebra $A$, assume that     
$A$ is formally smooth over $S$ and that $A$ is a free $S$-module.  
Then $(\wedge^{\cdot}\Theta_{A/S}, \delta) 
\to (\mathcal{C}^{\cdot}(A/S), d + \delta)$ is a quasi-isomorphism.} 
\vspace{0.15cm}

For the proof of Proposition 6, see Fresse [Fr 1], Proposition 1.4.9. 
\vspace{0.15cm}

{\bf Definition}. Let $T := \mathrm{Spec}(S)$ and $X$ a $T$-scheme. 
Then $(X, \{\; , \;\})$ is a Poisson scheme over $T$ if 
$\{\; , \;\}$ is an $\mathcal{O}_T$-linear map: 
$$ \{\; , \;\}: \wedge^2_{\mathcal{O}_T}\mathcal{O}_X \to 
\mathcal{O}_X$$ 
such that, for $a,b,c \in \mathcal{O}_X$, 
\begin{enumerate}
\item $\{a, \{b,c\}\} + \{b,\{c,a\}\} + \{c,\{a,b\}\} = 0$ 
\item $\{a,bc\} = \{a,b\}c + \{a,c\}b.$ 
\end{enumerate}

We assume that $X$ is a {\em smooth} Poisson scheme over $T$, 
where $T = \mathrm{Spec}(S)$ with a local Artinian {\bf C}-algebra 
$S$ with $S/m_S = \mathbf{C}$. Then 
the Lichnerowicz-Poisson complex can be globalized \footnote{The 
definition of the Poisson cochain complex is subtle because 
the sheafication of each component of the Harrison complex 
is not quasi-coherent (cf. [G-K])} to the 
complex on $X$ $$\mathcal{LC}^{\cdot}(X/T) := 
(\wedge^{\cdot}\Theta_{X/T}, \delta).$$  
 
We define the $i$-th {Poisson cohomology} as 
$$\mathrm{HP}^i(X/T) := \mathbf{H}^i(X, \mathcal{LC}^{\cdot}(X/T)).$$ 
 
{\bf Remark 7}. 
When $X = \mathrm{Spec}(A)$, 
$\mathrm{HP}^i(X/T) = \mathrm{HP}^i(A/S)$. In fact, 
there is a spectral sequence induced from the 
stupid filtration: 
$$E^{p,q}_1 := H^q(X, \mathcal{LC}^p(X/T)) => 
\mathrm{HP}^i(X/T).$$ 
Since each $\mathcal{LC}^p(X/T)$ is quasi-coherent 
on the affine scheme $X$, $H^q(X, \mathcal{LC}^p) = 0$ 
for $q > 0$. Therefore, this spectral sequence 
degenerate at $E_2$-terms and we have 
$$\mathrm{HP}^i(X/T) = H^i(\Gamma (X, \mathcal{LC}^{\cdot})),$$ 
where the right hand side is nothing but $HP^i(A/S)$ by 
Proposition 6.   
\vspace{0.15cm} 

One can generalize Proposition 5 to smooth Poisson schemes. 
Let $S$ be an Artinian $\mathbf{C}$-algebra and put 
$T := \mathrm{Spec}(S)$. 
Let $X$ be a Poisson $T$-scheme which is smooth over $T$.  
We put $T[\epsilon] := \mathrm{Spec}S[\epsilon]$ with 
$\epsilon^2 = 0$.     
A Poisson deformation $\mathcal{X}$ of $X$ 
over $T[\epsilon]$ is a Poisson $T[\epsilon]$-algebra 
such that $\mathcal{X}$ is flat over $T[\epsilon]$ 
and there is a Poisson isomorphism 
$\mathcal{X}\times_{T[\epsilon]}T \cong X$ over 
$T$. Two Poisson deformations $\mathcal{X}$ and 
$\mathcal{X}'$ are equivalent if there is an 
isomorphism $\mathcal{X} \cong \mathcal{X}'$ as 
Poisson $T[\epsilon]$-schemes such that it 
induces the identity map of $X$ over $T$. 
Denote by $PD(X/T, T[\epsilon])$ the set of 
equivalence classes of Poisson deformations of 
$X$ over $T[\epsilon]$. For a Poisson deformation 
$\mathcal{X}$ of $X$ over $T[\epsilon]$, we denote 
by $Aut(\mathcal{X}, T)$ the set of all automorphisms 
of $\mathcal{X}$ as a Poisson $T[\epsilon]$-scheme 
such that they induce the identity map of $X$ 
over $T$.  
\vspace{0.15cm}

{\bf Proposition 8}. (1) {\em There is a one-to-one 
correspondence between $\mathrm{HP}^2(X/T)$ and 
$\mathrm{PD}(X/T, T[\epsilon])$.} 
\vspace{0.12cm}

(2) {\em For a Poisson deformation $\mathcal{X}$ of 
$X$ over $T[\epsilon]$, there is a one-to-one correspondence 
between $\mathrm{HP}^1(X/T)$ and $Aut(\mathcal{X}, T)$.} 
\vspace{0.12cm}

{\em Proof}. We only prove (1).  For an affine open covering 
$\mathcal{U} := \{U_i\}_{i \in I}$ of $X$, construct 
a double complex $\Gamma(\mathcal{LC}^{\cdot}(\mathcal{U}, X/T))$ 
as follows: 
 
\begin{equation}
\begin{CD}
@A{\delta}AA  @A{\delta}AA @.  \\ 
\prod_{i_0}\Gamma(\mathcal{LC}^2(U_{i_0}/T)) @>>> 
\prod_{i_0, i_1}\Gamma(\mathcal{LC}^2(U_{i_0i_1}/T)) @>>> 
{...} \\ 
@A{\delta}AA  @A{\delta}AA @.  \\ 
\prod_{i_0}\Gamma(\mathcal{LC}^1(U_{i_0}/T)) @>>> 
\prod_{i_0,i_1}\Gamma(\mathcal{LC}^1(U_{i_0,i_1}/T)) 
@>>> {...} 
\end{CD}
\end{equation}
  
Here the horizontal maps are \v{C}ech coboundary maps. 
Since each $\mathcal{LC}^p$ is quasi-coherent, 
one can calculate the Poisson cohomology by 
the total complex associated with this double 
complex: 
$$\mathrm{HP}^i(X/T) = H^i(\Gamma(\mathcal{LC}^{\cdot}(\mathcal{U},X/T))).$$ 
An element $\zeta \in \mathrm{HP}^2(X/T)$ corresponds to 
a 2-cocycle 
$$ (\prod \zeta_{i_0}, \prod \zeta_{i_0,i_1}) 
\in \prod_{i_0}\Gamma(\mathcal{LC}^2(U_{i_0}/T)) \oplus 
\prod_{i_0,i_1}\Gamma(\mathcal{LC}^1(U_{i_0i_1}/T)).$$ 
By Proposition 5, (1), $\zeta_{i_0}$ determines 
a Poisson deformation $\mathcal{U}_{i_0}$ of 
$U_{i_0}$ over $T[\epsilon]$. 
Moreover, $\zeta_{i_0i_1}$ determines a Poisson 
isomorphism $\mathcal{U}_{i_0}\vert_{U_{i_0i_1}} \cong 
\mathcal{U}_{i_1}\vert_{U_{i_0i_1}}$. 
One can construct a Poisson deformation of $\mathcal{X}$ 
of $X$ by patching together $\{\mathcal{U}_{i_0}\}$. 
Conversely, a Poisson deformation $\mathcal{X}$ is 
obtained by patching together local Poisson deformations 
$\mathcal{U}_i$ of $U_i$ for an affine open covering 
$\{U_i\}_{i \in I}$ of $X$. Each $\mathcal{U}_i$ determines 
$\zeta_i \in  \Gamma(\mathcal{LC}^2(U_i/T))$, and each Poisson 
isomorphism $\mathcal{U}_i\vert_{U_{ij}} \cong 
\mathcal{U}_j\vert_{U_{ij}}$ determines $\zeta_{ij} \in 
\Gamma(\mathcal{LC}^1(U_{ij}))$. Then 
$$ (\prod \zeta_i, \prod \zeta_{ij}) \in 
\prod_i \Gamma(\mathcal{LC}^2(U_i/T)) \oplus 
\prod_{i,j} \Gamma(\mathcal{LC}^1(U_{ij}/T))$$ 
is a 2-cocycle: hence gives an element of 
2-nd \v{C}ech cohomology.  
 
\section{Symplectic varieties} 
Assume that $X_0$ is a non-singular variety over 
$\mathbf{C}$ of dimension $2d$. Then $X_0$ is called a 
{\em symplectic manifold} if there is a 2-form $\omega_0 \in 
\Gamma(X_0, \Omega^2_{X_0})$ such that $d\omega_0 = 0$ 
and $\wedge^d\omega_0$ is a nowhere-vanishing section 
of $\Omega^{2d}_{X_0}$. The 2-form $\omega_0$ is called 
a {\em symplectic form}, and it gives an identification 
$\Omega^1_{X_0} \cong \Theta_{X_0}$. For a local section 
$f$ of $\mathcal{O}_{X_0}$, the 1-form $df$ corresponds to 
a local vector field $H_f$ by this identification. 
We say that $H_f$ is the {\em Hamiltonian vector field} 
for $f$. If we put $\{f,g\} := \omega(H_f, H_g)$, then 
$X_0$ becomes a Poisson scheme over $\mathrm{Spec}(\mathbf{C})$. 
Now let us consider a Poisson deformation $X$ of $X_0$ 
over $T := \mathrm{Spec}(S)$ with a local Artinian 
$\mathbf{C}$-algebra $S$ with $S/m_S = \mathbf{C}$. 
The Poisson bracket $\{\; , \;\}$ on $X$ can be written as 
$\{f, g\} = \Theta(df \wedge dg)$ for a relative bi-vector 
(Poisson bi-vector) 
$\Theta \in \Gamma(X, \wedge^2\Theta_{X/T})$. The restriction of 
$\Theta$ to the central fiber $X$ is nothing but the 
Poisson bi-vector for the original Poisson structure, which 
is non-degenerate because it is defined via the symplectic 
form $\omega_0$. Hence $\Theta$ is also a non-degenerate 
relative bi-vector. It gives an identification of $\Theta_{X/T}$ 
with $\Omega^1_{X/T}$. Hence $\Theta \in \Gamma(X, \wedge^2 \Theta_{X/T})$ 
defines an element $\omega \in \Gamma(X, \Omega^2_{X/T})$ that 
restricts to $\omega_0$ on $X_0$. One can define the 
Hamiltonian vector field $H_f \in \Theta_{X/T}$ for $f \in 
\mathcal{O}_X$. \vspace{0.15cm}

{\bf Proposition 9}. {\em Assume that $X$ is a Poisson deformation 
of a symplectic manifold $X_0$  
over an Artinian base $T$. Then $\mathcal{LC}^{\cdot}(X/T)$ is 
quasi-isomorphic to the truncated De Rham complex 
$(\Omega^{\geq 1}_{X/T}, d)$.} 
\vspace{0.12cm} 

{\em Proof}. By the symplectic form $\omega$, we have an 
identification $\phi: \Theta_{X/T} \cong \Omega^1_{X/T}$; hence, for 
each $i \geq 1$, we get $\wedge^i\Theta_{X/T} \cong \Omega^i_{X/T}$, 
which we denote also by $\phi$ (by abuse of notation). 
We shall prove that $\phi\circ\delta(f) = d\phi(f)$ for 
$f \in \wedge^i\Theta_{X/T}$. In order to do that, it suffices 
to check this for the $f$ of the form: $f = \alpha f_1 \wedge ... 
\wedge f_i$ with $\alpha \in \mathcal{O}_X$, $f_1$, ..., $f_i \in 
\Theta_{X/T}$. It is enough to check that 
$$d\phi(f)(H_{a_1}\wedge ... \wedge H_{a_{i+1}}) = 
\delta f(da_1 \wedge ... \wedge d{a_{i+1}}).$$ 

We shall calculate the left hand side. In the following, for 
simplicity, we will not write the $\pm$ signature exactly as 
$(-1)^{...}$, but only write $\pm$ because it does not cause 
any confusion. We have \vspace{0.3cm}

$(L.H.S.) = d(\alpha\omega(f_1, \cdot)\wedge ... 
\wedge\omega(f_i, \cdot))(H_{a_1}\wedge ... \wedge H_{a_{i+1}})$ 
\vspace{0.2cm}

$$ = \sum_{1 \le j \le i+1}(-1)^{j+1}(\sum_{\{l_1, ..., l_i\} = 
\{1, ..., \breve{j}, ..., i+1\}}\pm H_{a_j}(\alpha\omega(f_1, H_{a_{l_1}})
\cdot\cdot\cdot \omega(f_i,H_{a_{l_i}}))$$ 

$$ + \sum_{1 \le j < k \le i+1}(-1)^{j+k}
(\sum_{\{l_1, ..., \breve{l}, ..., l_i\} = \{1, ..., \breve{j}, ..., 
\breve{k}, ..., i+1\}}\pm\alpha\omega(f_1, H_{a_{l_1}})\times ... $$ 
 
$$... \times \omega(f_l, [H_{a_j}, H_{a_k}])\times ...
\times \omega(f_i, H_{a_{l_i}}))$$

$$ = \sum_{1 \le j \le i+1}(-1)^{j+1}(\sum\pm H_{a_j}(\alpha
f_1(da_{l_1})\cdot\cdot\cdot f_i(da_{l_i}))) $$
 
$$+ \sum_{1 \le j < k \le i+1}(-1)^{j+k}(\sum\pm \alpha f_1(da_{l_1})
\cdot\cdot\cdot 
f_l(d\{a_j, a_k\})\cdot\cdot\cdot f_i(da_{l_i}))$$

$$= \sum_{1 \le j \le i+1}(-1)^{j+1}H_{a_j}(\alpha
f(da_1 \wedge ... \wedge \breve{da_j} \wedge ... \wedge da_{i+1}))$$
  
$$+ \sum_{1 \le j < k \le i+1}(-1)^{j+k}\alpha
f(d\{a_j,a_k\} \wedge da_1 \wedge ... \wedge \breve{da_j} \wedge 
... \wedge \breve{da_k}\wedge ... \wedge da_{i+1})$$

$$ = \sum_{1 \le j \le i+1}(-1)^{j+1}\{a_j, \alpha f(da_1 \wedge ... 
\wedge \breve{da_j} \wedge ... \wedge da_{i+1})\}$$ 

$$+ \sum_{1 \le j < k \le i+1}(-1)^{j+k}\alpha f(d\{a_j, a_k\} \wedge 
... \wedge \breve{da_j} \wedge ... \wedge \breve{da_k} \wedge 
... \wedge da_{i+1}) $$

$= (R.H.S.)$. 
\vspace{0.3cm}

{\bf Corollary 10}. {\em Assume that $X$ is a Poisson deformation 
of a symplectic manifold $X_0$ over an Artinian base $T$.  
If $H^1(X, \mathcal{O}_X) 
= H^2(X, \mathcal{O}_X) = 0$, then $\mathrm{HP}^2(X/T) 
= H^2((X_0)^{an}, S)$, where $(X_0)^{an}$ is a complex analytic space 
associated with $X_0$ and $S$ is the constant sheaf with value in $S$.} 
\vspace{0.12cm}

{\em Proof}. By the distinguished triangle 
$$\Omega^{\geq 1}_{X/T} \to \Omega^{\cdot}_{X/T} \to 
\mathcal{O}_X \stackrel{[1]}\to \Omega^{\geq 1}_{X/T}[1]$$ 
we have an exact sequence 
$$\to \mathrm{HP}^i(X/T) \to \mathbf{H}^i(\Omega^{\cdot}_{X/T}) 
\to H^i(\mathcal{O}_X) \to .$$ 
Here $\mathbf{H}^i(X, \Omega^{\cdot}_{X/T}) \cong H^i((X_0)^{an}, 
S)$; from this we obtain 
the result. 
We prove this by an induction of $\mathrm{length}_{\mathbf{C}}(S)$. 
We take $t \in S$ such that $t\cdot m_S = 0$. For the exact sequence 
$$0 \to \mathbf{C} \stackrel{t}\to S \to \bar{S} \to 0,$$ 
define $\bar{X} := X \times_T \bar{T}$, where $\bar{T} := 
\mathrm{Spec}(\bar{S})$. Then we 
obtain a commutative diagrams of exact sequences: 

\begin{equation}
\begin{CD}
@>>> \mathbf{H}^i(X_0, \Omega^{\cdot}_{X_0}) @>>> 
\mathbf{H}^i(X, \Omega^{\cdot}_{X/T}) @>>> 
\mathbf{H}^i(\bar{X}, \Omega^{\cdot}_{\bar{X}/\bar{T}})  \\ 
@. @V{\cong}VV @VVV @V{\cong}VV  \\ 
@>>>  \mathbf{H}^i((X_0)^{an}, \Omega^{\cdot}_{(X_0)^{an}}) @>>> 
\mathbf{H}^i(X^{an}, \Omega^{\cdot}_{X^{an}/T}) @>>> 
\mathbf{H}^i(\bar{X}^{an}, \Omega^{\cdot}_{\bar{X}^{an}/\bar{T}})  
\end{CD}
\end{equation}

By a theorem of Grothendieck [G], the first vertical maps are isomorphisms 
and the third vertical maps are isomorphisms 
by the induction. Hence the middle vertical maps are also isomorphisms. 
By the Poincare lemma (cf. [De]). we know that 
$\mathbf{H}^i(X^{an}, \Omega^{\cdot}_{X^{an}/T}) \cong H^i((X_0)^{an}, S)$.  
\vspace{0.15cm}

{\bf Example 11}. When $f: X \to T$ is a proper smooth morphism 
of $\mathbf{C}$-schemes, by GAGA, we have 
$$\mathbf{R}^if_*\Omega^{\cdot}_{X/T} \otimes_{\mathcal{O}_T}
\mathcal{O}_{T^{an}} \cong R^i(f^{an})_*\mathbf{C}\otimes_{\mathbf{C}} 
\mathcal{O}_{T^{an}}$$ 
without the Artinian condition for $T$. But when $f$ is not proper, 
the structure of $\mathbf{R}^if_*\Omega^{\cdot}_{X/T}$ is 
complicated. For example, Put $X := \mathbf{C}^2\setminus\{xy = 1\}$, where 
$x$ and $y$ are standard coordinates of $\mathbf{C}^2$. Let 
$f: X \to T:=\mathbf{C}$ be the map defined by $(x,y) \to x$.      
Set $\hat{T} := \mathrm{Spec}\mathbf{C}[[x]]$ and 
$T_n := \mathrm{Spec}\mathbf{C}[x]/(x^{n+1})$. Define 
$\hat{X} := X \times_T \hat{T}$ and define $\hat{f}$ to be 
the natural map from $\hat{X} \to \hat{T}$. Finally put 
$X_n := X \times_T T_n$. 
Then 
\begin{enumerate}
\item $\mathbf{R}^1f_*\Omega^{\cdot}_{X/T}$ is a quasi-coherent 
sheaf on $T$, and $\mathbf{R}^1f_*\Omega^{\cdot}_{X/T}\vert_{T\setminus
\{0\}}$ is an invertible sheaf.  

\item $\mathrm{proj.lim}\;H^1(X_n, \Omega^{\cdot}_{X_n/T_n}) = 0.$
    
\end{enumerate}
  
{\bf Definition}. Let $X_0$ be a normal variety of 
dimenson $2d$ over $\mathbf{C}$ 
and let $U_0$ be its regular part. 
Then $X_0$ is a symplectic variety if $U_0$ admits a 
2-form $\omega_0$ such that 
\begin{enumerate} 
\item $d\omega_0 = 0$,  
\item $\wedge^d \omega_0$ is a nowhere-vanishing 
section of $\wedge^d\Omega^1_{U_0}$, 
\item for any resolution $\pi: Y_0 \to X_0$ of $X_0$ 
with $\pi^{-1}(U_0) \cong U_0$, $\omega_0$ extends to 
a (regular) 2-form on $Y_0$. 
\end{enumerate}  
\vspace{0.12cm} 

If $X_0$ is a symplectic variety, then $U_0$ becomes 
a Poisson scheme. Since $\mathcal{O}_{X_0} = 
(j_0)_*\mathcal{O}_{U_0}$, the Poisson bracket 
$\{\;, \;\}$ on $U_0$ uniquely extends to that on 
$X_0$. Thus $X_0$ is a Poisson scheme.   
By definition, its Poisson 
bi-vector $\Theta_0$ is non-degenerate over $U_0$. 
The $\Theta_0$ identifies $\Theta_{U_0}$ with 
$\Omega^{1}_{U_0}$; by this identification, $\Theta_0\vert_{U_0}$ 
corresponds to $\omega_0$.  A symplectic variety $X_0$ has rational 
Gorenstein singularities; in other words, $X$ has canonical 
singularities of index 1. When $X_0$ has only terminal singularities, 
$\mathrm{Codim}(\Sigma_0 \subset X_0) \geq 4$ for $\Sigma_0 
:= \mathrm{Sing}(X_0)$. \vspace{0.12cm}

{\bf Definition}. Let $X_0$ be a symplectic variety. 
Then $X_0$ is {\em convex} if there is a birational 
projective morphism from $X_0$ to an affine normal variety $Y_0$. 
In this case, $Y_0$ is isomorphic to $\mathrm{Spec}\Gamma(X_0, 
\mathcal{O}_{X_0})$. 
\vspace{0.12cm}

{\bf Lemma 12}. 
{\em Let $X_n$ be a Poisson deformation of a convex symplectic 
variety $X_0$ over $T_n := \mathrm{Spec}(S_n)$ with $S_n 
:= \mathbf{C}[t]/(t^{n+1})$. We define $U_n \subset X_n$ 
to be locus where $X_n \to S_n$ is smooth.   
Assume that $X_0$ has only terminal singularities. Then 
$\mathrm{HP}^2(U_n/T_n) \cong H^2((U_0)^{an}, S_n)$, where 
$S_n$ is the constant sheaf over $(U_0)^{an}$ with value in $S_n$.} 
\vspace{0.12cm}

{\em Proof}. Since $X_0$ has terminal singularities, $X_0$ is 
Cohen-Macaulay and $\mathrm{Codim}(\Sigma_0 \subset X_0) \geq 4$. 
Similarly, $X_n$ is Cohen-Macaulay and $\mathrm{Codim}(\Sigma_n 
\subset X_n) \geq 4$ for $\Sigma_n := \mathrm{Sing}(X_n \to T_n)$.  
The affine normal variety $Y_0$ has symplectic singularities; 
hence $Y_0$ has rational singularities. This implies that 
$H^i(X_0, \mathcal{O}_{X_0}) = 0$ for $i > 0$. 
Since $X_0$ is Cohen-Macaulay and $\mathrm{Codim}
(\Sigma_0 \subset X_0) \geq 4$, we see that 
$H^1(U_0, \mathcal{O}_{U_0}) = 
H^2(U_0, \mathcal{O}_{U_0}) = 0$ 
by the depth argument. By using the exact 
sequences $$0 \to \mathcal{O}_{U_0} \stackrel{t^k}\to 
\mathcal{O}_{U_k} \to \mathcal{O}_{U_{k-1}} \to 0$$ 
inductively, we conclude that $H^1(\mathcal{O}_{U_n}) 
= H^2(\mathcal{O}_{U_n}) = 0$. Then, by 
Corollary 10, we have $\mathrm{HP}^2(U_n/T_n) 
\cong H^2((U_0)^{an}, S_n)$.     
\vspace{0.15cm}

Let $X_n$ be the same as Lemma 12. 
Put $T_n[\epsilon] := \mathrm{Spec}(S_n[\epsilon])$ 
with $\epsilon^2 = 0$. As in Proposition 8, we 
define $\mathrm{PD}(X_n/T_n, T_n[\epsilon])$ to 
be the set of equivalence classes of the Poisson 
deformations of $X_n$ over $T_n[\epsilon]$. 
Let $\mathcal{X}_n$ be a Poisson deformation of 
$X_n$ over $T_n[\epsilon]$. Then we denote by 
$\mathrm{Aut}(\mathcal{X}_n, T_n)$ the set of all 
automorphisms of $\mathcal{X}_n$ as a Poisson 
$T_n[\epsilon]$-scheme such that they induce 
the identity map of $X_n$ over $T_n$. Then 
we have: 
\vspace{0.15cm}

{\bf Proposition 13}. 

(1) {\em There is a one-to-one 
correspondence between $\mathrm{HP}^2(U_n/T_n)$ 
and \\
$\mathrm{PD}(X_n/T_n, T_n[\epsilon])$.} 

(2) {\em There is a one-to-one correspondence 
between $\mathrm{HP}^1(U_n/T_n)$ and \\
$\mathrm{Aut}(\mathcal{X}_n. T_n)$.}  
\vspace{0.15cm} 

{\em Proof}. Assume that $\mathcal{U}_n$ is 
a Poisson deformation of $U_n$ over $T_n[\epsilon]$. 
Since $\mathrm{Codim}(\Sigma_n \subset X_n) \geq 3$ and 
$X_n$ is Cohen-Macaulay, by [K-M, 12.5.6],  
$$\mathrm{Ext}^1(\Omega^1_{X_n/T_n}, \mathcal{O}_{X_n}) 
\cong \mathrm{Ext}^1(\Omega^1_{U_n/T_n}, \mathcal{O}_{U_n}).$$ 
This implies that, over $T_n[\epsilon]$, 
$\mathcal{U}_n$ extends uniquely to an 
$\mathcal{X}_n$ so that it gives a flat deformation of $X_n$. 
Let us denote by $j: \mathcal{U}_n \to \mathcal{X}_n$ 
the inclusion map. Then, by the depth argument, 
we see that $\mathcal{O}_{\mathcal{X}_n} = 
j_*\mathcal{O}_{\mathcal{U}_n}$. Therefore, the 
Poisson structure on $\mathcal{U}_n$ also extends 
uniquely to that on $\mathcal{X}_n$. Now Proposition 
8 implies (1). As for (2), let $\mathcal{U}_n$ be the 
locus of $\mathcal{X}_n$ where $\mathcal{X}_n \to T_n[\epsilon]$ 
is smooth. Then, we see that $$\mathrm{Aut}(\mathcal{U}_n, T_n) 
= \mathrm{Aut}(\mathcal{X}_n, T_n),$$ which implies (2) 
again by Proposition 8. 
\vspace{0.15cm}  

Let $X$ be a convex symplectic variety with terminal 
singularities. We regard $X$ as a Poisson scheme by 
the natural Poisson structure $\{\; , \;\}$ induced by the symplectic form 
on the regular locus $U := (X)_{\mathrm{reg}}$. 
For a local Artinian $\mathbf{C}$-algebra 
$S$ with $S/m_S = \mathbf{C}$, we define 
$\mathrm{PD}(S)$ to be the set of equivalence classes 
of the pairs of Poisson deformations $\mathcal{X}$ of $X$ over 
$\mathrm{Spec}(S)$ and Poisson isomorphisms 
$\phi: \mathcal{X}\times_{\mathrm{Spec}(S)}\mathrm{Spec}(\mathbf{C}) 
\cong X$. Here $(\mathcal{X}, \phi)$ and 
$(\mathcal{X}', \phi')$ are equivalent if there is a Poisson 
isomorphism $\varphi: \mathcal{X} \cong \mathcal{X}'$ over 
$\mathrm{Spec}(S)$ which induces the identity map of $X$ over 
$\mathrm{Spec}(\mathbf{C})$ via $\phi$ and $\phi'$. 
We define the {\em Poisson deformation functor}:  
$$\mathrm{PD}_{(X, \{\;, \;\})}: (\mathrm{Art})_{\mathbf{C}} 
\to (\mathrm{Set})$$ by $\mathrm{PD}(S)$ for 
$S \in (\mathrm{Art})_{\mathbf{C}}$.      
\vspace{0.15cm}

{\bf Theorem 14}. {\em Let $(X, \{\; , \;\})$ be a Poisson 
scheme associated with a convex symplectic varieties 
with terminal singularities. Then $\mathrm{PD}_{(X, \{\;, \;\})}$ 
has a pro-reprentable hull in the sense of Schlessinger. 
Moreover $\mathrm{PD}$ is pro-representable.} 
\vspace{0.12cm}

{\em Proof}. We have to check Schlessinger's conditions [Sch] 
for the existence of a hull. By Proposition 13, $\mathrm{PD}
(\mathbf{C}[\epsilon]) = H^2(U^{an}, \mathbf{C}) < \infty$. 
Other conditions are checked in a similar way as the case 
of usual deformations. For the last statement, we have 
to prove the following. Let $\mathcal{X}$ be a Poisson 
deformation of $X$ over an Artinian base $T$, 
and let $\bar{\mathcal{X}}$ be its restriction over 
a closed subscheme $\bar{T}$ of $T$. Then, any Poisson 
automorphism of $\bar{\mathcal{X}}$ over $\bar{T}$ inducing 
the identity map on $X$,  
extends to a Poisson automorphism of $\mathcal{X}$ 
over $T$. Let $R$ be the pro-representable hull 
of $\mathrm{PD}$ and put $R_n := R/(m_R)^{n+1}$. Take a formal versal 
Poisson deformation $\{\mathcal{X}_n\}$ over $\{R_n\}$. Note that, if we 
are given an Artinian local $R$-algebra $S$ with residue field $\mathbf{C}$, 
then we get a Poisson deformation $X_S$ of $X$ over $\mathrm{Spec(S)}$. 
We then define $\mathrm{Aut}(S)$ to be the set of all 
Poisson automorphisms of $X_S$ over $\mathrm{Spec}(S)$ which 
induce the identity map of $X$. Let 
$$\mathrm{Aut} : (\mathrm{Art})_R \to (\mathrm{Set})$$ 
be the covariant functor defined in this manner. We want to 
prove that $\mathrm{Aut}(S) \to \mathrm{Aut}(\bar{S})$ is 
surjective for any surjection $S \to \bar{S}$. It is enough to 
check this only for a small extension $S \to \bar{S}$, that is, 
the kernel $I$ of $S \to \bar{S}$ is generated by an 
element $a$ such that $am_S = 0$. For each small extension 
$S \to \bar{S}$, one can define the {\em obstruction map} 
$$ \mathrm{ob}: \mathrm{Aut}(\bar{S}) \to a\cdot\mathrm{HP}^2(U)$$ 
in such a way that any element $\phi \in \mathrm{Aut}(\bar{S})$ 
can be lifted to an element of $\mathrm{Aut}(S)$ if and only if 
$\mathrm{ob}(\phi) = 0$. The obstruction map is constructed as 
follows. For $\phi \in \mathrm{Aut}(\bar{S})$, we have two 
Poisson extensions $X_{\bar{S}} \to X_S$ and $X_{\bar{S}} 
\stackrel{\phi}\to  X_{\bar{S}} \to X_S$. This gives an 
element of $a\cdot\mathrm{HP}^2(U)$ (cf. Proposition 13 
\footnote{Exactly, one can prove the following. Let $T := 
\mathrm{Spec}(S)$ with a local Artinian $\mathbf{C}$-algebra $S$ 
with $S/m_s = \mathbf{C}$. 
Let $X \to T$ be a Poisson deformation of a 
convex symplectic variety $X_0$ with only  
terminal singularities. 
Assume that $T$ is a closed subscheme of $T'$ defined by 
the ideal sheaf $I = (a)$ such that $a\cdot m_{S'} = 0$. 
Denote by $\mathrm{PD}(X/T, T')$ the set of equivalence 
classes of Poisson deformations of $X$ over $T'$. 
If $\mathrm{PD}(X/T, T') \ne \emptyset$, then 
$\mathrm{HP}^2(U_0) \cong \mathrm{PD}(X/T, T').$}). 
Obviously, if this element 
is zero, then these two extensions are equivalent and 
$\phi$ extends to a Poisson automorphism of $X_S$. 
\vspace{0.12cm}

{\bf Case 1} ($S = S_{n+1}$ and $\bar{S} := S_n$): 
We put $S_n := \mathrm{C}[t]/(t^{n+1})$. We shall prove that 
$\mathrm{Aut}(S_{n+1}) \to \mathrm{Aut}(S_n)$ is surjective.  
Taking Proposition 13, (2) into consideration, 
we say that $X$ has $T^0$-lifting property 
if, for any Poisson deformation $X_n$ of $X$ over $T_n := \mathrm{Spec}(S_n)$ 
and its restriction $X_{n-1}$ over $T_{n-1} := \mathrm{Spec}(S_{n-1})$, 
the natural map $\mathrm{HP}^1(U_n/T_n) \to \mathrm{HP}^1(U_{n-1}/T_{n-1})$ 
is surjective. \vspace{0.12cm} 

{\bf Claim}. {\em $X$ has $T^0$-lifting property.} 
\vspace{0.12cm}

{\em Proof}. Note that $X_n$ is Cohen-Macaulay. 
Let $U_n$ be the locus of $X_n$ where 
$X_n \to T_n$ is smooth. We put   
$$K_n := \mathrm{Coker}[H^0(U^{an}, 
S_n) \to H^0(U_n, \mathcal{O}_{U_n})]. $$ 
By the proof of Corollary 10, there is an 
exact sequence 
$$0 \to K_n \to \mathrm{HP}^1(U_n/T_n) \to 
H^1(U^{an}, S_n) \to 0.$$ 
Since 
$H^1(U, \mathcal{O}_U) = 0$, the restriction map 
$H^0(U_n, \mathcal{O}_{U_n}) \to 
H^0(U_{n-1}, \mathcal{O}_{U_{n-1}})$ is surjective. 
Hence the map $K_n \to K_{n-1}$ is surjective. 
On the other hand, $H^1(U^{an}, S_n) \to 
H^1(U^{an}, S_{n-1})$ is also surjective; hence 
the result follows. 
\vspace{0.12cm}

Note that $t \to t + \epsilon$ induces the commutative 
diagram of exact sequences: 

\begin{equation} 
\begin{CD} 
0 @>>> (t^{n+1}) @>>> S_{n+1} @>>> S_n @>>> 0 \\ 
@. @V{\cong}VV @VVV @VVV @. \\ 
0 @>>> (t^n \epsilon ) @>>> S_n[\epsilon] @>>> 
S_{n-1}[\epsilon]\times_{S_{n-1}}S_n @>>> 0 
\end{CD} 
\end{equation} 

Applying $\mathrm{Aut}$ to this diagram, we obtain 

\begin{equation}
\begin{CD} 
\mathrm{Aut}(S_{n+1}) @>>> \mathrm{Aut}(S_n) @>{\mathrm{ob}}>> 
t^{n+1}\cdot\mathrm{HP}^2(U) \\ 
@VVV @VVV @V{\cong}VV \\ 
\mathrm{Aut}(S_n[\epsilon]) @>>> 
\mathrm{Aut}(S_{n-1}[\epsilon]\times_{S_{n-1}}S_n) @>{\mathrm{ob}}>> 
t^n\epsilon \cdot \mathrm{HP}^2(U) 
\end{CD}
\end{equation}  

The $T^0$-lifting property implies that 
the map $\mathrm{Aut}(S_n[\epsilon]) \to 
\mathrm{Aut}(S_{n-1}[\epsilon]\times_{S_{n-1}}S_n)$ 
is surjective. Hence, by the commutative diagram, we 
see that $\mathrm{Aut}(S_{n+1}) \to 
\mathrm{Aut}(S_n)$ is surjective. 
\vspace{0.12cm}

{\bf Case 2} (general case): For any small extension 
$S \to \bar{S}$, one can find the following commutative 
diagram for some $n$: 

\begin{equation}
\begin{CD} 
0 @>>> aS @>>> S @>>> \bar{S} @>>> 0 \\ 
@. @V{\cong}VV @VVV @VVV @. \\
0 @>>> (t^{n+1}) @>>> S_{n+1} @>>> S_n @>>> 0 \\ 
\end{CD} 
\end{equation} 

Applying $\mathrm{Aut}$ to this diagram, we get: 

\begin{equation} 
\begin{CD} 
\mathrm{Aut}(S) @>>> \mathrm{Aut}(\bar{S}) @>{\mathrm{ob}}>> 
a\cdot \mathrm{HP}^2(U) \\ 
@VVV @VVV @V{\cong}VV \\ 
\mathrm{Aut}(S_{n+1}) @>>> \mathrm{Aut}(S_n) @>{\mathrm{ob}}>> 
t^{n+1}\cdot \mathrm{HP}^2(U) 
\end{CD} 
\end{equation}  

By Case 1, we already know that $\mathrm{Aut}(S_{n+1}) 
\to \mathrm{Aut}(S_n)$ is surjective. By the commutative 
diagram we see that $\mathrm{Aut}(S) \to 
\mathrm{Aut}(\bar{S})$ is surjective. 
\vspace{0.15cm}

{\bf Corollary 15}. {\em Let $(X, \{\;, \;\})$ be the same 
as Theorem 14. Then} \vspace{0.12cm}

(1) {\em $X$ has $T^1$-lifting property.} (cf. [Kaw, Na 5]) 
\vspace{0.12cm}

(2) {\em $\mathrm{PD}_{(X, \{\;, \;\})}$ is unobstructed.} 
\vspace{0.12cm}

{\em Proof}. (1): We put $S_n := \mathbf{C}[t]/(t^{n+1})$ and 
$T_n := \mathrm{Spec}(S_n)$.  
Let $X_n$ be a Poisson deformation of $X$ over $T_n$ and 
let $X_{n-1}$ be its restriction over $T_{n-1}$. By Proposition 13,(1), 
we have to prove that $\mathrm{HP}^2(U_n/T_n) \to 
\mathrm{HP}^2(U_{n-1}/T_{n-1})$ is surjective. By Lemma  
12, $\mathrm{HP}^2(U_n/T_n) \cong H^2(U^{an}, S_n)$. 
Since $H^2(U^{an}, S_n) = H^2(U^{an}, \mathbf{C})\otimes_{\mathbf{C}}
S_n$, we conclude that this map is surjective. 
\vspace{0.12cm}

(2): By Theorem 14, $\mathrm{PD}$ has a pro-representable 
hull $R$. Denote by $h_R: (\mathrm{Art})_{\mathbf{C}} 
\to (\mathrm{Set})$ the covariant functor defined by 
$h_R(S) := \mathrm{Hom}_{\mathrm{local \; {\bf C}-alg.}}
(R, S)$. Since $\mathrm{PD}$ is pro-representable by 
Theorem 14, $h_R = \mathrm{PD}$. We write $R$ as 
$\mathbf{C}[[x_1, ..., x_r]]/J$ with $r := \dim_{\mathbf{C}}
m_R/(m_R)^2$. Let $S$ and $S_0$ be the objects of 
$(\mathrm{Art})_{\mathbf{C}}$ such that $S_0 = S/I$ with an 
ideal $I$ such that $Im_S = 0$. Then we have an exact 
sequence (cf. [Gr, (1.7)]) 
$$ h_R(S) \to h_R(S_0) \stackrel{ob}\to 
(J/m_RJ)^{*}\otimes_{\mathbf{C}}I.$$ 

By sending $t$ to $t + \epsilon$, we have the commutative 
diagram of exact sequences: 

\begin{equation} 
\begin{CD} 
0 @>>> (t^{n+1}) @>>> S_{n+1} @>>> S_n @>>> 0 \\ 
@. @V{\cong}VV @VVV @VVV @. \\ 
0 @>>> (t^n \epsilon ) @>>> S_n[\epsilon] @>>> 
S_{n-1}[\epsilon]\times_{S_{n-1}}S_n @>>> 0 
\end{CD} 
\end{equation} 

Applying $h_R$ to this diagram, we obtain 

\begin{equation}
\begin{CD} 
h_R(S_{n+1}) @>>> h_R(S_n) @>{\mathrm{ob}}>> 
t^{n+1}\otimes(J/m_RJ)^{*} \\ 
@VVV @VVV @V{\cong}VV \\ 
h_R(S_n[\epsilon]) @>>> 
h_R(S_{n-1}[\epsilon]\times_{S_{n-1}}S_n) @>{\mathrm{ob}}>> 
t^n\epsilon \otimes (J/m_RJ)^{*} 
\end{CD}
\end{equation}  
  
By (1), we see that $h_R(S_n[\epsilon]) \to  
h_R(S_{n-1}[\epsilon]\times_{S_{n-1}}S_n)$ 
is surjective. Then, by the commutative diagram, 
we conclude that 
$h_R(S_{n+1}) \to h_R(S_n)$ is surjective. 
\vspace{0.15cm}

{\bf Twistor deformations} (cf. [Ka 1]): Let $X$ be a 
convex symplectic variety with terminal singularities. 
We put $U := X_{reg}$. Let $\{\;, \;\}$ be the natural 
Poisson structure on $X$ defined by the symplectic 
form $\omega$ on $U$. Fix a line bundle $L$ on $X^{an}$. 
Define a class $[L]$ of $L$ as the image of $L$ by the map 
$$H^1(U^{an}, \mathcal{O}^*_{U^{an}}) \to 
\mathbf{H}^2(U^{an}, \Omega^{\cdot}_{U^{an}}) 
\cong H^2(U^{an}, \mathbf{C}).$$ 
We put $S_n := \mathbf{C}[t]/(t^{n+1})$ and 
$T_n := \mathrm{Spec}(S_n)$. By Proposition 13, (1), the element 
$[L] \in H^2(U^{an}, \mathbf{C})$ determines 
a Poisson deformation $X_1$ of $X$ 
over $T_1$. We shall construct Poisson deformations $X_n$ 
over $T_n$ inductively. Assume that we already have a 
Poisson deformation $X_n$ over $T_n$. Define $X_{n-1}$ 
to be the restriction of $X_n$ over $T_{n-1}$. Since 
$H^1(X^{an}, \mathcal{O}_{X^{an}}) = 
H^2(X^{an}, \mathcal{O}_{X^{an}}) = 0$, 
$L$ extends uniquely to a line bundle $L_n$ on $(X_n)^{an}$. 
Denote by $L_{n-1}$ the restriction of $L_n$ to $(X_{n-1})^{an}$.  
Consider the map $S_n \to S_{n-1}[\epsilon]$ defined 
by $t \to t + \epsilon$. This map induces 
$$\mathrm{PD}(S_n) \to \mathrm{PD}(S_{n-1}[\epsilon]).$$
The class $[L_{n-1}] \in H^2(U^{an}, S_{n-1})$ determines  
a Poisson deformation $(X_{n-1})'$ of $X_{n-1}$ over 
$T_{n-1}[\epsilon]$. Assume that $X_n$ satisfies the 
condition \vspace{0.12cm}

$(*)_n$ : $[X_n] \in \mathrm{PD}(S_n)$ is sent to 
$[(X_{n-1})'] \in \mathrm{PD}(S_{n-1}[\epsilon])$.       
\vspace{0.12cm}

Note that $X_1$ actually has this property. We shall 
construct $X_{n+1}$ in such a way that $X_{n+1}$ satisfies 
$(*)_{n+1}$. 
Look at the commutative diagram: 

\begin{equation}
\begin{CD}
\mathrm{PD}(S_{n+1}) @>>> \mathrm{PD}(S_n) \\ 
@VVV @VVV \\ 
\mathrm{PD}(S_n[\epsilon]) @>>> 
\mathrm{PD}(S_{n-1}[\epsilon]\times_{S_{n-1}}S_n)
\end{CD}
\end{equation}  

Note that we have an element 
$$[X_n \leftarrow X_{n-1} \rightarrow (X_{n-1})')] 
\in \mathrm{PD}(S_{n-1}[\epsilon]\times_{S_{n-1}}S_n).$$  
Identifying $\mathrm{HP}^2(U_n/T_n)$ with $H^2(U^{an}, 
S_n)$, $[L_n]$ is sent to $[L_{n-1}]$ by the map 
$$\mathrm{HP}^2(U_n/T_n) \to \mathrm{HP}^2(U_{n-1}/T_{n-1}).$$ 
Now, by Proposition 13,(1), we get a lifting $[(X_n)'] \in 
\mathrm{PD}(S_n[\epsilon])$ of  
$$[X_n \leftarrow X_{n-1} \rightarrow (X_{n-1})')] 
\in \mathrm{PD}(S_{n-1}[\epsilon]\times_{S_{n-1}}S_n)$$ 
corresponding to $[L_n]$. 
By the standard argument used in $T^1$-lifting 
principle (cf. proof of Corollary 15, (2)), one 
can find a Poisson deformation $X_{n+1}$ such that 
$[X_{n+1}] \in \mathrm{PD}(S_{n+1})$ is sent to 
$[(X_n)'] \in \mathrm{PD}(S_n[\epsilon])$ in the 
diagram above. Moreover, since $\mathrm{PD}$ is 
pro-representable, such $[X_{n+1}]$ is unique. 
By the construction, $X_{n+1}$ satisfies $(*)_{n+1}$. 
This construction do not need the sequence of line 
bundles $L_n$ on $(X_n)^{an}$; we only need the sequence 
of line bundles on $(U_n)^{an}$. For example, if we are 
given a line bundle $L^0$ on $U^{an}$. Then, since 
$H^i(U^{an}, \mathcal{O}_{U^{an}}) = 0$ for $i= 1,2$, we have 
a unique extension $L^0_n \in \mathrm{Pic}((U_n)^{an})$. 
By using this, one can construct a formal deformation 
of $X$. \vspace{0.15cm}

{\bf Definition}.    
(1) When $L \in \mathrm{Pic}(X^{an})$, we call the 
formal deformation $\{X_n\}_{n\geq 1}$ 
the {\em twistor deformation} of $X$ associated with 
$L$. \vspace{0.15cm} 

(2) More generally, for $L^0 \in \mathrm{Pic}(U^{an})$, 
we call, the formal deformation $\{X_n\}_{n \geq 1}$ similarly 
constructed, the {\em quasi-twistor deformation} of 
$X$ associated with $L^0$. When $L^0$ extends to a line 
bundle $L$ on $X^{an}$, the corresponding quasi-twistor deformation  
coincides with the twistor deformation associated with $L$. 
\vspace{0.15cm}    

We next define the {\em Kodaira-Spencer class} of 
the formal deformation $\{X_n\}$. As before, we denote by 
$U_n$ the locus of $X_n$ where $f_n: X_n \to T_n$ is 
smooth. We put $f^0_n := f_n\vert_{U_n}$. 
The extension class 
$\theta_n \in H^1(U, \Theta_{U_{n-1}/T_{n-1}})$
of the exact sequence 
$$0 \to (f^0_n)^*\Omega^1_{T_n/\mathbf{C}} \to \Omega^1_{U_n/\mathbf{C}} 
\to \Omega^1_{U_n/T_n} \to 0$$ 
is the Kodaira Spencer class for $f_n: X_n \to T_n$. 
Here note that $\Omega^1_{T_n} \cong \mathcal{O}_{T_{n-1}}dt$. 
\vspace{0.15cm}

{\bf Lemma 16}. {\em Let $\{X_n\}$ be the twistor deformation 
of $X$ associated with $L \in \mathrm{Pic}(X^{an})$. 
Write $L_n \in \mathrm{Pic}(X^{an}_n)$ for the extension of $L$ to 
$X_n$. Let $\omega_n \in \Gamma(U_n, \Omega^2_{U_n/T_n})$ be 
the symplectic form defined by the Poisson $T_n$-scheme $X_n$. 
Then} 
$$ \imath(\theta_{n+1})(\omega_n) = [L_n] \in 
H^1(U, \Omega^1_{U_n/T_n}),$$ 
{\em where the left hand side is the interior product.} 
\vspace{0.15cm}

{\em Proof}. We use the same notation in the definition of 
a twistor deformation. By the commutative diagram 

\begin{equation}
\begin{CD}
(X_n)' @>>> X_{n+1} \\
@VVV @VVV \\ 
T_n[\epsilon] @>>> T_{n+1}
\end{CD} 
\end{equation}

we get the commutative diagram of exact sequences: 

\begin{equation}
\begin{CD}
0 @>>> \mathcal{O}_{U_n}d\epsilon @>>> 
\Omega^1_{(U_n)'/T_n}\vert_{U_n} @>>> \Omega^1_{U_n/T_n} 
@>>> 0 \\ 
@. @A{\cong}AA @A{\cong}AA @A{\cong}AA @. \\ 
0 @>>> \mathcal{O}_{U_n}dt @>>> 
\Omega^1_{U_{n+1}/\mathbf{C}}\vert_{U_n} @>>> 
\Omega^1_{U_n/T_n} @>>> 0 
\end{CD}
\end{equation}

The second exact sequence is the Kodaira-Spencer's 
sequence where the first term is 
$(f^0)^*\Omega^1_{T_{n+1}/\mathbf{C}}$ and the 
third term is $\Omega^1_{U_{n+1}/T_{n+1}}\vert_{U_n}$. 
Let $\eta \in H^1(U, \Theta_{U_n/T_n})$ be the 
extension class of the first exact sequence. By the 
definition of $(X_n)'$, we have $i(\eta)(\omega_n) 
= [L_n]$. On the other hand, the extension class 
of the second exact sequence is $\theta_{n+1}$. 
Hence $\eta = \theta_{n+1}$. 
\vspace{0.15cm}

Let $\{X_n\}$ be the twistor deformation of $X$ 
asociated with $L \in \mathrm{Pic}(X)$. For each $n$, 
we put $Y_n := \mathrm{Spec}\Gamma(X_n , \mathcal{O}_{X_n})$. 
$Y_n$ is an affine scheme over $T_n$. 
Since $H^1(X, \mathcal{O}_X) = 0$, $\Gamma(X_n, 
\mathcal{O}_{X_n}) \to \Gamma(X_{n-1}, \mathcal{O}_{X_{n-1}})$ 
is surjective. Define $$Y_{\infty} := \mathrm{Spec}
(\lim_{\leftarrow}\Gamma(X_n, \mathcal{O}_{X_n})).$$ Note that 
$Y_{\infty}$ is an affine variety over $T_{\infty} := 
\mathrm{Spec}\mathbf{C}[[t]]$. Fix an ample line bundle 
$A$ on $X$. Since $H^1(X, \mathcal{O}_X) 
= H^2(X, \mathcal{O}_X) = 0$, $A$ extends uniquely 
to ample line bundles $A_n$ on $X_n$. Then, by 
[EGA III, Th\'{e}or\`{e}me (5.4.5)], there is 
an {\em algebraization} $X_{\infty}$ of $\{X_n\}$ such that 
$X_{\infty}$ is a projective scheme over $Y_{\infty}$ and 
$X_{\infty}\times_{Y_{\infty}}Y_n \cong X_n$ for all $n$. 
By [ibid, Theoreme 5.4.1], the algebraization $X_{\infty}$ is 
unique. 
We denote by $g_{\infty}$ the projective morphism 
$X_{\infty} \to Y_{\infty}$.  
\vspace{0.2cm}

{\bf Theorem 17}. {\em Let $X$ be 
a convex symplectic variety with terminal 
singularities. Let $(X, \{\;, \;\})$ 
be the Poisson structure induced by 
the symplectic form on the regular part. 
Assume that $X^{an}$ is 
{\bf Q}-factorial \footnote{ 
Since $X$ is convex, there is a projective birational 
morphism $f$ from $X$ to an affine variety $Y$. 
 Take a reflexive sheaf $F$ on $X^{an}$ of rank 1. 
The direct image $f^{an}_*F^*$ of the dual sheaf $F^*$ 
is a coherent sheaf on the Stein variety $Y^{an}$. Hence 
$f^{an}_*F^*$ has a non-zero global section; in other words, 
there is an injection $\mathcal{O}_{X^{an}} \to 
F^*$. By taking its dual, $F$ is embedded in $\mathcal{O}_{X^{an}}$.   
Thus, $F = \mathcal{O}(-D)$ for an analytic effective 
divisor $D$. So, for any reflexive sheaf $F$ of rank $1$, the 
double dual sheaf $(F^{\otimes m})^{**}$ becomes an invertible 
sheaf for some $m$.}. Then any Poisson 
deformation of $(X, \{\; , \;\})$ is 
locally trivial as a flat deformation 
(after forgetting Poisson structure).} 
\vspace{0.12cm}

{\em Proof}. Define a subfunctor $$\mathrm{PD}_{lt}: 
(\mathrm{Art})_{\mathbf{C}} \to (\mathrm{Set})$$  
of $\mathrm{PD}$ by setting  
$\mathrm{PD}_{lt}(S)$ to be the set of 
equivalence classes of Poisson deformations of 
$(X, \{\;, \;\})$ over $\mathrm{Spec}(S)$ which 
are {\em locally trivial} as usual flat deformations. 
One can check that $\mathrm{PD}_{lt}$ 
has a pro-representable hull. Let $X_n \to T_n$ 
be an object of $\mathrm{PD}_{lt}(S_n)$, where 
$S_n := \mathrm{C}[t]/(t^{n+1})$ and $T_n := 
\mathrm{Spec}(S_n)$. Write $T^1_{X_n/T_n}$ 
for $\underline{\mathrm{Hom}}(\Omega^1_{X_n/T_n}, 
\mathcal{O}_{X_n})$. By Proposition 13, we have 
a natural map $$\mathrm{HP}^2(U_n/T_n) 
\to \mathrm{Ext}^1(\Omega^1_{X_n/T_n}, \mathcal{O}_{X_n}).$$ 
Define $T(X_n/T_n)$ to be 
the kernel of the composite 
$$\mathrm{HP}^2(U_n/T_n) \to 
\mathrm{Ext}^1(\Omega^1_{X_n/T_n}, \mathcal{O}_{X_n}) 
\to H^0(X_n, T^1_{X_n/T_n}).$$ 
Let $\mathrm{PD}_{lt}(X_n/T_n; T_n[\epsilon])$ be 
the set of equivalence classes of Poisson deformations 
of $X_n$ over $T_n[\epsilon]$ which are locally trivial 
as usual deformations. Here $T_n[\epsilon] := 
\mathrm{Spec}(S_n[\epsilon])$ and $S_n[\epsilon] 
= \mathbf{C}[t, \epsilon]/(t^{n+1}, \epsilon^2)$. 
Two Poisson deformations of $X_n$ over $T_n[\epsilon]$ 
are equivalent  
if there is a Poisson $T_n[\epsilon]$-somorphisms 
between them which induces the identity of $X_n$. 
Then there is a one-to-one correspondence 
between $T(X_n/T_n)$ and $\mathrm{PD}_{lt}(X_n/T_n; 
T_n[\epsilon])$.  
\vspace{0.15cm}

{\bf Lemma 18}. $T(X_n/T_n) = \mathrm{HP}^2(U_n/T_n).$ 
\vspace{0.12cm}

{\em Proof}. Since $H^0(X, T^1_{X_n/T_n}) \subset 
H^0(X^{an}, T^1_{X^{an}_n/T_n})$, it suffices to 
prove that $\mathrm{HP}^2(U_n/T_n) \to 
H^0(X^{an}, T^1_{X^{an}_n/T_n})$ is the zero map. 
In order to do this, for $p \in \Sigma (=\mathrm{Sing}(X))$, 
take a Stein open neighborhood $X^{an}_n(p)$ of $p \in X_n$, and 
put $U^{an}_n(p) := X^{an}(p) \cap U^{an}_n$. We have to 
prove that $H^2(U^{an}, S_n) \to H^2(U^{an}_n(p), S_n)$ 
is the zero map. In fact, on one hand, $$\mathrm{HP}^2(U_n/T_n) 
\cong H^2(U^{an}, S_n)$$ by Lemma 12. On the other hand, 
$H^0(X^{an}_n(p), T^1_{X^{an}_n/T_n}) \cong H^1(U^{an}_n(p), 
\Theta_{U^{an}_n(p)})$ (cf. the proof of [Na, Lemma 1]). 
By the symplectic form $\omega_n \in \Gamma(U_n, \Omega^2_{X_n/T_n})$, 
$\Theta_{U^{an}_n(p)}$ is identified with $\Omega^1_{U^{an}_n(p)}$. 
Hence $$H^0(X^{an}_n(p), T^1_{X^{an}_n/T_n}) \cong 
H^1(U^{an}_n(p), \Omega^1_{U^{an}_n(p)}).$$ 
By these identifications, the map 
$\mathrm{HP}^2(U_n/T_n) \to H^0(X^{an}_n(p), T^1_{X^{an}_n/T_n})$ 
coincides with the composite 
$$H^2(U^{an}, S_n) \to H^2(U^{an}(p), S_n) \to H^1(U^{an}(p), 
\Omega^1_{U^{an}_n(p)}),$$  
where the second map is induced by the spectral sequence 
$$E^{p,q}_1 := H^q(U^{an}_n(p), \Omega^p_{U^{an}_n(p)}) 
=> H^{p+q}(U^{an}, S_n)$$ 
(for details, see the proof of [Na, Lemma 1]).     
Let us consider the commutative diagram 
\begin{equation}
\begin{CD} 
\mathrm{Pic}(X^{an})\otimes_{\mathbf{Z}}S_n @>>> 
\mathrm{Pic}(U^{an})\otimes_{\mathbf{Z}}S_n @>{\cong}>> 
H^2(U^{an}, S_n) \\ 
@VVV @VVV @VVV \\ 
\mathrm{Pic}(X^{an}(p))\otimes_{\mathbf{Z}}S_n @>>> 
\mathrm{Pic}(U^{an}(p))\otimes_{\mathbf{Z}}S_n @>{\cong}>> 
H^2(U^{an}(p), S_n) 
\end{CD} 
\end{equation} 

Here the second map on the first row is an isomorphism 
because \\$H^1(U^{an}, \mathcal{O}_{U^{an}}) = 
H^2(U^{an}, \mathcal{O}_{U^{an}}) = 0$. 
Since $\mathrm{Codim}(\Sigma \subset X) \geq 3$, 
any line bundle on $U^{an}$ extends to a coherent 
sheaf on $X^{an}$. Thus, by the {\bf Q}-factoriality 
of $X^{an}$, the first map on the first row is 
surjective. If we take $X^{an}(p)$ small enough, 
then $\mathrm{Pic}(X^{an}(p)) = 0$. Now, by the 
commutative diagram above, we conclude that $H^2(U^{an}, S_n) 
\to H^2(U^{an}(p), S_n)$ is the zero map. This 
completes the proof of Lemma 18.  
\vspace{0.15cm}

Let us return to the proof of Theorem 17. 
The functor $\mathrm{PD}$ has 
$T^1$-lifting property by Corollary 15. 
By Lemma 18, $\mathrm{PD}_{lt}$ also has 
$T^1$-lifting property. Let $R$ and $R_{lt}$ 
be the pro-representable hulls of $\mathrm{PD}$ 
and $\mathrm{PD}_{lt}$ respectively. 
Then these are both regular local $\mathbf{C}$-algebra. 
There is a surjection $R \to R_{lt}$ because 
$\mathrm{PD}_{lt}$ is a sub-functor of 
$\mathrm{PD}$. By Lemma 18, the cotangent spaces 
of $R$ and $R_{lt}$ coincides. Hence 
$R \cong R_{lt}$.      
\vspace{0.15cm}

{\bf Theorem 19}. {\em Let $X$ be a convex symplectic 
variety with terminal singularities. Let $L$ 
be a (not necessarily ample) line bundle on $X^{an}$. 
Then the twistor deformation $\{X_n\}$ of $X$ associated 
with $L$ is locally trivial 
as a flat deformation.} 
\vspace{0.15cm} 
 
{\em Proof}.  Define $U_n \subset X_n$ to 
be the locus where $X_n \to T_n$ is smooth. We put 
$\Sigma := \mathrm{Sing}(X)$. For each point $p \in \Sigma$, 
we take a Stein open neighborhood $p \in X_n(p)$ in $(X_n)^{an}$, and 
put ${U}^{an}_n(p) := X_n(p) \cap U^{an}_n$. 
Let $L_n \in \mathrm{Pic}(X_n)$ be the (unique) extension 
of $L$ to $X_n$. 
We shall show that $[L_n] \in H^2(U^{an}, S_n)$ is sent to 
zero by the map 
$$H^2(U^{an}, S_n) \to H^2(U^{an}(p), S_n).$$ This is 
enough for us to prove that the twistor deformation 
$\{X_n\}$ is locally trivial. In fact, we have to show that 
the {\em local} Kodaira-Spencer class $\theta^{loc}_{n+1}(p) 
\in H^1(U^{an}_n(p), \Theta_{U^{an}_n(p)})$ is zero. 
By the same argument as Lemma 16, one can show that 
$$\iota (\theta^{loc}_{n+1}(p))(\omega_n) = 
[L_n\vert_{U^{an}(p)}] 
\in H^1(U^{an}_n(p), \Omega^1_{U^{an}_n(p)}).$$ 
Now let us consider 
the commutative diagram induced from the Hodge spectral sequences: 
  
\begin{equation}
\begin{CD} 
H^2(U^{an}, S_n)  @>>> H^1(U_n, \Omega^1_{U_n/T_n}) \\ 
@VVV @VVV \\ 
H^2(U^{an}(p), S_n) @>>> 
H^1(U^{an}_n(p), \Omega^1_{U^{an}_n/T_n}) 
\end{CD} 
\end{equation} 
 
For the existence of the first horizontal map, 
we use Grothendieck's theorem [G] and the fact 
$H^i(U_n, \mathcal{O}_{U_n}) = 0$ ($i = 1,2$) (cf. Lemma 12). 
Since $X^{an}_n$ is Cohen-Macaulay and $\mathrm{Codim}(\Sigma 
\subset X) \geq 4$, we have  
$H^i(U^{an}_n(p), \mathcal{O}_{U^{an}_n(p)}) = 0$ for 
$i = 1,2$, by the depth 
argument. This assures the existence of the second horizontal map. 
The vertical map on the right-hand side is just the 
composite of the maps 
$$ H^1(U_n, \Omega^1_{U_n/T_n}) \to 
H^1(U^{an}_n, \Omega^1_{U^{an}_n/T_n}) 
\to H^1(U^{an}_n(p), \Omega^1_{U^{an}_n(p)/T_n}).$$ 
If $[L_n] \in H^2(U^{an}, S_n)$ is sent to zero by 
the map 
$$H^2(U^{an}, S_n) \to H^2(U^{an}(p), S_n),$$  
then, by the diagram, 
$[L_n\vert_{U^{an}(p)}] = 0$. Thus, the local 
Kodaira-Spencer class  
$\theta^{loc}_{n+1}(p)$ vanishes.  
Let us consider the same diagram in the proof of 
Theorem 17.  

\begin{equation}
\begin{CD} 
\mathrm{Pic}({X}^{an})\otimes_{\mathbf{Z}}S_n @>>> 
\mathrm{Pic}({U}^{an})\otimes_{\mathbf{Z}}S_n @>{\cong}>> 
H^2(U^{an}, S_n) \\ 
@VVV @VVV @VVV \\ 
\mathrm{Pic}(X^{an}(p))\otimes_{\mathbf{Z}}S_n @>>> 
\mathrm{Pic}(U^{an}(p))\otimes_{\mathbf{Z}}S_n @>{\cong}>> 
H^2(U^{an}(p), S_n) 
\end{CD} 
\end{equation} 

Since $L_n$ is a line bundle of $(X_n)^{an}$ and 
$\mathrm{Pic}((X_n)^{an}) \cong \mathrm{Pic}(X^{an})$, 
$[L_n] \in H^2(U^{an}, S_n)$ comes from $\mathrm{Pic}(X^{an})$. 
If we take $X^{an}(p)$ small enough, then 
$\mathrm{Pic}(X^{an}(p)) = 0$. Hence, by the 
commutative diagram, we see that $[L_n] \in 
H^2(U^{an}, S_n)$ is sent to zero by the map 
$H^2(U^{an}, S_n) \to H^2(U^{an}(p), S_n)$. 
\vspace{0.2cm} 

\section{Symplectic varieties with good 
$\mathbf{C}^*$-actions} 

Let $X$ be a convex symplectic variety with 
terminal singularities and, in addition, with a 
$\mathbf{C}^*$-action. We put $Y := 
\mathrm{Spec}\;\Gamma (X, \mathcal{O}_X)$. 
Then the natural morphism $g: X \to Y$ 
is a $\mathbf{C}^*$-equivariant 
morphism. We assume that $Y$ has a good 
$\mathbf{C}^*$-action with a unique fixed point 
$0 \in Y$.  
By definition, $V := Y - \mathrm{Sing}(Y)$ 
admits a symplectic 2-form $\omega$; hence 
it gives a Poisson structure $\{\; , \;\}$ on $Y$. 
We assume that this Poisson structure has 
a positive weight $l > 0$ with respect to the 
$\mathbf{C}^*$-action, that is, 
$$ \mathrm{deg}\{a,b\} = \mathrm{deg}(a) + 
\mathrm{deg}(b) - l $$ 
for all homogenous elements $a, b \in \mathcal{O}_Y$. 
Now let us consider the Poisson deformation 
functor $\mathrm{PD}_{(X, \{ , \})}$ (cf. Section 3). 
By Theorem 14, it is pro-represented by a certain 
complete regular local $\mathbf{C}$-algebra 
$R = \lim R_n$ and a universal formal 
Poisson deformation $\{X^{univ}_n\}$ of $X$ over it. 
\vspace{0.15cm} 

{\bf Lemma 20}. 
{\em The $\mathbf{C}^*$-action on $X$ naturally induces 
a $\mathbf{C}^*$-action on $R$ and $\{X^{univ}_n\}$.} 
\vspace{0.12cm} 

{\em Proof}. Take an infinitesimal Poisson 
deformation $(X_S, \{\; , \;\}_S; \iota)$ of $X$ over 
$S = \mathrm{Spec}(A)$  
with $A/m = \mathbf{C}$. 
By definition, $\iota: 
X_S \otimes_A A/m \cong X$ is an identification 
of the central fiber with $X$.  
Since $X$ is a $\mathbf{C}^*$-variety, 
for each $\lambda \in \mathbf{C}^*$, we get an isomorphism  
$\phi_{\lambda}: X \to X$. By the assumption, $\phi_{\lambda}^*
\{\; ,\;\} = \lambda^l\{\; , \;\}$.  
Then $(X_S, \lambda^l\{\; , \;\}_S ; \phi_{\lambda}\circ \iota)$  
gives another Poisson deformation of $X$ over $S$. 
This operation naturally gives a $\mathbf{C}^*$-action 
on $R$ and $\{X_n\}$. Q.E.D.
\vspace{0.15cm}

We shall investigate the $\mathbf{C}^*$-action of $R$. 
In order to do that, take $L \in 
H^1(U^{an}, \mathcal{O}^*_{U^{an}})$ and consider 
the corresponding quasi-twistor deformation of $X$.  
We define a $\mathbf{C}^*$-action on $\mathbf{C}[[t]]$ 
so that $t$ has weight $l$. This induces a 
$\mathbf{C}^*$-action on each quotient ring $S_n := 
\mathbf{C}[t]/(t^{n+1})$. We put $T_n := \mathrm{Spec}(S_n)$. 
\vspace{0.15cm}

{\bf Lemma 21}. {\em Any quasi-twistor deformation 
$\{X_n\}$ of $X$ 
has a $\mathbf{C}^*$-action so that $\{X_n\} \to 
\{T_n\}$ is $\mathbf{C}^*$-equivariant.} 
\vspace{0.15cm} 

{\em Proof}. Let   
$R \to \mathbf{C}[[t]]$ be the surjection determined by 
our quasi-twistor 
deformation. We shall prove this map is 
$\mathbf{C}^*$-equivariant. 
For $\lambda \in \mathbf{C}^*$, let  
$\lambda^l : T_n \to T_n$  be the morphism 
induced by $t \to \lambda^l t$. We shall lift  
$\mathbf{C}^*$-actions of $X_n$ inductively. 
More explicitly, for each $\lambda \in 
\mathbf{C}^*$, we shall construct an isomorphism 
$\phi_{\lambda, n} : X_n \to X_n$ 
in such a way that:  

(i) the following diagram commutes 

\begin{equation}
\begin{CD} 
X_n @>{\phi_{\lambda, n}}>> X_n \\ 
@VVV @VVV \\ 
T_n @>{\lambda^l}>> T_n 
\end{CD} 
\end{equation} 
       
(ii) $(\phi_{\lambda, n})^* \{\; , \; \}_n = 
\lambda^l \{\; , \; \}_n $, and 
 
(iii) the collection $\{\phi_{\lambda, n}\}$, 
$\lambda \in \mathbf{C}^*$ gives a 
$\mathbf{C}^*$-action of $X_n$.
      
Suppose that it can be achieved. As in Lemma 20, 
let us fix an original identification 
$\iota: X_n \times_{T_n}T_0 \cong X$. 
Let 
$h_n : R \to S_n$  and $h^{\lambda}_n : 
R \to S_n$ be the maps determined by 
 $(X_n, \{\; , \;\}_n; \iota)$ 
and $(X_n, \lambda^l \{\; , \;\}_n; \phi_{\lambda}
\circ \iota)$ respectively.  
Let $\lambda \in \mathbf{C}^*$ act on 
$R_n$ as $\psi_{\lambda,n} : R_n \to R_n$. 
By definition, $h_n \circ \psi_{\lambda, n} 
= h^{\lambda}_n$.    
Then the 
existence of $\phi_{\lambda, n}$ implies that 
there is a commutative diagram  
 
\begin{equation}
\begin{CD} 
R_n @>{h_n}>> S_n \\ 
@V{\psi_{\lambda,n}}VV @V{\lambda^l}VV \\ 
R_n @>{h_n}>> S_n 
\end{CD} 
\end{equation} 
 
The construction of $\phi_{\lambda, n}$ goes as 
follows. We assume that $\phi_{\lambda, n-1}$  
already exist. Let $U_{n-1} \subset X_{n-1}$ be the 
locus where $X_{n-1} \to T_{n-1}$ is smooth. 
Let $\omega_{n-1} \in 
\Gamma (U_{n-1}, \Omega^2_{U_{n-1}/T_{n-1}})$ be the symplectic 
2-form corresponding to the Poisson structure $\{\; ,\; \}_{n-1}$. 
By the assumption, $X_{n-1} \to T_{n-1}$ is a $\mathbf{C}^*$-
equivariant morphism. The symplectic 2-form $\omega_{n-1}$ 
has weight $l$ with the induced $\mathbf{C}^*$-action on 
$U_{n-1}$. Let $L_{n-1} \in \mathrm{Pic}(U^{an}_{n-1})$ 
be the (unique) extension of $L \in \mathrm{Pic}(U^{an})$.  
By $T^1$-lifting principle, the extension of $X_{n-1}$ to 
$X_n$ is determined by an element 
$\theta_n \in \mathbf{H}^2(U_{n-1}, 
\wedge^{\cdot}\Theta_{U_{n-1}/T_{n-1}})$, where 
$\wedge^{\cdot}\Theta_{U_{n-1}/T_{n-1}}$ is the 
Lichnerowicz-Poisson complex defined in \S 2. 
The symplectic 2-form $\omega_{n-1}$ gives an 
identification (cf. \S 2): 
$$  \mathbf{H}^2(U_{n-1}, 
\wedge^{\cdot}\Theta_{U_{n-1}/T_{n-1}}) 
\cong 
\mathbf{H}^2(U_{n-1}, \Omega^{\geq 1}_{U_{n-1}/T_{n-1}}).$$ 
By definition of the twistor deformation, $\theta_n$ is sent to 
$[L_{n-1}]$. Since $[L_{n-1}]$ and  
$\omega_{n-1}$ have respectively weights $0$ and $l$ for  
the $\mathbf{C}^*$-action,  
$\theta_n$ should have weight $-l$. 
This is what we want. Q.E.D.  
\vspace{0.15cm}

Since any 1-st order Poisson deformation of $X$ 
is realized in a suitable quasi-twistor deformation, 
$\mathbf{C}^*$  has only weight $l$ on 
the maximal ideal $m_R$ of $R$.   
Thus, one can write $R$ as 
$\mathbf{C}[[t_1, ..., t_m]]$, where $t_i$ are 
all eigen-elements with weight $l$.   
Let $\hat{\mathcal{Y}} := \mathrm{Spec}\lim  
\Gamma (X^{univ}_n, \mathcal{O}_{X^{univ}_n})$. 
The $\mathbf{C}^*$-action 
on $\{X^{univ}_n\}$ induces a $\mathbf{C}^*$-action on 
$\hat{\mathcal{Y}}$. 
Let $\hat{\mathcal{X}} \to \hat{\mathcal{Y}}$ be the 
algebraization of $\{X^{univ}_n\}$ over $\hat{\mathcal{Y}}$. 
Since $Y$ and $R$ are both 
positively weighted, $\hat{\mathcal{Y}}$ is also positively weighted. 
The $\mathbf{C}^*$-action of the formal scheme $\{X^{univ}_n\}$ 
induces a $\mathbf{C}^*$-action of $\hat{\mathcal{X}}$  
in such a way that $\hat{\mathcal{X}} \to \hat{\mathcal{Y}}$ 
becomes $\mathbf{C}^*$-equivariant. 
\vspace{0.15cm}
 
{\bf Lemma 22}. {\em There is a a projective birational 
morphism of algebraic 
varieties with $\mathbf{C}^*$-actions 
$$\mathcal{X} \to \mathcal{Y}$$ over $\mathrm{Spec}\;\mathbf{C}
[t_1, ..., t_m]$ which is an algebraization of $\hat{\mathcal{X}} 
\to \hat{\mathcal{Y}}$. Moreover, $\mathcal{X}$ and 
$\mathcal{Y}$ 
admit natural Poisson structures over $\mathrm{Spec}\;\mathbf{C}
[t_1, ..., t_m]$.} 
\vspace{0.15cm}

{\em Proof}. Let $A$ be the completion of the coordinate 
ring of the affine scheme $\hat{\mathcal{Y}}$ at the 
origin. Then $A$ becomes a complete local ring with a 
good $\mathbf{C}^*$-action. The $\mathbf{C}^*$-equivariant 
projective morphism 
$\hat{\mathcal{X}} \to \hat{\mathcal{Y}}$ induces 
a $\mathbf{C}^*$-equivariant projective morphism 
$\hat{\mathcal{X}}_A 
\to \mathrm{Spec}(A)$. By Lemma A.8, 
there is a $\mathbf{C}^*$-linealized  
ample line bundle on $\hat{\mathcal{X}}_A$. 
By Lemma A.2, there is a $\mathbf{C}$-algebra $R$ 
of finite type such that $\hat{R} = A$. 
Put $\mathcal{Y} := \mathrm{Spec}(R)$. Since $R$ is 
generated by eigen-vectors (homogenous elements) 
of $\mathbf{C}^*$-action, 
$R$ contains $t_i$. So $R$ is a ring over 
$\mathbf{C}[t_1, ..., t_m]$. By Proposition A.5, 
there is a $\mathbf{C}^*$-equivariant projective 
morphism $\mathcal{X} \to \mathcal{Y}$ which algebraizes 
$\hat{\mathcal{X}}_A \to \mathrm{Spec}(A)$. 
This automatically algebraizes $\hat{\mathcal{X}} \to 
\hat{\mathcal{Y}}$.  
The complete local ring $A$ admits a Poisson structure 
over $\mathbf{C}[[t_1, ..., t_m]]$ 
induced by that of $\Gamma(\hat{\mathcal{Y}}, 
\mathcal{O}_{\hat{\mathcal{Y}}})$. This Poisson 
structure induces a Poisson structure of $R$ over 
$\mathbf{C}[t_1, ..., t_m]$ because, if $a, b \in A$ 
are homogenous, then $\{a,b\} \in A$ is again homogenous. 
The corresponding relative Poisson bi-vector $\Theta$ 
of $\mathcal{Y}$ is non-degenerate on the smooth part. 
Hence it defines a relative symplectic 2-form on the 
smooth part of $\mathcal{Y}$. This relative symplectic 
2-form is pulled back to $\mathcal{X}$ and defines 
a relative Poisson structure of $\mathcal{X}$. 
\vspace{0.2cm}

Let us fix an algebraic line bundle $L$ on $X$. 
Since $H^i(X, \mathcal{O}_X) = 0$ for $i = 1,2$, 
there is a unique line bundle $\hat{L} \in 
\mathrm{Pic}(\hat{\mathcal{X}})$ extending $L$. 
Let  
$\hat{L}_A \in \mathrm{Pic}(\hat{\mathcal{X}}_A)$  
be the pull-back of $\hat{L}$ to $\hat{\mathcal{X}}_A$. 
Since $\hat{L}$ is fixed by the $\mathbf{C}^*$-action 
of $\hat{\mathcal{X}}$, $\hat{L}_A$ is fixed by 
the $\mathbf{C}^*$-action of $\hat{\mathcal{X}}_A$.  
By Lemma A.8, for some 
$k > 0$, $(\hat{L}_A)^{\otimes k}$ is $\mathbf{C}^*$-linearized. 
By Proposition A.6, there is a $\mathbf{C}^*$-linearized 
line bundle on $\mathcal{X}$ extending $(\hat{L}_A)^{\otimes k}$. 
Thus, by replacing $L$ by its suitable multiple, 
we may assume that $L$ extends to a line bundle on 
$\mathcal{X}$.    
Let $U$ be the regular part of $X$ and let $[L] \in 
H^2(U^{an}, \mathbf{C})$ be the associated class with 
$L\vert_U$. Let us denote by $M$ the maximal ideal of 
$\mathbf{C}[t_1, ..., t_m]$ and identify $(M/M^2)^*$ with 
$H^2(U^{an}, \mathbf{C})$. Then $[L]$ can be written as a linear 
combination 
$$ a_1 t^*_1 + a_2 t^*_2 + ... + a_m t^*_m $$ 
with the dual base $\{t^*_i\}$ of $\{t_i\}$.  
Take a base change of $$\mathcal{X} \to \mathrm{Spec}\;\mathbf{C}
[t_1, ..., t_n]$$ by the map 
$$\mathrm{Spec}\;\mathbf{C}[t] \to  \mathrm{Spec}\;\mathbf{C}
[t_1, ..., t_n]$$ with $t_i = a_it$. Then we have a 1-parameter 
deformation $\mathcal{X}^L$ of $X$ over 
$T := \mathrm{Spec}\;\mathbf{C}[t]$. 
As we have shown in Lemma 21, this 
deformation gives an algebraization 
of the twistor deformation $X^L_{\infty} \to 
\mathrm{Spec}\;\mathbf{C}[[t]]$.  
We put $\mathcal{Y}^L := \mathrm{Spec}\;\Gamma (\mathcal{X}^L, 
\mathcal{O}_{\mathcal{X}^L})$. Now let us consider the 
birational projective morphism  
$$g_T: \mathcal{X}^L \to \mathcal{Y}^L$$ 
over $\mathrm{Spec}\;\mathbf{C}[t]$. Let $\eta \in 
\mathrm{Spec}\mathbf{C}[t]$ be the generic point and 
let $X^L_{\eta}$ and $Y^L_{\eta}$ be the generic fibers. 
Then we get a birational projective morphism 
$$ g_{\eta} : X^L_{\eta} \to Y^L_{\eta}.$$

{\bf Proposition 23} (Kaledin) {\em Assume that 
$X$ is smooth and $L$ is 
ample. Then 
$g_{\eta}: X^L_{\eta} \to Y^L_{\eta}$ is an isomorphism.} 
\vspace{0.12cm}

{\em Proof}. 
Denote by $T (\cong \mathrm{Spec}\mathbf{C}[t])$ 
the base space of our algebraized twistor deformation 
$\mathcal{X}^L$. Since $T$ has a good $\mathbf{C}^*$-action, 
$\mathcal{X}^L$ is smooth over $T$.  
The line bundle $L$ on $X$ uniquely extend a line bundle 
$\mathcal{L}$ on $\mathcal{X}^L$. 
Moreover, $\mathcal{X}^L$ is a Poisson $T$-scheme 
extending the original Poisson scheme $X$; thus, the 
symplectic 2-form $\omega$ on $X$ extends to a relative 
symplectic 2-form $\omega_T \in 
\Gamma (\mathcal{X}^L, \Omega^2_{\mathcal{X}^L/T})$. 
Let $\theta_T \in 
H^1(\mathcal{X}^L, \Theta_{\mathcal{X}^L/T})$ 
be the extension class (Kodaira-Spencer class) of the exact sequence 
$$0 \to (f_T)^*\Omega^1_{T/\mathbf{C}} 
\to \Omega^1_{\mathcal{X}^L / \mathbf{C}} \to 
\Omega^1_{\mathcal{X}^L/T} \to 0.$$  
By Lemma 16, we see that, in $H^1(\mathcal{X}^L, 
\Omega^1_{\mathcal{X}^L/T})$,    
$i(\theta_T)(\omega_T) = [\mathcal{L}]$.
We put $\omega_{\eta} := \omega_T\vert_{X^L_{\eta}}$,  
$\theta_{\eta} := \theta_T\vert_{X^L_{\eta}}$ 
and $L_{\eta} := \mathcal{L}\vert_{X_{\eta}}$. 
Then, in $H^1(X^L_{\eta}, \Omega^1_{X^L_{\eta}/k(\eta)})$, 
we have an equality: 
$$ i(\theta_{\eta})(\omega_{\eta}) = [L_{\eta}].$$ 
Since $g_{\eta}$ is a proper birational morphism, 
we only have to show that $X_{\eta}$ does not contain 
a proper curve defined over $k(\eta)$. Now let 
$\iota: C \to X_{\eta}$ be a morphism from a proper 
regular curve $C$ defined over $k(\eta)$ 
to $X_{\eta}$. We shall prove that $\iota(C)$ is a point.  
Let $\theta_C \in H^1(C, \Theta_{C/k(\eta)})$ be the 
Kodaira-Spencer class for $h: C \to \mathrm{Spec}k(\eta)$. 
In other words, $\theta_C$ is the extension class of the 
exact sequence 
$$ 0 \to h^*\Omega^1_{k(\eta)/\mathbf{C}} \to 
\Omega^1_{C/\mathbf{C}} \to \Omega^1_{C/k(\eta)} \to 0.$$ 
Then, by the compatibility of Kodaira-Spencer classes, we 
have 
$$i(\theta_C)(\iota^*\omega_{\eta}) = 
\iota^*(i(\theta_{\eta})(\omega_{\eta})).$$ 
The left hand side is zero because $\iota^*\omega_{\eta} 
= 0$. On the other hand, the right hand side is 
$\iota^*[L_{\eta}]$. Since $L$ is ample, $L_{\eta}$ is also ample. 
If $\iota (C)$ is not a point, 
then $\iota^*[L_{\eta}] \ne 0$, which is a contradiction.  
\vspace{0.15cm}

The following is a generalization of Proposition 23  
to the singular case. 
\vspace{0.15cm} 

{\bf Proposition 24}. {\em Assume that $X$ has only symplectic 
singularities and $L \in \mathrm{Pic}(X)$.} 

(a) {\em If $L$ is ample, then  $g_{\eta}: X^L_{\eta} \to Y^L_{\eta}$ 
is an isomorphism.} \vspace{0.12cm}

(b) {\em Let $X^+$ be another convex symplectic 
variety over $Y$ with terminal singularities and assume 
that $L$ becomes the proper transform of an ample line 
bundle $L^+$ on $X^+$. Then   
$g_{\eta}$ is a small birational morphism; in other words, 
$\mathrm{codim} \mathrm{Exc}(g_{\eta}) \geq 2$.} 
\vspace{0.12cm}

{\em Proof}. (i) We shall use the same notation as 
the proof of Proposition 23. We note that the Kodaira-Spencer 
class $\theta_T \in \mathrm{Ext}^1(\Omega^1_{\mathcal{X}^L/T}, 
\mathcal{O}_{\mathcal{X}^L})$ is contained in 
$H^1(\mathcal{X}^L, \Theta_{\mathcal{X}^L/T})$ because the 
twistor deformation is locally trivial by   
Theorem 19. Let $\mathcal{U} \subset \mathcal{X}^L$ be the 
locus where $\mathcal{X}^L \to T$ is smooth. Denote 
by $U_{\eta}$ the generic fiber of $\mathcal{U} \to T$. 
Let $(\theta_T)^0 \in H^0(\mathcal{U}, \Theta_{\mathcal{U}/T})$ 
be the restriction of $\theta_T$ to $\mathcal{U}$. 
The relative Poisson structure on $\mathcal{X}^L$ over $T$ gives an 
element $(\omega_T)^0 \in 
H^0(\mathcal{U}, \Omega^2_{\mathcal{U}/T})$. 
Note that, in general,  $(\omega_T)^0$ cannot 
extend to a global section of $\Omega^2_{\mathcal{X}^L/T}$.
Let $[\mathcal{L}]^0 \in H^1(\mathcal{U}, 
\Omega^1_{\mathcal{U}/T})$ be the class corresponding to 
a restricted line bundle $\mathcal{L}\vert_{\mathcal{U}}$.   
Then, $(\theta_T)^0$, $(\omega_T)^0$ and 
$[\mathcal{L}]^0$ defines respectively the 
classes $$\theta^0_{\eta} \in 
H^1(U_{\eta},\Theta^1_{U_{\eta}/k(\eta)}),$$ 
$$\omega^0_{\eta} \in H^0(U_{\eta}, \Omega^2_{U_{\eta}/k(\eta)})$$ 
and $$[L_{\eta}]^0 \in H^1(U_{\eta}, \Omega^1_{U_{\eta}/k(\eta)}).$$ 
We then have 
$$i(\theta^0_{\eta})(\omega^0_{\eta}) = [L_{\eta}]^0.$$

(ii)(Construction of a good resolution): 
We shall construct a good equivariant resolution 
of $\mathcal{X}^L$. In order to do that, first take   
an equivariant resolution $\pi_0: \tilde{X} \to X$ of  
$X$, that is, $(\pi_0)_*\Theta_{\tilde{X}} 
= \Theta_{X}$. Here $\Theta_X := \underline{\mathrm{Hom}}
(\Omega^1_X, \mathcal{O}_X)$. By Theorem 19, our twistor 
deformation gives us a sequence of locally trivial 
formal deformations of $X$: 
$$X \to X_1 \to ... \to X_n \to ... $$ 
We shall construct resolutions $\pi_n: \tilde{X}_n \to X_n$ 
inductively so that 
there is an affine open cover 
$X_n = \cup_{i \in I}U_{n,i}$ such that 
$(\pi_n)^{-1}(U_{n,i}) \cong (\pi_0)^{-1}(U)\times_{T_0}T_n$. 
Note that, if this could be done, then 
$(\pi_n)_*\Theta_{\tilde{X}_n/T_n} = 
\Theta_{X_n/T_n}$. Moreover, if we let 
$\tilde{\theta}_n \in H^1(\tilde{X}_{n-1}, \Theta_{\tilde{X}_{n-1}/T_{n-1}})$ 
be the Kodaira-Spencer class of $\tilde{X}_n \to T_n$, then 
$\tilde{\theta}_n$ coincides with the Kodaira-Spencer 
class $\theta_n \in H^1(X_{n-1}, \Theta_{X_{n-1}/T_{n-1}})$ of 
$X_n \to T_n$ because 
$\tilde{\theta}_n$ is mapped to zero by the map 
$$H^1(\tilde{X}_{n-1}, \Theta_{\tilde{X}_{n-1}/T_{n-1}}) 
\to H^0(X_n, R^1(\pi_{n-1})_*\Theta_{\tilde{X}_{n-1}/T_{n-1}}).$$   
Now assume that we are given such a resolution 
$\pi_n: \tilde{X}_n \to X_n$. Take the affine 
open cover $\{U_{n,i}\}_{i \in I}$ of $X_n$ as above. 
We put $\tilde{U}_{n,i} := (\pi_n)^{-1}(U_{n,i})$. 
For $i,j \in I$, there is an identification 
$U_{n,i}\vert_{U_{ij}} \cong U_{n,j}\vert_{U_{ij}}$ 
determined by $X_n$. For each $i \in I$, let $\mathcal{U}_{n,i}$ and 
$\tilde{\mathcal{U}}_{n,i}$ be trivial deformations 
of $U_{n,i}$ and $\tilde{U}_{n,i}$ over $T_{n+1}$ 
respectively. For each $i,j \in I$, take a $T_{n+1}$-isomorphism 
$$g_{ji}: \mathcal{U}_{n,i}\vert_{U_{ij}} \to 
\mathcal{U}_{n,j}\vert_{U_{ij}}$$ such that 
$g_{ji}\vert_{T_n} = id$. Then 
$$h_{ijk} := g_{ij}\circ g_{jk} \circ g_{ki}$$ 
gives an automorphism of $\mathcal{U}_{n,i}\vert_{U_{ijk}}$  
over $T_{n+1}$ such that $h_{ijk}\vert_{T_n} = id$. 
Since $\pi_n : \tilde{X}_n \to X_n$ is an equivariant 
resolution, $g_{ij}$ extends uniquely to 
$$\tilde{g}_{ij}: \tilde{\mathcal{U}}_{n,i}\vert_{U_{ij}} 
\cong \tilde{\mathcal{U}}_{n,j}\vert_{U_{ij}}.$$ 
One can consider $\{h_{ijk}\}$ as a 2-cocycle of 
the \v{C}ech cohomology of $\Theta_X$; hence it gives an 
element $ob \in H^2(X, \Theta_X)$. But, since $X_n$ 
extends to $X_{n+1}$, $ob = 0$. Therefore, by modifying 
$g_{ij}$ to $g'_{ij}$ suitably, one can get 
$$g'_{ij}\circ g'_{jk} \circ g'_{ki} = id.$$  
Then $$\tilde{g}'_{ij}\circ \tilde{g}'_{jk} \circ 
\tilde{g}'_{ki} = id.$$ 
Now $\tilde{X}_n$ also extends to $\tilde{X}_{n+1}$ 
and the following diagram commutes: 

\begin{equation}
\begin{CD} 
\tilde{X}_n  @>>> \tilde{X}_{n+1} \\ 
@VVV @VVV \\ 
X_n @>>> X_{n+1} 
\end{CD} 
\end{equation} 

By Th\'{e}or\`{e}me (5.4.5) of [EGA III], one has an 
algebraization $\tilde{X}^L_{\infty} \to 
Y_{\infty}$ of 
$\{\tilde{X}_n \to Y_n\}$. Moreover, the morphism 
$\{\pi_n: \tilde{X}_n \to X_n\}$ induces 
$\pi_{\infty}: \tilde{X}^L_{\infty} \to X^L_{\infty}$. 
By the construction, the $\mathbf{C}^*$-action 
on $X^L_{\infty}$ lifts to $\tilde{X}^L_{\infty}$. 
Then $\tilde{X}^L_{\infty} \to Y_{\infty}$ 
is algebraized 
to a $\mathbf{C}^*$-equivariant projective morphism  
$\tilde{\mathcal{X}}^L \to \mathcal{Y}^L$ in such a 
way that it factors through $\mathcal{X}^L$.       

(iii) Let $\pi : \tilde{\mathcal{X}}^L \to \mathcal{X}^L$ 
be the equivariant resolution constructed in (ii). 
Let us denote by $\tilde{X}_{\eta}$ the generic fiber of 
$\tilde{\mathcal{X}}^L \to T$. 
This resolution gives 
an equivariant resolution $\pi_{\eta}: \tilde{X}_{\eta} 
\to X_{\eta}$. In particular, 
$(\pi_{\eta})_*\Theta_{\tilde{X}/k(\eta)} = 
\Theta_{X_{\eta}/k(\eta)}$. 
Let $\tilde{\theta}_{\eta} \in 
H^1(\tilde{X}_{\eta}, \Theta_{\tilde{X}_{\eta}/k(\eta)})$ 
be the Kodaira-Spencer class for $\tilde{X}_{\eta} 
\to \mathrm{Spec}k(\eta)$. Then the Kodaira-Spencer 
class $\theta_{\eta} \in H^1(X_{\eta}, \Theta_{X_{\eta}/k(\eta)})$ 
for $X_{\eta}$ coicides with $\tilde{\theta}_{\eta}$ by 
the natural injection $H^1(X_{\eta}, \Theta_{X_{\eta}/k(\eta)}) 
\to H^1(\tilde{X}_{\eta}, \Theta_{\tilde{X}_{\eta}/k(\eta)})$. 
Let  $i_{\eta}: U_{\eta} \to X_{\eta}$ be the embedding 
of the regular part. Since $(i_{\eta})_*\Omega^2_{U_{\eta}/k(\eta)} 
\cong (\pi_{\eta})_*\Omega^2_{\tilde{X}_{\eta}/k(\eta)}$, by 
[Fl], $\omega^0_{\eta}$ extends to 
$$ \omega_{\eta} \in \Gamma (X_{\eta},  
(\pi_{\eta})_*\Omega^2_{\tilde{X}_{\eta}/k(\eta)}).$$

(iv) We have a pairing map: 
$$H^0(X_{\eta}, (\pi_{\eta})_*\Omega^2_{\tilde{X}_{\eta}/k(\eta)}) 
\times H^1(X_{\eta}, (\pi_{\eta})_*\Theta_{\tilde{X}_{\eta}/k(\eta)}) 
\to H^1(X_{\eta}, (\pi_{\eta})_*\Omega^1_{\tilde{X}_{\eta}/k(\eta)}).$$ 
Denote by $i(\theta_{\eta})(\omega_{\eta})$ the image of $(\omega_{\eta}, 
\theta_{\eta})$ by this pairing map. By pulling back $L_{\eta}$ by 
$\pi_{\eta}$, one can define a class 
$[L_{\eta}] \in H^1(X_{\eta}, 
(\pi_{\eta})_*\Omega^1_{\tilde{X}_{\eta}k(\eta)})$. 
Let us consider the exact sequence 
$$H^1_{\Sigma}(X_{\eta}, (\pi_{\eta})_*\Omega^1_{\tilde{X}_{\eta}/k(\eta)}) 
\to H^1(X_{\eta}, (\pi_{\eta})_*\Omega^1_{\tilde{X}_{\eta}/k(\eta)}) 
\to H^1(U_{\eta}, \Omega^1_{U_{\eta}/k(\eta)}),$$ 
where $\Sigma := X_{\eta}\setminus U_{\eta}$. 
Since $(\pi_{\eta})_*\Omega^1_{\tilde{X}_{\eta}/k(\eta)} \cong 
(i_{\eta})_*\Omega^1_{U_{\eta}/k(\eta)}$ by [Fl], it is a reflexive 
sheaf. A reflexive sheaf on $X_{\eta}$ is locally written as 
the kernel of a homomorphism from a free sheaf to a torsion 
free sheaf. Since $X_{\eta}$ is Cohen-Macaulay and 
$\mathrm{Codim}(\Sigma \subset X_{\eta}) \geq 2$, we have 
$H^1_{\Sigma}(X_{\eta}, (\pi_{\eta})_*\Omega^1_{\tilde{X}_{\eta}/k(\eta)}) 
= 0$. We already know in (i) that $[L_{\eta}]^0 = 
i(\theta^0_{\eta})(\omega^0_{\eta})$ in 
$H^1(U_{\eta}, \Omega^1_{U_{\eta}/k(\eta)})$. Therefore, by 
the exact sequence, we see that 
$$[L_{\eta}] = i(\theta_{\eta})(\omega_{\eta}).$$  
   
(v) Consider the pairing map 
$$H^0(\tilde{X}_{\eta}, \Omega^2_{\tilde{X}_{\eta}/k(\eta)}) 
\times H^1(\tilde{X}_{\eta}, \Theta_{\tilde{X}_{\eta}/k(\eta)}) 
\to H^1(\tilde{X}_{\eta}, \Omega^1_{\tilde{X}_{\eta}/k(\eta)}).$$ 
By the construction of $\tilde{X}_{\eta}$, 
the Kodaira-Spencer class $\tilde{\theta}_{\eta}$ of 
$\tilde{X}_{\eta} \to \mathrm{Spec}k(\eta)$ coincides with 
the Kodaira-Spencer class $\theta_{\eta}$. Hence 
$(\omega_{\eta}, \tilde{\theta}_{\eta})$ is sent to 
$(\pi_{\eta})^*[L_{\eta}]$ by the pairing map. 

Now we shall prove (a). 
For a proper regular curve $C$ defined over $k(\eta)$, 
assume that there is a $k(\eta)$-morphism 
$\iota: C \to \tilde{X}_{\eta}$. We shall prove that 
$(\pi_{\eta})\circ\iota(C)$ is a point. 
By the compatibility of the Kodaira-Spencer classes, 
we have 
$$i(\theta_C)(\iota^*\omega_{\eta}) = 
\iota^*(i(\omega_{\eta})(\tilde{\theta}_{\eta})).$$ 
The left hand side is zero because 
$\iota^*\omega_{\eta} = 0$. The right hand side 
is $\iota^*(\pi_{\eta})^*[L_{\eta}]$ as we just 
remarked above. If $\pi_{\eta}\circ\iota(C)$ is 
not a point, then this is not zero because $L_{\eta}$ 
is ample; but this is a contradiction. 

Next we shall prove (b). 
We shall derive a contradiction assuming that 
$g_T: \mathcal{X}^L \to \mathcal{Y}^L$ is a divisorial 
birational contraction. We put $\mathcal{E} := 
\mathrm{Exc}(g_T)$. By the assumption, there 
is another convex symplectic variety $X^+$ over $Y$, and 
$X$ and $X^+$ are isomorphic in codimension one over $Y$. 
Let $F \subset X$ (resp. $F^+ \subset X^+$) be the locus 
where the birational map $X --\to X^+$ is not an 
isomorphism. Then    
$\mathrm{codim}(F \subset X) \geq 2$. We shall prove that 
$(L, \bar{C}) > 0$ for any proper irreducible curve $\bar{C}$ which 
is not contained in $F$. Let $\bar{C}$ be such a curve. Take a 
common resolution $\mu : Z \to X$ and $\mu^+: Z \to X^+$. 
We may assume that $\mathrm{Exc}(\mu)$ is a union of irreducible 
divisors, say $\{E_i\}$. 
Since $X$ and $X^+$ are isomorphic in codimension one, 
$\mathrm{Exc}(\mu) = \mathrm{Exc}(\mu^+)$.  
On can write $(\mu^+)^*L^+ = \mu^*L - \Sigma a_iE_i$ with 
{\em non-negative} integers $a_i$. In fact, if $a_j < 0$ for 
some $j$, then $E_j$ should be a fixed component of the 
linear system $\vert (\mu^+)^*L^+ \vert$; but this is a 
contradiction since $L^+$ is (very) ample. One can find 
a proper curve $D$ on $Z$ such that $\mu(D) = \bar{C}$ and 
such that $C^+ := \mu^+(D)$ is an irreducible curve on $X^+$. 
(i.e. $\mu^+(D)$ is not reduced to a point.) Note that $D$ is 
not contained in any $E_i$. 
Then 
$$(L, \bar{C}) = (\mu^*L, D) = ((\mu^+)^*L^+ + \Sigma a_iE_i, D) 
> 0.$$ 
Let us consider  
all effective 1-cycles on $X$ which are contracted to 
points by $g$ and are obtained as the limit of 
effective 1-cycles on $X^L_{\eta}$. Since $\mathcal{E}$ has 
codimension 1 in $\mathcal{X}^L$, one can find such an   
effective 1-cycle whose support intersects $F$ at most 
in finite points. In other words, there is a flat family 
$\mathcal{C} \to T$ 
of proper curves in $\mathcal{X}^L/T$ in such a way 
that any irreducible component of 
$\mathcal{C}_0 := \mathcal{C} \cap X$ is not contained in $F$.  
Let $\bar{C}_{\eta}$ be the generic fiber of 
$\mathcal{C}\to T$. Take a regular proper curve 
$C$ over $k(\eta)$ and a $k(\eta)$-morphism 
$\iota : C \to \tilde{X}_{\eta}$ so that 
$\pi_{\eta}\circ\iota(C) = \bar{C}_{\eta}$. 
By the definition of $\mathcal{C}$, 
$(L_{\eta}, \bar{C}_{\eta}) > 0$. 
Now one can get a contradiction by using this curve $C$ 
in the similar way to (a).     
\vspace{0.15cm}

{\bf Corollary 25}. {\em Let $Y$ be an affine 
symplectic variety with a good $\mathbf{C}^*$-action 
and assume that the Poisson structure of $Y$ is positively 
weighted.  
Let 
$$ X \stackrel{f}\rightarrow Y \stackrel{f'}\leftarrow 
X'$$ 
be a diagram such that, } 
\begin{enumerate}
\item $f$ (resp. $f'$) is a crepant, birational, projective 
morphism. 
\item $X$ (resp. $X'$) has only terminal singularities. 
\item $X$ (resp. $X'$) is {\bf Q}-factorial. 
\end{enumerate} 
{\em Then both $X$ and $X'$ have locally trivial deformations 
to an affine variety $Y_t$ obtained as a Poisson  
deformation of $Y$. In particular, $X$ and $X'$ have the same 
kind of singularities.} 
\vspace{0.12cm}

{\em Proof}. (i) By Step 1 of the proof of Proposition A.7, 
the $\mathbf{C}^*$-action of $Y$ lifts to $X$ and $X'$. 
So we are in the situation of section 4. 
Since $Y$ is a symplectic variety, outside 
certain locus at least of codimension 4 (say $\bar{\Sigma}$), 
its singularity is locally isomorphic to the product 
$(\mathbf{C}^{n-2}, 0) \times (S,0)$ (as an analytic 
space). Here $(S,0)$ is the germ of a rational double 
point singularity of a surface (cf. [Ka 2]). 
We put $\bar{V} := Y - \bar{\Sigma}$. Since $f$ and $f'$ 
are both (unique) minimal resolutions of rational double 
points over $\bar{V}$, $f^{-1}(\bar{V}) \cong (f')^{-1}(\bar{V})$.    

(ii) Fix an ample line bundle $L$ of $X$ and 
let $\{X_n\}$ be the twistor deformation associated with 
$L$.  
This induces a formal deformation $\{Y_n\}$ of $Y$. 
Let $L'$ be the proper transform of $L$ by 
$X --\to X'$. Since $X'$ is {\bf Q}-factorial, we may 
assume that $L'$ is a line bundle of $X'$.\footnote{The 
twistor deformation associated with $L$ and the one associated 
with $L^{\otimes m}$ are essentially the same. The latter 
one is obtained from the first one just by changing  
the parameters $t$ by $mt$.}   Let $\{X'_n\}$ be the 
twistor deformation of $X'$ associated with $L'$. 
This induces a formal deformation $\{Y'_n\}$ of $Y$.  
\vspace{0.2cm}

{\bf Lemma 26}. {\em The formal deformation 
$\{Y'_n\}$ coincides with $\{Y_n\}$.} 
\vspace{0.12cm}

{\em Proof}. The formal deformation $\{X_n\}$ 
of $X$ induces a formal deformation of 
$W := f^{-1}(\bar{V})$, say $\{W_n\}$. 
The deformation induces a formal deformation 
$\{\bar{V}_n\}$ of $\bar{V}$ by $\bar{V}_n := 
\mathrm{Spec}\Gamma 
(W_n, \mathcal{O}_{W_n})$ because 
$R^1(f\vert_W)_*\mathcal{O}_{W} = 0$ and $(f\vert_W)_*\mathcal{O}_W 
= \mathcal{O}_{\bar{V}}$ (cf. [Wa]). 
Since $\bar{V} = Y - \bar{\Sigma}$ with $\mathrm{codim}(\bar{\Sigma} 
\subset Y) \geq 4$, the formal deformation $\{\bar{V}_n\}$ of $\bar{V}$ 
extends uniquely to that of $Y$ (cf. Proposition 13, (1)). 
This extended deformation is nothing but $\{Y_n\}$. 
On the other hand, the formal deformation $\{X'_n\}$ of $X'$ 
induces a formal deformation of $W' := (f')^{-1}(\bar{V})$, say 
$\{W'_n\}$. As remarked in (i), $W \cong W'$. Moreover, by 
Corollary 10, the 
Poisson deformations of $W$ (resp. $W'$) are controlled by 
the cohomology 
$H^2(W^{an}, \mathbf{C})$ (resp. $H^2((W')^{an}, \mathbf{C})$) 
because $H^i(\mathcal{O}_W) = 0$ for $i = 1, 2$ (resp. 
$H^i(\mathcal{O}_{W'}) = 0$ for $i = 1, 2$).    
Since $L'$ is the proper transform of $L$, 
$[L\vert_W]$ is sent to $[L'\vert_{W'}]$ by the natural 
identification $H^2(W^{an}, \mathbf{C}) \cong 
H^2({W'}^{an}, \mathbf{C})$. This implies that $\{W_1\}$ 
and $\{W'_1\}$ coincide. By the construction of $\{X_n\}$ 
(resp. $\{X'_n\}$), $L$ (resp. $L'$) extends uniquely to 
$L_n$ (resp. $L'_n$). Then $[L_1\vert_{W}] \in 
H^2(W^{an}, S_1)$ is sent to $[L'_1\vert_{W'}] 
\in H^2(((W')^{an}, S_1)$, which implies that 
$W_2$ and ${W'}_2$ coincide. By the similar inductive 
process, one concludes that $\{W_n\}$ and $\{W'_n\}$ 
coincide.    
The formal deformation $\{W'_n\}$ of $W'$ 
induces a formal deformation $\{\bar{V}'_n\}$ of $\bar{V}$, 
which coincides with $\{\bar{V}_n\}$. So the extended 
deformation $\{Y'_n\}$ also coincides with $\{Y_n\}$.  
\vspace{0.15cm} 

(iii)  
Let 
$$ \mathcal{X}^L \rightarrow \mathcal{Y} \leftarrow 
(\mathcal{X}')^L$$ 
be the algebraizations of 
$$ \{X_n\} \rightarrow \{Y_n\} \leftarrow 
\{X'_n\}$$ 
over $T$. Let $\eta \in T$ 
be the generic point. Then, by Proposition 24, (a),  
$X^L_{\eta} \cong Y^L_{\eta}$. Since $X$ is 
{\bf Q}-factorial, we have: 
\vspace{0.12cm} 

{\bf Lemma 27}. {\em $X^L_{\eta}$ is also 
{\bf Q}-factorial.} 
\vspace{0.12cm}

{\em Proof}. Let $D$ be a Weil divisor of $X^L_{\eta}$. 
One can extend $D$ to a Weil 
divisor $\bar{D}$ of $\mathcal{X}^L$ by taking its closure. 
The restriction of $\bar{D}$ to $X$ defines a Weil 
divisor $\bar{D}\vert_X$. Note that the support of 
$\bar{D}\vert_X$ is $\bar{D} \cap X$ and the multiplicity 
on each irreducible component is well determined 
because $\bar{D}$ is a Cartier divisor at a regular 
point of $X$. Let $m >0$ be an integer such that 
$m(\bar{D}\vert_X)$ is a Cartier divisor. 
Let $\mathcal{O}(m\bar{D})$ be the 
reflexive sheaf associated with $m\bar{D}$ 
and let $\mathcal{O}(m\bar{D}\vert_X)$ be the 
line bundle associated with $m\bar{D}\vert_X$.  
By [K-M, Lemma (12.1.8)], $$\mathcal{O}(m\bar{D})
\otimes_{\mathcal{O}_{\mathcal{X}^L}}\mathcal{O}_X 
= \mathcal{O}({m\bar{D}\vert_X}).$$ 
In particular, $\mathcal{O}(m\bar{D})$ is 
a line bundle around $X$. Therefore, $m\bar{D}$ is a 
Cartier divisor on some Zariski open neighborhood 
of $X \subset \mathcal{X}^L$. Let $Z$ be the non-Cartier 
locus of $m\bar{D}$. Since $\mathcal{O}(m\bar{D})$ 
is fixed by the $\mathbf{C}^*$-action on $\mathcal{X}^L$, 
$Z$ is stable under the $\mathbf{C}^*$-action. 
Since $f_T: \mathcal{X}^L  
\to \mathcal{Y}$ is a projective morphism, $f_T(Z)$ 
is a closed subset of $\mathcal{Y}$. 
Since $Y \cap f_T(Z) = \emptyset$ and $Y$ 
has a good $\mathbf{C}^*$-action, $f_T(Z)$ should 
be empty; hence $Z$ should be also empty.  Q.E.D. 
\vspace{0.15cm}

Since $X$ and $X'$ are both crepant partial 
resolutions (with terminal singularities) of $Y$, 
they are isomorphic in codimension one. Now one 
can apply Proposition 24, (b) to the twistor deformation 
$\{X'_n\}$ of $X'$. 
Then we conclude that $(X')^{L'}_{\eta} 
\to Y_{\eta}$ is a {\em small} birational projective 
morphism.  On the other hand, $Y_{\eta} (\cong X^L_{\eta})$ 
is {\bf Q}-factorial by Lemma 27. 
These imply that $(X')^L_{\eta} \cong Y_{\eta}$. 
By Theorem 19, $\mathcal{X}^L \to T$ and 
$(\mathcal{X}')^{L'} \to T$ are locally trivial 
deformations of $X$ and $X'$ respectively.    
\vspace{0.2cm}

{\bf Corollary 28}. {\em Let $Y$ be an affine 
symplectic variety with a good $\mathbf{C}^*$-action. 
Assume that the Poisson structure of $Y$ is positively weighted, and 
$Y$ has only terminal singularities.  
Let 
$ f: X \rightarrow Y $ 
be a crepant, birational, projective 
morphism such that $X$ has only terminal singularities and 
such that $X$ is {\bf Q}-factorial. 
Then the following are equivalent.} 
\vspace{0.12cm} 

(a) {\em $X$ is non-singular.} 

(b) {\em $Y$ is smoothable by a Poisson 
deformation.} 
\vspace{0.12cm}

{\em Proof}. 
First of all, the $\mathbf{C}^*$-action of $Y$ lifts 
to $X$ by Step 1 of the proof of Proposition A.7. 
Secondly, by Corollary A.10, $X^{an}$ is $\mathbf{Q}$-factorial. 
We regard 
$X$ and $Y$ as Poisson schemes. 
The Poisson deformation functors $\mathrm{PD}_X$ and  
$\mathrm{PD}_Y$ have 
pro-representable hulls $R_X$ and $R_Y$  
respectively (Theorem 14). We put 
$U := (X)_{reg}$ and  
$V := Y_{reg}$. Then, by Lemma 12, 
$\mathrm{HP}^2(U) = H^2(U^{an}, \mathbf{C})$ 
and $\mathrm{HP}^2(V) = H^2(V^{an}, \mathbf{C})$. 
Note that, by Proposition 13, they coincide with 
$\mathrm{PD}_X(\mathbf{C}[\epsilon])$ and 
$\mathrm{PD}_Y(\mathbf{C}[\epsilon])$, respectively. 
By the proof of [Na, Proposition 2], we see that 
$$(*): H^2(U^{an}, \mathbf{C}) \cong H^2(V^{an}, \mathbf{C}).$$
Let $l >0$ be the weight of Poisson structure on $Y$. 
Then one can get universal $\mathbf{C}^*$-equivariant 
Poisson deformations 
$\mathcal{X}$ and $\mathcal{Y}$ over the same affine 
base $B := \mathrm{Spec}\;\mathbf{C}[t_1, ..., t_m]$, where 
$m = h^2(U^{an}, \mathbf{C})$ and each $t_i$ has weight 
$l$. 
By Theorem 17 $\mathcal{X} \to B$ is a locally trivial 
deformations of $X$. 
The birational projective morphism $f$ induces a 
birational projective $B$-morphism 
$$ f_B: \mathcal{X} \to \mathcal{Y}. $$ 
Let $\eta$ be the generic point of $B$. Take the generic fibers 
over $\eta$.  
Then we have 
$$\mathcal{X}_{\eta} \stackrel{f_{\eta}}\rightarrow 
\mathcal{Y}_{\eta}.$$  
Every twistor deformation of $X$ associated 
with an ample line bundle $L$ determines a 
(non-closed) point 
$\zeta_L \in B$. By Proposition 24, 
$f_{\zeta_L}$ is an isomorphism. This 
implies that $f_{\eta}$ is an isomorphism. 
Therefore, $\mathcal{Y}_{\eta}$ is regular 
if and only if $X$ is non-singular.  
\vspace{0.2cm} 

\section{General cases} 
Let $X$ be a convex symplectic variety with 
terminal singularities. 
Let $\{X_n\}$ be a twistor deformation for
 $L \in \mathrm{Pic}(X)$.  We put 
$Y_n := \mathrm{Spec} \Gamma (X_n, \mathcal{O}_{X_n})$ 
and $Y^L_{\infty} := \mathrm{Spec} \lim \Gamma (X_n, 
\mathcal{O}_{X_n})$. As in \S 3, $\{X_n\}$ is 
algebraized to $g_{\infty}: X^L_{\infty} \to 
Y^L_{\infty}$ over 
$T_{\infty}$, where $T_{\infty} := 
\mathrm{Spec}\; \mathbf{C}[[t]]$.   
We do not know, however, as in \S 3, that $\{X_n\}$ can be 
algebraized to $g_T : \mathcal{X}^L \to 
\mathcal{Y}^L$ over $T := \mathrm{Spec}\; \mathbf{C}[t]$.  
Let $\eta_{\infty} \in \mathrm{Spec}\; \mathbf{C}[[t]]$ 
be the generic point and let $g_{\eta_{\infty}} : 
X_{\eta_{\infty}} \to Y_{\eta_{\infty}}$ be the morphism 
between the generic fibers induced by $g_{\infty}$. 
\vspace{0.15cm}

{\bf Proposition 29}.   
(a) {\em If $L$ is ample, then  $g_{\eta_{\infty}}: 
X^L_{\eta_{\infty}} \to Y^L_{\eta_{\infty}}$ 
is an isomorphism.} \vspace{0.12cm} 

(b) {\em Let $X^+$ be another convex symplectic 
variety over $Y$ with terminal singularities and assume 
that $L$ becomes the proper transform of an ample line 
bundle $L^+$ on $X^+$. Then 
$g_{\eta_{\infty}}$ is a small birational morphism; in other words, 
$\mathrm{codim} \mathrm{Exc}(g_{\eta_{\infty}}) \geq 2$.} 
\vspace{0.15cm} 

{\em Proof}. The idea of the proof is the same as Proposition 
24. But we need more delicate argument because 
neither $X^L_{\infty}$ or $Y^L_{\infty}$ is of finite type 
over $T_{\infty}$. 
First of all, we should replace the usual differential sheaves 
$\Omega^i_{X^L_{\infty}/T_{\infty}}$  $(i \geq 1)$,  
$\Omega^1_{X^L_{\infty}/\mathbf{C}}$, 
and $\Omega^1_{T_{\infty}/\mathbf{C}}$ respectively 
by 
$\hat{\Omega}^i_{X^L_{\infty}/T_{\infty}}$, 
$\hat{\Omega}^1_{X^L_{\infty}/\mathbf{C}}$, 
and $\hat{\Omega}^1_{T_{\infty}/\mathbf{C}}$. 
Here 
$\hat{\Omega}^i_{X^L_{\infty}/T_{\infty}}$ 
is a coherent sheaf on $X^L_{\eta_{\infty}}$ determined as 
the limit of the formal sheaves  
$\{\Omega^i_{X_n/T_n}\}$, 
$\hat{\Omega}^1_{X^L_{\infty}/\mathbf{C}}$ is 
a coherent sheaf on $X^L_{\infty}$ determined 
as the limit of $\{\Omega^1_{X_{n+1}/\mathbf{C}}\vert_{X_n}\}$ 
and $\hat{\Omega}^1_{T_{\infty}/\mathbf{C}}$ is a 
coherent sheaf on $T_{\infty}$ determined as the 
limit of $\{\Omega^1_{T_{n+1}/\mathbf{C}}\vert_{T_n}\}$. 
Now the Kodaira-Spencer class $\theta_{T_{\infty}} \in 
\mathrm{Ext}^1(\hat{\Omega}^1_{X^L_{\infty}/T_{\infty}}, 
\mathcal{O}_{X^L_{\infty}})$ for $X^L_{\infty} \to T_{\infty}$ 
is the extension class of the exact sequence 
$$ 0 \to (f_{\infty})^*\hat{\Omega}^1_{T_{\infty}/\mathbf{C}} 
\to \hat{\Omega}^1_{X^L_{\infty}/\mathbf{C}} 
\to \hat{\Omega}^1_{X^L_{\infty}/T_{\infty}} \to 0.$$ 
Then, as in Proposition 24,  
we can construct a good resolution $\pi_{\infty}: 
\tilde{X}^L_{\infty} \to X^L_{\infty}$ of $X^L_{\infty}$. 
Let $E_{\infty}$ be the exceptional 
locus of $g_{\infty}$. Assume that 
$f_{\infty}(E_{\infty})$ 
contains a generic point $\eta_{\infty} \in T_{\infty}$. 
By cutting $E_{\infty}$ by  
$g_{\infty}$-very ample divisors and by the pull-back of 
suitable divisors on $Y_{\infty}$, we can find an integral subscheme 
$\bar{C}_{\infty} \subset X^L_{\infty}$ of dimension 2 
such that $g_{\infty}(\bar{C}_{\infty}) 
\to T_{\infty}$ is a finite surjective morphism. 
Note that $\bar{C}_{\infty} \to T_{\infty}$ is a flat projective 
morphism with fiber dimension 1. Take a desingularization 
$C_{\infty} \to \bar{C}_{\infty}$ which factors through 
$\tilde{X}^L_{\infty}$.  We put $C_n := 
C_{\infty}\times_{T_{\infty}}T_n$, and $C_{\eta_{\infty}} := 
C_{\infty}\times_{T_{\infty}}\mathrm{Spec}\; k(\eta_{\infty})$. 
Then $C_{\eta_{\infty}}$ 
is a proper regular curve over $k(\eta_{\infty})$. 
Moreover, one can define 
$\hat{\Omega}^1_{C_{\infty}/\mathbf{C}}$ as the limit of 
the formal sheaf $\{\Omega^1_{C_n/\mathbf{C}}\}$. Then, 
the Kodaira-Spencer class $\theta_{C_{\eta_{\infty}}}$ for 
$C_{\eta_{\infty}} 
\to \mathrm{Spec}\; k(\eta_{\infty})$ is well-defined as an element of 
$H^1(C_{\eta_{\infty}}, 
\Theta_{C_{\eta_{\infty}}/k(\eta_{\infty})})$. Then, the final 
argument in the proof of Proposition 24 is valid in 
our case.  
\vspace{0.2cm}

The same argument of Proposition 25 now yields: 
\vspace{0.12cm}

{\bf Corollary 30}. {\em Let $Y$ be an affine 
symplectic variety. 
Let 
$$ X \stackrel{f}\rightarrow Y \stackrel{f'}\leftarrow 
X'$$ 
be a diagram such that, } 
\begin{enumerate}
\item $f$ (resp. $f'$) is a crepant, birational, projective 
morphism. 
\item $X$ (resp. $X'$) has only terminal singularities. 
\item $X$ (resp. $X'$) is {\bf Q}-factorial. 
\end{enumerate} 
{\em Then, there is a 
flat deformation 
$$ X_{\infty} \rightarrow Y_{\infty} \leftarrow 
X'_{\infty}$$ over $T_{\infty} := \mathrm{Spec} \mathbf{C}[[t]]$ 
of the original diagram $X \rightarrow Y \leftarrow X'$ 
such that} 

(i){\em  $X_{\infty} \to T_{\infty}$ and $X'_{\infty} \to T_{\infty}$ 
are both locally trivial deformations, and}  

(ii) {\em the generic fibers are all isomorphic: 
$$ X_{\eta} \cong Y_{\eta} \cong X'_{\eta} $$ 
for the generic point $\eta \in T_{\infty}$.} 
\vspace{0.2cm}

In Corollary 30, $X_{\eta}$ (res. $X'_{\eta}$) is 
not of finite type over $k(\eta)$. So, at this moment, 
it is not clear how the singularities of $X$ are 
related to those of $X'$. However, one can say more 
when $X$ is smooth: 
\vspace{0.2cm} 

{\bf Corollary 31}. {\em With the same assumption as 
Corollary 30, if $X$ is non-singular, then $X'$ is 
also non-singular.} 
\vspace{0.15cm}

{\em Proof}. Since $X$ is non-singular, $X_{\infty}$ 
is formally smooth over $\mathbf{C}$. Since 
$X_{\eta} \cong Y_{\eta}$, $Y_{\infty}$ is formally 
smooth over $\mathbf{C}$ outside $Y$. By [Ar 2, Theorem 3.9] 
(see also [Hi], [Ri], [Ka 3]), for each closed point 
$p \in Y$, there is an etale map $Z_{\infty} \to Y_{\infty}$ 
whose image contains $p \in Y_{\infty}$, and $Z_{\infty} \to 
T_{\infty}$ is algebraized to $\mathcal{Z} \to T$. Here 
$T = \mathrm{Spec}\; \mathbf{C}[t]$. 
The completion $\hat{Z}$ of $\mathcal{Z}$ along 
the closed fiber coincides with $Z_{\infty}$. 
The diagram 
$$ X_{\infty} \to Y_{\infty} \leftarrow X'_{\infty}$$   
is pulled back by the map  
$Z_{\infty} \to Y_{\infty}$ to  
$$ X_{\infty} \times_{Y_{\infty}}Z_{\infty} \to 
Z_{\infty} \leftarrow X'_{\infty} \times_{Y_{\infty}}Z_{\infty}.$$  
Take generic fibers of this diagram over $T_{\infty}$. 
Then three generic fibers are all isomorphic. Hence, the formal 
completion of the diagram along the closed fibers (over $0 \in 
T_{\infty}$) gives two ``formal modifications" in the sense of 
[Ar 3].  By [Ar 3], there exists a 
diagram of algebraic spaces of finite type over $\mathbf{C}$: 
$$ \mathcal{X} \to \mathcal{Z} \leftarrow \mathcal{X}'$$ 
which extends such formal modifications. Take the closed fibers 
of this diagram over $0 \in T$. Then $\mathcal{X}_0$ is 
non-singular since $\mathcal{X}_0$ is etale over $X$. On the 
other hand, $\mathcal{X}_0$ and $\mathcal{X}'_0$ both have 
locally trivial deformations to a common affine variety 
$\mathcal{Z}_t$ $(t \ne 0)$ by the diagram. Therefore, 
$\mathcal{X}'_0$ is non-singular. Since $\mathcal{X}'_0$ is 
etale over $X'$ and the image of this etale map contains 
$(f')^{-1}(p)$ by the construction, $X'$ 
is non-singular at every point $q \in X'$ with $f'(q) = p$. 
Since $p \in Y$ is an arbitrary closed point, $X'$ is non-singular.     
\vspace{0.2cm}

\section{Examples}

{\bf Example 32}. 
Assume that $\mathcal{O}_x \subset 
\mathfrak{sl}(n)$ is the orbit containing an nilpotent 
element $x$ of Jordan type 
$\mathbf{d} := [d_1, ..., d_k]$. 
Let $[s_1, ..., s_m]$ be the 
dual partition of $\mathbf{d}$, that is, 
$s_i := \sharp \{j; d_j \geq i\}$. 
Let $P \subset SL(n)$ be the parabolic 
subgroup of flag type $(s_1, ..., s_m)$. 
Define $F : = 
SL(n)/P$. Note that $h^1(F, \Omega^1_F) = m-1$.  
Let $$ \tau_1  \subset \cdots \subset  \tau_{m-1} \subset
\cit^n \otimes_{\cit}\mathcal{O}_{F_{\sigma}}$$ 
be the universal subbundles on $F_{\sigma}$. 
A point of the cotangent bundle $T^*F$ of $F$ 
is expressed as a pair 
$(p, \phi)$ of $p \in F$ and $\phi \in 
\mathrm{End}(\cit^n)$ such that 
$$ \phi(\cit^n) \subset \tau_{m-1}(p), \cdots,  
\phi(\tau_2(p)) \subset 
\tau_1(p), \phi(\tau_1(p)) = 0.$$ 
The Springer resolution 
$$ s: T^*F \to \bar{\mathcal{O}_x}$$ 
is defined as $s_{\sigma}((p, \phi)) := \phi$. 
Therefore, $T^*F$ is a smooth convex symplectic 
variety. 
Let $\mathcal{E}$ be the universal extension 
of $\mathcal{O}_F$ by $\Omega^1_F$. In other 
words, $\mathcal{E}$ fits in the exact sequence 
$$ 0 \to \Omega^1_F \to \mathcal{E} 
\stackrel{\eta}\to 
\mathcal{O}_{F}^{m-1} \to 0,$$ and 
the induced map $H^0(F, \mathcal{O}_F^{m-1}) \to 
H^1(F, \Omega^1_F)$ is an isomorphism. 
The locally free sheaf $\mathcal{E}$ can be 
constructed as follows. 
For $p \in F_{\sigma}$, we can choose a basis 
of $\mathbf{C}^n$ such that $\Omega^1_F(p)$ 
consists of the matrices of the following form 
 
$$ 
\begin{pmatrix} 
0 & * & \cdots &* \\
0 & 0 &  \cdots & *\\
\cdots & & & \cdots \\
0 & 0 &  \cdots& 0 \\ 
\end{pmatrix}.
$$  

Then $\mathcal{E}(p)$ is the vector 
subspace of $\mathfrak{sl}(n)$ consisting of 
the matrices $A$ of the following 
form
$$ 
\begin{pmatrix} 
a_1 & * & \cdots & * \\
0 & a_2 &\cdots  & *\\
\cdots & & & \cdots \\
0 & 0 & \cdots  &a_m \\ 
\end{pmatrix},$$
where $a_i := a_i I_{s_i}$ and $I_{s_i}$ is 
the identity matrix of 
the size $s_i \times s_i$. 
Since $A \in \mathfrak{sl}(n)$, 
$\Sigma_i s_i a_i = 0$. Here we define the 
map $\eta (p): \mathcal{E}(p) 
\to {\cit}^{\oplus m-1}$ as 
$\eta (p)(A) := (a_1, a_2, \cdots, a_{m-1}).$ 
Let $\mathbf{A}(\mathcal{E}^{*}) := 
\mathrm{Spec}_F \mathrm{Sym}^{\cdot}(\mathcal{E}^*)$ 
be the vector bundle over $F$ associated with 
$\mathcal{E}$. Then we have an exact sequence of 
vector bundles 
$$ 0 \to T^*F \to \mathbf{A}(\mathcal{E}^*) \to 
F \times \mathbf{C}^{n-1} \to 0.$$ 
The last homomorphism in the exact sequence 
gives a map 
$$ f: \mathbf{A}(\mathcal{E}^*) \to \mathbf{C}^{m-1},$$ 
where $f^{-1}(0) = T^*F$. 
This is a universal Poisson deformation of the 
Poisson scheme $T^*F$ (with respect to the canonical 
symplectic 2-form). In fact, by Proposition 1.4.14 
of [C-G], there is a relative symplectic 2-form 
of $f$ extending the canonical symplectic 2-form on 
$T^*F$; hence $f$ is a Poisson deformation. 
Let $p: T^*F \to F$ be the canonical projection. 
Then we have a commutative diagram of exact sequences: 

\begin{equation}
\begin{CD}
0 @>>> p^*\Omega^1_F @>>> 
p^*\mathcal{E} @>>> p^*\mathcal{O}_F^{m-1} 
@>>> 0 \\ 
@. @VVV @VVV @VVV @. \\ 
0 @>>> \Theta_{T^*F} @>>> 
\Theta_{\mathbf{A}(\mathcal{E}^*)}\vert_{T^*F}
 @>>> 
N_{T^*F/\mathbf{A}(\mathcal{E}^*)} @>>> 0 
\end{CD}
\end{equation}
 
Let $T$ be the tangent space of the base space $\mathbf{C}^{m-1}$ 
of $f$ at $0 \in \mathbf{C}^{m-1}$. The Kodaira-Spencer map $\theta_f$
of $f$ is given as the composite 
$$ T \to H^0(T^*F, N_{T^*F/\mathbf{A}(\mathcal{E}^*)}) 
\to H^1(T^*F, \Theta_{T^*F}).$$ 
On the other hand, if one identifies $T$ with 
$H^0(F, \mathcal{O}_F^{m-1})$, then one has a map 
$$ T \cong H^0(F, \mathcal{O}_F^{m-1}) \to 
H^1(F, \Omega^1_F).$$ By the construction, the Kodaira-Spencer map 
is factored by this map: 
$$  T \to H^1(F, \Omega^1_F) \to H^1(T^*F, \Theta_{T^*F}). $$ 
The first map is an isomorphism by the definition of 
$\mathcal{E}$. The second map is an injection. In fact, let 
$S \subset T^*F$ be the zero section. Then $N_{S/T^*F} \cong 
\Omega^1_S$ and the composite $H^1(F, \Omega^1_F) 
\to H^1(T^*F, \Theta_{T^*F}) \to H^1(S, \Omega^1_S)$ is an 
isomorphism. Therefore, the Kodaira-Spencer map 
$\theta_f$ is an injection. Since $f$ is a Poisson deformation 
of $T^*F$, the Kodaira-Spencer map $\theta_f$ is factored by 
the ``Poisson Kodaira-Spencer map " $\theta^P_f$: 
$$ T \stackrel{\theta^P_f}\to H^2(T^*F, \mathbf{C}) \to 
H^1(T^*F, \Omega^1_{T^*F}).$$ Hence $\theta^P_f$ is also 
injective. Since $\dim T = h^2(T^*F, \mathbf{C}) = m-1$, 
$\theta^P_f$ is actually an isomorphism.  

More generally, let $G$ be a complex simple Lie 
group and $\mathcal{O}$ be a nilpotent orbit 
in $\g := \mathrm{Lie}(G)$. Assume that the closure 
$\bar{\mathcal{O}}$ of $\mathcal{O}$ admits a Springer 
resolution $\mu: T^*(G/P) \to \bar{\mathcal{O}}$ for 
some parabolic subgroup $P \subset G$. One can identify 
$T^*(G/P)$ with the adjoint bundle $G \times^P n(P)$, 
where $n(P)$ is the nilradical of $\mathfrak{p} := \mathrm{Lie}(P)$.      
Let $r(P)$ be the solvable radical of 
$\mathfrak{p}$ and let $m(P)$ be the Levi-factor 
of $\mathfrak{p}$. We put $\mathfrak{k}(P) := \g^{m(P)}$, 
where $$\g^{m(P)} := \{x \in \g; [x,y] = 0, y \in m(P)\}.$$ 
In [Na 4, \S 7], we have defined a flat deformation of 
$T^*(G/P)$ as 
$$ G \times^P r(P) \to \mathfrak{k}(P). $$ 
Then this becomes a universal Poisson deformation of 
$T^*(G/P)$. 
\vspace{0.2cm}

{\bf Example 33}. Let $\mathcal{O}$ be the nilpotent 
orbit in $\mathfrak{sl}(3)$ of Jordan type $[1,2]$. 
Then the closure $\bar{\mathcal{O}}$ has two different 
Springer resolutions 
$$ T^*(SL(3)/P_{1,2}) \to \bar{\mathcal{O}} \leftarrow 
T^*(SL(3)/P_{2,1}), $$ 
where $P_{1,2}$ and $P_{2,1}$ are parabolic subgroups 
of $SL(3)$ of flag type $(1,2)$ and $(2,1)$ respectively. 
We put $X^+ := T^*(SL(3)/P_{1,2})$ and 
$X^- := T^*(SL(3)/P_{2,1})$. Then $X^+$ and $X^-$ are both 
isomorphic to the cotangent bundle of $\mathbf{P}^2$. 
We call the diagram a Mukai flop. Let $G \subset SL(3)$ 
be the finite group of order $3$ generated by
     
$$ 
\begin{pmatrix} 
1 & 0 & 0 \\
0 & \zeta &  0 \\
0 & 0 & \zeta^2 \\ 
\end{pmatrix},
$$

where $\zeta$ is a primitive 3-rd root of unity. 
Then $G$ acts on $\bar{\mathcal{O}}$ by the adjoint 
action. Since the Kostant-Kirillov 2-form on 
$\mathcal{O}$ is $SL(3)$-invariant, the $G$-action 
lifts to symplectic actions on $X^+$ and $X^-$. 
Divide $\bar{\mathcal{O}}$, $X^+$ and $X^-$ by 
these $G$-action, we get the diagram of a singular 
flop: 
$$ X^+/G \to \bar{\mathcal{O}}/G \leftarrow 
X^-/G. $$ 
Here $X^+/G$ (resp. $X^-/G$) has 3 isolated 
quotient (terminal) singularities. This is a typical 
example of Corollary 25.

\section{Appendix} 
 
(A.1) Let $Y := \mathrm{Spec}\: R$ be an 
affine variety over {\bf C}. A $\mathbf{C}^*$-action 
on $Y$ is a homomorphism $\mathbf{C}^* \to 
\mathrm{Aut}_{\mathbf{C}}(R)$ induced from 
a $\mathbf{C}$-algebra homomorphism 
$$ R \to R \otimes_{\mathbf{C}}\mathbf{C}[t, 1/t].$$ 
More exactly, a $\mathbf{C}$-valued point 
of $\mathbf{C}^*$ is regarded as a surjection of 
$\mathbf{C}$-algebras: 
$$\sigma : \mathbf{C}[t, 1/t] \to \mathbf{C}.$$ 
Then 
$$ R \to R \otimes_{\mathbf{C}}\mathbf{C}[t, 1/t] 
\stackrel{id \otimes \sigma}\to R $$ 
is an element of $\mathrm{Aut}_{\mathbf{C}}(R)$.  
If this correspondence gives a homomorphism 
$\mathbf{C}^* \to \mathrm{Aut}_{\mathbf{C}}(R)$, 
we say that $R$ (or $Y$) has a $\mathbf{C}^*$-action.  
A $\mathbf{C}^*$-action on $Y$ is 
called good if there is a maximal ideal $m_R$ of $R$ 
fixed by the action and if $\mathbf{C}^*$ has only 
positive weight on $m_R$.   
Next let us consider the case where $Y$ is the 
spectrum of a local complete $\mathbf{C}$-algebra $R$ with 
$R/m_R = \mathbf{C}$. A $\mathbf{C}^*$-action on $Y$ is then a 
homomorphism $\mathbf{C}^* \to 
\mathrm{Aut}_{\mathbf{C}}(R)$ induced from 
a $\mathbf{C}$-algebra homomorphism 
$$ R \to R \hat{\otimes}_{\mathbf{C}} \mathbf{C}[t, 1/t],$$ 
where $R \hat{\otimes}_{\mathbf{C}}\mathbf{C}[t, 1/t]$ 
is the completion of $R \otimes_{\mathbf{C}}\mathbf{C}[t, 1/t]$ 
with respect to the ideal $m_R(R \otimes_{\mathbf{C}}
\mathbf{C}[t, 1/t])$. Then the $\mathbf{C}^*$-action is 
called good if $\mathbf{C}^*$ has only positive weight 
on the maximal ideal of $R$. 
\vspace{0.12cm}

{\bf Lemma (A.2)}. {\em Let $(A, m)$ be a complete 
local $\mathbf{C}$-algebra with a good 
$\mathbf{C}^*$-action. Assume that $A/m = 
\mathbf{C}$. Let $R$ be the {\bf C}-vector subspace of $A$ 
spanned by all eigen-vectors in $A$. Then 
$R$ is a finitely generated $\mathbf{C}$-algebra 
with a good $\mathbf{C}^*$-action. Moreover, 
$\hat{R} = A$ where $\hat{R}$ is the completion 
of $R$ with the maximal ideal $m_R$.}  
\vspace{0.12cm}

{\em Proof}. 
Since $A/m^k$ ($k \geq 1$) are finite dimensional 
$\mathbf{C}$-vector spaces,  
they are direct sum of eigen-spaces with non-negative weights: 
$$A/m^k = \oplus_w (A/m^k)^w.$$ 
The natural maps $(A/m^k)^w \to (A/m^{k-1})^w$ 
are surjections for all $k$. Since $m/m^2$ is also 
decomposed into the direct sum of eigen-spaces, one 
can take eigen-vectors 
$\bar{\phi}_i$, $(i = 1,2, ..., l)$ as 
a generator of $m/m^2$. We put $w_i := wt(\bar{\phi}_i) > 0$. 
One can lift $\bar{\phi}_i$ to 
$\phi_i \in \lim (A/m^k)^{w_i}$ by the surjections above. 
Since $A$ is complete, $\phi_i \in A$ and $wt(\phi_i) = w_i$. 
Put $w_{min} := min\{w_1, ..., w_l\} > 0$. 
We shall prove that $R = \mathbf{C}[\phi_1, ..., \phi_l]$. 
Let $\psi \in A$ be an eigen-vector with weight $w$. 
Take an integer $k_0$ such that $\psi \in m^{k_0}$ and 
$\psi \notin m^{k_0+1}$. Since every element of 
$m^{k_0}/m^{k_0+1}$ can be written as a homogenous 
polynomial of $\phi = (\bar{\phi}_1, ..., 
\bar{\phi}_l)$ of degree $k_0$, we see that 
$$ \psi \equiv f_{k_0}(\phi_1, ..., \phi_l) \: 
(\mathrm{mod} \: m^{k_0 +1}) $$ 
for some homogenous polynomial $f_{k_0}$ of 
degree $k_0$. We continue the similar approximation 
by replacing $\psi$ with $\psi - f_{k_0}(\phi)$. 
Finally, for any given $k$, we have an 
approximation 
$$ \psi \equiv f_{k_0}(\phi) + ... + f_{k-1}(\phi) 
\: 
(\mathrm{mod} \:  m^k).$$ 
Assume here that $k > w/w_{min}$.  
We set $\psi' := \Sigma_{k_0 \leq i \leq k-1}f_i(\phi)$.  
Assume that $\psi - \psi' 
\in m^r$ and $\psi - \psi' \notin 
m^{r+1}$ with some $r \geq k$. 
Since $\psi - \psi'$ has weight $w$, 
$[\phi - \phi'] \in m^r/m^{r+1}$ also has 
weight $w$. On the other hand, every non-zero eigen-vector 
in $m^r/m^{r+1}$ has weight at least $rw_{min}$. 
Hence $w \geq rw_{min}$, but this contradicts 
that $r \geq k > w/w_{min}$. 
Therefore, $\psi = \psi' \; \mathrm{mod} \; m^r$ 
for any $r$. Thus, 
$$ \psi = f_{k_0}(\phi) + ... + f_{k-1}(\phi).$$ 
This implies that $R = \mathbf{C}[\phi_1, ..., \phi_l]$. 
Let $m_R \subset R$ be the maximal ideal generated by 
$\phi_i$'s.  
Let $R_k$ be the {\bf C}-vector subspace of $m^k \;(\subset A)$ 
spanned by the 
eigen-vectors. The argument above shows that 
$R_k = (m_R)^k$. Since $R_k  = 
\oplus_w \lim(m^k/m^{k+i})^w$, we conclude 
that $(m_R)^k = \oplus_w \lim (m^k/m^{k+i})^w$. 
We now have 
$$ R/(m_R)^k = \oplus_w \lim (A/m^i)^w/\oplus_w 
\lim (m^k/m^{k+i})^w = $$ 
$$ \oplus_w \{\lim (A/m^i)^w/\lim (m^k/m^{k+i})^w\} 
= \oplus_w (A/m^k)^w = A/m^k. $$ 
Here the 2-nd last equality holds because 
$\{(m^k/m^{k+i})^w\}_i$ satisfies the Mittag-Leffler 
condition. This implies that 
$\hat{R} = A$. 
\vspace{0.2cm}

(A.3) Let $R$ be a integral domain 
finitely generated over {\bf C} or a 
complete local {\bf C}-algebra with 
residue field {\bf C}. Assume that 
$R$ has a good $\mathbf{C}^*$-action. 
Let $M$ be a finite $R$-module. We say that 
$M$ has an equivariant 
$\mathbf{C}^*$-action if, for 
each $\sigma \in \mathbf{C}^*$, we 
are given a map 
$$ \phi_{\sigma}: M \to M $$ 
with the following properties: 

(1) $\phi_{\sigma}$ is a {\bf C}-linear 
map. 

(2) $\phi_{\sigma}(rx) = \sigma (r) 
\phi_{\sigma}(x)$ for $r \in R$ and 
$x \in M$. 

(3) $\phi_{\sigma\tau} = 
\phi_{\sigma}\circ \phi_{\tau}$ 
for $\sigma, \tau \in \mathbf{C}^*$. 

(4) $\phi_{1} = id$. 
\vspace{0.15cm}  

We say that a non-zero element $x \in M$ 
is an eigen-vector 
if there exists an integer $w$ such that 
$\phi_{\sigma}(x) = \sigma^w x$ for all 
$\sigma \in G$.  

Let $M$ and $N$ be $R$-modules with 
equivariant $\mathbf{C}^*$-actions. 
Then an $R$-homomorphism $f: M \to N$ 
is an equivariant map if $f$ 
is compatible with both $\mathbf{C}^*$-actions. 
\vspace{0.15cm}

{\bf Lemma (A.4)}. {\em Let $A$ and $R$ be 
the same as Lemma (A.2). Let $M$ be a finite 
$A$-module with an equivariant $\mathbf{C}^*$-action. 
Define $M_R$ to be the {\bf C}-vector 
subspace of $M$ spanned by 
the eigen-vectors of $M$. Then $M_R$ is a 
finite $R$-module with an equivariant 
$\mathbf{C}^*$-action. Moreover, $M_R \otimes_R 
A = M$.} 
\vspace{0.12cm}

{\em Proof}. The idea is the same as Lemma (A.2). 
The finite dimensional $\mathbf{C}$-vector 
space $M/m^kM$ is the direct sum of eigen-spaces. 
Thus, for each weight $w$, 
$$ (M/m^kM)^w \to (M/m^{k-1}M)^w $$ 
is surjective. Let $\bar{x}_i$ ($i = 1, ..., r$) 
be the eigen-vectors which generate $M/mM$. 
We lift $\bar{x}_i$ to $x_i \in \hat{M}$ by 
the surjections above. Since $\hat{M} = M$, 
$x_i$ are eigen-vectors of $M$. We set 
$u_i := wt(x_i)$ and $u_{min} := 
min \{u_1, ..., u_r\}$. We shall prove that 
$M_R$ is generated by $\{x_i\}$ as an 
$R$-module. Let $y \in M$ be an eigen-vector 
with weight $u$. Take an integer $k_0$ such 
that $y \in m^{k_0}M$ and $y \notin m^{k_0 + 1}M$. 
Let us consider the surjection 
$$ m^{k_0}/m^{k_0 + 1} \otimes M/mM 
\to m^{k_0}M/m^{k_0 +1}M.$$ 
As in the proof of Lemma (A.2), every element 
of $m^{k_0}/m^{k_0 + 1}$ is written as a 
homogenous polynomial of $\phi_1$, ..., $\phi_l$ 
of degree $k_0$, where $\phi_i$ are certain 
eigen-vectors contained in $R$.  
We put $w_{min} := min \{wt(\phi_1), ..., 
wt(\phi_l)\} > 0$. On the other 
hand, $M/mM$ is spanned by $x_i$'s. 
Thus, 
$$ y \equiv \Sigma r_i(\phi)x_i \; \mathrm{mod} 
\; m^{k_0+1}M,$$ 
where $r_i$ are homogenous polynomials of degree $k_0$ 
such that $wt(r_i(\phi)) + u_i = u$. 
We write $g_{k_0}$ for the right-hand side for short. 
Now, we have $y - g_{k_0} \in m^{k_0 + 1}M$. 
By replacing $y$ with $y - g_{k_0}$, we continue the 
similar approximation. Finally, for any $k$, we 
have an approximation: 
$$ y \equiv g_{k_0} + g_{k_0 + 1} + ... + g_{k-1} \; 
\mathrm{mod} \; m^kM. $$ 
By the construction, $y - \Sigma_{k_0 \leq i 
\leq k-1} g_i$ is an eigen-vector with weight 
$u$. In particular, $[y - \Sigma_{k_0 \leq i \leq k-1}g_i] 
\in m^kM/m^{k+1}M$ has weight $u$. 
On the other hand, every non-zero eigen-vector 
of $m^kM/m^{k+1}M$ has weight at least 
$kw_{min} + u_{min}$. If we take $k$ sufficiently large, 
then $kw_{min} + u_{min} > u$. This implies that 
$[y - \Sigma_{k_0 \leq i \leq k-1}g_i] = 0$. 
Repeating the same, we conclude that, for any $r > k$, 
$$ y \equiv \Sigma_{k_0 \leq i \leq k-1}g_i \; \mathrm{mod} 
\; m^rM.$$ 
This implies that, in $M$,  $$y = \Sigma_{k_0 \leq i 
\leq k-1}g_i.$$ 
Thus, $M_R$ is generated by $\{x_i\}$ as an 
$R$-module.      
Let $M_k$ be the subspace of $M$ spanned by 
the eigen-vectors in $m^kM$. Then the argument 
above shows that $M_k = (m_R)^kM_R$. 
Then, by the same argument as Lemma (A.2), 
$M_R/(m_R)^kM_R = M/m^kM$; hence $M_R 
\otimes_R A = M$. In order to prove that 
$M_R$ has an equivariant $\mathbf{C}^*$-action, 
we have to check that $\phi_{\sigma}(M_R) 
\subset M_R$ for all $\sigma \in 
\mathbf{C}^*$ (cf. (A.3)); but it is 
straightforward.  
\vspace{0.2cm}

{\bf Proposition (A.5)}. {\em Let $A$ be 
a local complete $\mathbf{C}$-algebra with 
residue field {\bf C} and with a good 
$\mathbf{C}^*$-action. Let $f: X \to 
\mathrm{Spec}(A)$ be a $\mathbf{C}^*$-eqivariant 
projective morphism 
and let $L$ be an $f$-ample, $\mathbf{C}^*$-linearized 
line bundle. Let $R$ be the same as Lemma (A.2). 
Then there is a $\mathbf{C}^*$-equivariant 
projective morphism $f_R : X_R \to \mathrm{Spec}(R)$ 
and a $\mathbf{C}^*$-linearized, $f_R$-ample line bundle 
$L_R$, such that 
$X_R \times_{\mathrm{Spec}(R)}\mathrm{Spec}(A) \cong X$ 
and $L_R \otimes_R A \cong L$.} 
\vspace{0.2cm}

{\em Proof}. We put $A_i := \Gamma (X, L^{\otimes i})$ for 
$i \geq 0$. Then, $X = \mathrm{Proj}_A \oplus_{i \geq 0}A_i.$  
If necessary, by taking a suitable 
multiple $L^{\otimes m}$, we may assume that 
$A_* := \oplus_{i \geq 0}A_i$ is generated by $A_1$ as an 
$A_0 (=A)$-algebra. By Lemma (A.4), we take a finite 
$R$-module $A_{i,R}$ such that $A_{i,R}\otimes_R A 
= A_i$. The multiplication map 
$A_i \otimes_{A_0} A_j \to A_{i+j}$ induces 
a map $A_{i,R} \otimes_R A_{j,R} \to A_{i+j, R}$; 
hence $(A_*)_R := \oplus_{i \geq 0}A_{i,R}$ becomes 
a graded $R$-algebra. We shall check that $(A_*)_R$ 
is a finitely generated $R$-algebra. In order to do 
this, we only have to prove that $(A_R)_*$ is 
generated by $A_{1,R}$ as an $R$-algebra since 
$A_{1,R}$ is a finite $R$-module. 
Let us consider the $n$-multiplication map 
$$ m_n: 
\overbrace{A_{1,R} \otimes_R ... \otimes_R A_{1,R}}^n 
\to A_{n,R}.$$ 
Let $M$ be the cokernel of this map. 
Since $m_n$ is a $\mathbf{C}^*$-equivariant map, 
$M$ is a finite $R$-module with an equivariant 
$\mathbf{C}^*$-action. Taking the tensor 
product $\otimes_R A$ with $m_n$, we get the 
$n$-multiplication map for $A_*$; but this is 
surjective by the assumption. Therefore, 
$\hat{M} := M \otimes_R A = 0$. The support of 
$M$ is a closed subset of $\mathrm{Spec}(R)$, stable 
under the $\mathbf{C}^*$-action. Since $\hat{M} = 0$, 
this closed subset does not contain the origin 
$0 \in R$; hence it must be empty because $R$ has a 
good $\mathbf{C}^*$-action. Finally it is clear 
that $(A_R)_* \otimes_R A = A_*$ by the construction. 
\vspace{0.2cm}

{\bf Proposition (A.6)}. {\em Let $f: X \to 
\mathrm{Spec}(A)$ and $f_R: X_R \to \mathrm{Spec}(R)$ 
be the same as Lemma (A.5). Let $F$ be a coherent 
sheaf of $X$ with a $\mathbf{C}^*$-linearization. 
Then there is a $\mathbf{C}^*$-linearized coherent 
sheaf $F_R$ of $X_R$ such that $F_R \otimes_R A 
= F$.}  
\vspace{0.2cm}

{\em Proof}. We put $\mathcal{O}_X(1) := 
\widetilde{(\oplus_{i \geq 0}A_i)[1]}$. Then 
the coherent sheaf $F$ can be written as 
$$ F = \widetilde{\oplus_{i \geq 0} \Gamma (X, F(i))}.$$ 
Let us write $M_i$ for $\Gamma (X, F(i))$. 
By Lemma (A.4), there is a finite $R$-module 
$M_{i,R}$ such that $M_{i,R} \otimes_R A 
= M_i$. We define 
$$ F_R := \widetilde{\oplus_{i \geq 0} M_{i,R}}. $$ 
We shall prove that $(M_*)_R := 
\oplus_{i \geq 0}M_{i,R}$ is a finite $(A_*)_R$-module. 
There is an integer $n_0$ such that, for 
any $i \geq n_0$, and for any $j \geq 0$, 
the multiplication map 
$$ \overbrace{A_1 \otimes_{A_0} ... \otimes_{A_0} A_1}^j 
\otimes_{A_0} M_i \to M_{i+j} $$ 
is surjective. For the same $i$, $j$, let us 
consider the $R$-linear map  
$$ \overbrace{A_{1,R} \otimes_R ... \otimes_R A_{1,R}}^j 
\otimes_R 
M_{i,R} \to M_{i+j, R}.$$ 
Let $N$ be the cokernel of this map. 
Since this $R$-linear map is compatible with the 
$\mathbf{C}^*$-action on $R$, $N$ is a finite 
$R$-module with an equivariant $\mathbf{C}^*$-action. 
By the choice of $i$ and $j$, $\hat{N} := N \otimes_R A$ 
is zero. This implies that $N = 0$.      
\vspace{0.2cm}

{\bf Proposition (A.7)} 
{\em Let $Y$ be an affine symplectic 
variety.   
Assume that $Y$ has a good $\mathbf{C}^*$-action 
with a fixed point $0 \in Y$.   
Assume that, in the analytic category, $Y^{an}$ admits a 
crepant, projective, partial resolution   
$\bar{f}: \mathcal{X} \to Y^{an}$ 
such that  $\mathcal{X}$ has only terminal singularities. 
Then, in the algebraic category, $Y$ admits a crepant, 
projective, partial resolution 
$f : X \to Y$ such that $X^{an} = \mathcal{X}$ and 
$f^{an} = \bar{f}$.} 
\vspace{0.2cm}

{\em Proof}. 
(STEP 1): We shall prove that the $\mathbf{C}^*$-action of 
$Y^{an}$ lifts to $\mathcal{X}$.  
Since $Y^{an}$ is symplectic, one can take a 
closed subset $\Sigma$ of $Y^{an}$, stable under the 
$\mathbf{C}^*$-action and $\mathrm{codim}(\Sigma 
\subset Y^{an}) \geq 4$, such that the singularities 
of $Y^{an} - \Sigma$ are local trivial 
deformations of two dimensional rational double points.  
We put $Y_0 := Y^{an} - \Sigma$. Since $\bar{f}$ is 
the minimal resolution over $Y_0$, the $\mathbf{C}^*$-action 
on $Y_0$ extends to $\mathcal{X}_0 := \bar{f}^{-1}(Y_0)$. 
Note that, in $\mathcal{X}$, 
$\mathcal{X} - \mathcal{X}_0$ has codimension at least two by 
the semi-smallness of $\bar{f}$ ([Na 3]). 
The $\mathbf{C}^*$-action defines a holomorphic map 
$$\sigma^0 : \mathbf{C}^* \times \mathcal{X}_0 
\to \mathcal{X}_0,$$ 
and this extends to a meromorphic map 
$$\sigma: \mathbf{C}^* \times \mathcal{X} --\to \mathcal{X}.$$ 
Let us prove that $\sigma_t : \mathcal{X} --\to \mathcal{X}$, which 
is an isomorphism in codimension one,  
is actually an isomorphism everywhere for each $t \in \mathbf{C}^*$. 
Let $\mathcal{L}$ be an $\bar{f}$-ample line bundle on 
$\mathcal{X}$. We put $\mathcal{L}^0_t := 
(\sigma^0)^*\mathcal{L}\vert_{\{t\} \times \mathcal{X}_0}$. 
Since $\mathrm{Pic}(\mathcal{X}_0)$ is discrete, 
$\mathcal{L}^0_t$ are all isomorphic to $\mathcal{L}\vert_{X_0}$. 
Since $\mathcal{L}\vert_{X_0}$ extends to the line bundle 
$\mathcal{L}$ on $\mathcal{X}$, $(\sigma^0)^*\mathcal{L}$ extends to a 
line bundle on $\mathbf{C}^* \times \mathcal{X}$, say 
$\sigma^*\mathcal{L}$ by abuse of notation. The line bundle 
$\mathcal{L}_t := \sigma^*\mathcal{L}\vert_{\{t\} \times \mathcal{X}}$ 
coincides with the proper transform of $\mathcal{L}$ 
by $\sigma_t$. Since $\mathcal{L}_0 (= \mathcal{L})$ is 
$\bar{f}$-ample and $\mathrm{Pic}(\mathcal{X}/Y^{an}) 
:= \mathrm{Pic}(\mathcal{X})/\bar{f}^*\mathrm{Pic}(Y^{an})$ 
is discrete, $\mathcal{L}_t$ are all $\bar{f}$-ample. 
This implies that $\sigma_t$ are all isomorphisms and 
$\sigma$ is a holomorphic map. One can check that $\sigma$ 
gives a $\mathbf{C}^*$-action because it already becomes a 
$\mathbf{C}^*$-action on $\mathcal{X}_0$. 
\vspace{0.12cm}

(STEP 2): Let $Y_n$ be the $n$-th infinitesimal 
neighborhood of $Y^{an}$ at $0$, which becomes an  
affine scheme with a unique point $0$. We put   
$\mathcal{X}_n := \mathcal{X} \times_{Y^{an}} Y^{an}_n$. 
By GAGA, there are projective schemes $X_n$ over $Y_n$ such 
that $(X_n)^{an} = \mathcal{X}_n$. Fix an $\bar{f}$-ample 
line bundle $\mathcal{L}$ on $\mathcal{X}$. Again by 
GAGA, it induces line bundles $L_n$ on $X_n$. 
The $\mathbf{C}^*$-action on $\mathcal{X}$ induces a 
$\mathbf{C}^*$-action on $\mathcal{X}_n$ for each $n$. 
This action induces an {algebraic} $\mathbf{C}^*$-action of $X_n$. 
In fact, the $\mathbf{C}^*$-action of $\mathcal{X}$ 
originally comes from an algeraic $\mathbf{C}^*$-action 
on $Y$, the holomorphic action map 
$$ \mathbf{C}^* \times \mathcal{X} \to \mathcal{X} $$ 
extends to a meromorphic map 
$$ \mathbf{P}^1 \times \mathcal{X} --\to \mathcal{X}.$$ 
Thus, the holomorphic action map 
$$ \mathbf{C}^* \times \mathcal{X}_n \to \mathcal{X}_n $$ 
extends to a meromorphic map 
$$ \mathbf{P}^1 \times \mathcal{X}_n --\to \mathcal{X}_n.$$ 
Thus, by GAGA, we have a rational map 
$$ \mathbf{P}^1 \times X_n --\to X_n $$ 
which resticts to an algebraic $\mathbf{C}^*$-action 
on $X_n$. Let us regard $\{X_n\}$ and $\{Y_n\}$ 
as formal schemes and $\{f_n : X_n \to Y_n\}$ as a projective 
equivariant morphism of formal schemes with 
$\mathbf{C}^*$-actions. Put $\hat{A} := 
\lim \mathcal{O}_{Y_n, 0}$ and 
$\hat{Y} := \mathrm{Spec}(\hat{A})$. Then, by 
[EGA III], Theoreme 5.4.5, the projective 
morphism 
of formal schemes can be algebraized to 
a projective equivariant morphism of schemes  
$$ \hat{f}: \hat{X} \to \hat{Y}.$$ 
The affine scheme $\hat{Y}$ admits a $\mathbf{C}^*$-action 
coming from the original $\mathbf{C}^*$-action on $Y$, 
which is compatible with the $\mathbf{C}^*$-action on 
$\{Y_n\}$. The $\mathbf{C}^*$-action on $\{X_n\}$ 
naturally lifts to $\hat{X}$ in such a way that 
$\hat{f}$ becomes a $\mathbf{C}^*$-equivariant morphism. 
In fact, let $$\sigma_n : \mathbf{C}^* \times X_n 
\to X_n$$ be the $\mathbf{C}^*$-action on $X_n$. 
Let us consider the morphism (of formal schemes): 
$$ id \times \{\sigma_n\}: \mathbf{C}^* \times 
\{X_n\} \to \mathbf{C}^* \times \{X_n\}. $$ 
Here we regard the first factor (resp. the 
second factor) as a $\mathbf{C}^* \times \{Y_n\}$-
formal scheme  
by $id \times (\{f_n\} \circ \{\sigma_n\})$ 
(resp. $id \times \{f_n\}$). Then the morphism 
above is a $\mathbf{C}^* \times \{Y_n\}$-morphism. 
By [ibid, Theoreme 5.4.1], this morphism of 
formal schemes extends to a $\mathbf{C}^* 
\hat{\times} \hat{Y}$-morphism 
\footnote{$\hat{\times}$ means the formal 
product. Let $B$ be the completion of 
$\hat{A}[t, 1/t]$ by the ideal 
$m\hat{A}[t, 1/t]$ where $m \subset  
\hat{A}$ is the maximal ideal. Then 
$\mathbf{C}^* \hat{\times} \hat{Y} = 
\mathrm{Spec}(B)$. The scheme 
$\mathbf{C}^* \hat{\times} \hat{X}$ is 
defined as the fiber product of 
$\mathbf{C}^* \times \hat{X} \to 
\mathbf{C}^* \times \hat{Y}$ and 
$\mathbf{C}^* \hat{\times} \hat{Y} 
\to \mathbf{C}^* \times \hat{Y}$. }  
$$ \mathbf{C}^* \hat{\times} \hat{X} \to 
\mathbf{C}^* \hat{\times} \hat{X}, $$ 
where the first factor (resp. the second factor) 
is regarded as a $\mathbf{C}^* \hat{\times} \hat{Y}$-scheme 
by $id \hat{\times} \hat{f}\circ\hat{\sigma}$ 
(resp. $id \hat{\times} \hat{f}$). 
The extended morphism 
gives a $\mathbf{C}^*$-action 
$$ \mathbf{C}^* \hat{\times} \hat{X} \to 
\mathbf{C}^* \hat{\times} \hat{X} \stackrel{p_2}\to 
\hat{X}.$$ 
Moreover, $\{L_n\}$ is algebraized to an $\hat{f}$-ample 
line bundle $\hat{L}$ on $\hat{X}$ ([ibid, Theoreme 5.4.5]). 
Since $\mathcal{L}$ is fixed by the $\mathbf{C}^*$-action 
on $\mathcal{X}$, $\hat{L}$ is also fixed by the 
$\mathbf{C}^*$-action on $\hat{X}$. 
\vspace{0.2cm}

{\bf Lemma (A.8)}. 
{\em Let $\hat{f} : \hat{X} \to \hat{Y}$ be 
a $\mathbf{C}^*$-equivariant projective morphism 
where $\hat{Y} = 
\mathrm{Spec}(\hat{A})$ with a complete 
local $\mathbf{C}$-algebra $\hat{A}$ with 
$\hat{A}/m = \mathbf{C}$. 
Assume that $\hat{f}_*\mathcal{O}_{\hat{X}} = 
\mathcal{O}_{\hat{Y}}$. 
Let $\hat{L}$ be an $\hat{f}$-ample line 
bundle on $\hat{X}$ fixed by the $\mathbf{C}^*$-action. 
Then $\hat{L}^{\otimes m}$ can be  
$\mathbf{C}^*$-linearized for some 
$m > 0$. Moreover, in this case, any $\mathbf{C}^*$-fixed 
line bundle $M$ on $\hat{X}$ is $\mathbf{C}^*$-linearized 
after taking a suitable multiple of $M$.} 
\vspace{0.15cm}

{\em Proof}. We only have to deal with an $\hat{f}$-ample 
line bundle $\hat{L}$.  
In fact, let $M$ be an arbitrary line bundle 
on $\hat{X}$ fixed by the $\mathbf{C}^*$-action. 
Then $M \otimes \hat{L}^{\otimes r}$ 
becomes $\hat{f}$-ample for a sufficiently large 
$r$. If we could prove the lemma for $\hat{f}$-ample 
line bundles, then $M^{\otimes m} \otimes \hat{L}^{\otimes 
rm}$ is $\mathbf{C}^*$-linearized. Since $\hat{L}^{\otimes 
rm}$ is also $\mathbf{C}^*$-linearized, $M^{\otimes m}$ 
is $\mathbf{C}^*$-linearized. 
We assume that $\hat{L}$ is $\hat{f}$-very ample and  
$\hat{X}$ is embedded into 
$\mathbf{P}_{\hat{A}}(H^0(\hat{X}, \hat{L}))$ 
as a $\hat{Y}$-scheme, where $H^0(\hat{X}, \hat{L})$ 
is a free $\hat{A}$-module of finite rank, say $n$. 
Since $\mathbf{C}^*$ acts on $\hat{A}$, we regard 
$\mathbf{C}^*$ as a subgroup 
of the automorphism group of the $\mathbf{C}$-algebra 
$\hat{A}$. Let $\sigma \in \mathbf{C}^*$ and let $M$ be 
an $\hat{A}$-module. Then a $\mathbf{C}$-linear map 
$\phi : M \to M$ is called a 
twisted $\hat{A}$-linear map if there exists 
$\sigma \in \mathbf{C}^*$ and $\phi (ax) = 
\sigma (a) \phi (x)$ for $a \in \hat{A}$ and 
for $x \in M$. Now let us consider the case 
$M = H^0(\hat{X}, \hat{L})$, which is a 
free $\hat{A}$-module of rank $n$.     
We define 
$G(n, \hat{A})$ to be the group of all 
twisted $\hat{A}$-linear bijective maps 
from $H^0(\hat{X}, \hat{L})$ onto itself. 
One can define a surjective 
homomorphisms 
$G(n, \hat{A}) \to \mathbf{C}^*$ 
by sending $\phi \in G(n, \hat{A})$ to 
the associated twisting element $\sigma 
\in \mathbf{C}^*$. Note that this 
homorphism admits a canonical splitting 
$\iota : \mathbf{C}^* \to G(n, \hat{A})$ 
defined by $\iota(\sigma)(x_1, ..., x_n) 
:= (\sigma (x_1), ..., \sigma (x_n))$. 
There is an exact sequence 
$$ 1 \to GL(n, \hat{A}) \to G(n, \hat{A}) 
\to \mathbf{C}^* \to 1. $$ 
Let us denote by $\hat{A}^*$ the multiplicative 
group of units of $\hat{A}$. 
One can embed $\hat{A}^*$ diagonally in 
$GL(n, \hat{A})$; hence in $G(n, \hat{A})$. 
Then $\hat{A}^*$ is a normal subgroup 
of $G(n, \hat{A})$. 
The group $PG(n, \hat{A}) := G(n, \hat{A})/\hat{A}^*$ acts 
faithfully on $H^0(\hat{X}, \hat{L}) - \{0\}/\hat{A}^*$. 
On the other hand, define $S(n, \hat{A})$ to 
be the subgroup of $G(n, \hat{A})$ generated 
by $SL(n, \hat{A})$ and $\iota (\mathbf{C}^*)$ 
There are two exact sequences 
$$ 1 \to SL(n, \hat{A}) \to S(n, \hat{A}) \to 
\mathbf{C}^* \to 1,$$ and  
$$ 1 \to PGL(n, \hat{A}) \to PG(n, \hat{A}) 
\to \mathbf{C}^* \to 1.$$ 
Since $\hat{A}$ is a complete local ring 
and its residue is an algebraically 
closed field with characteristic $0$, the canonical map 
$SL(n, \hat{A}) \to PGL(n, \hat{A})$ is a 
surjection; hence the composed map $S(n, \hat{A}) \to 
G(n, \hat{A}) \to PG(n, \hat{A})$ is surjective.   

Let us start the proof. Note that $H^0(\hat{X}, \hat{L})-\{0\}
/\hat{A}^*$ is identified with the space of Cartier 
divisors whose associated line bundle is $\hat{L}$.
Since $\hat{L}$ is fixed by the $\mathbf{C}^*$-action, 
the $\mathbf{C}^*$ action on $\hat{X}$ induces a $\mathbf{C}^*$ 
action on $H^0(\hat{X}, \hat{L})-\{0\}/\hat{A}^*$. 
This action gives a splitting 
$$ \alpha : \mathbf{C}^* \to PG(n, \hat{A}) $$ 
of the exact sequence above. 
We want to lift the map $\alpha$ to $S(n, \hat{A})$. 
We put $H := \varphi^{-1}(\alpha(\mathbf{C}^*))$, 
where $\varphi : S(n, \hat{A}) \to 
PG(n, \hat{A})$ is the quotient map. Since 
$\mathrm{Ker}(\varphi) = \mu_n$, $H$ is an 
etale cover of $\mathbf{C}^*$. Now $H$ acts on 
$H^0(\hat{X}, \hat{L})$. The $H$ naturally acts on 
the $n$-th  symmetric product 
$\mathrm{S}^n(H^0(\hat{X}, \hat{L}))$, where 
$\mu_n$ acts trivially. Therefore, we get a 
$\mathbf{C}^*$-action on $S^n(H^0(\hat{X}, \hat{L}))$. 
This $\mathbf{C}^*$-action induces a $\mathbf{C}^*$-linearization 
of $\mathcal{O}_{\mathbf{P}(H^0(\hat{X}, \hat{L})}(n)$. 
Since $\hat{L}^{\otimes n}$ is the pull-back of this line 
bundle by the $\mathbf{C}^*$-equivariant embedding 
$\hat{X} \to \mathbf{P}(H^0(\hat{X}, \hat{L}))$, 
$\hat{L}^{\otimes n}$ has a $\mathbf{C}^*$-linearization.   
\vspace{0.15cm}
   
By the lemma above, $\hat{L}^{\otimes m}$ is $\mathbf{C}^*$-lineraized 
for some $m > 0$. 
Now one can write 
$$ \hat{X} = \mathrm{Proj}_{\hat{A}}\oplus_{n \geq 0}
\hat{f}_*\hat{L}^{\otimes nm}, $$ 
where each $\hat{f}_*\hat{L}^{\otimes nm}$ is an $\hat{A}$-module 
with $\mathbf{C}^*$-action. Since $Y$ has a 
good $\mathbf{C}^*$-action, there exists a projective 
$\mathbf{C}^*$-equivariant morphism $f: X \to Y$ 
such that $X \times_Y \hat{Y} = \hat{X}$ by 
Proposition (A.5).    
\vspace{0.12cm}

(STEP 3): We shall finally show that $X^{an} = \mathcal{X}$ 
and $f^{an} = \bar{f}$. 
The formal neighborhoods of $X^{an}$ and $\mathcal{X}$ 
along $f^{-1}(0)$ are the same.  
By [Ar], the bimeromorphic map $X^{an} --\to \mathcal{X}$ 
is an isomorphism over a small open neighborhood $U$ of 
$0 \in Y^{an}$. But, since $Y^{an}$ has a good 
$\mathbf{C}^*$-action and this action lifts to both 
$X^{an}$ and $\mathcal{X}$, the bimeromorphic map 
must be an isomorphism over $Y^{an}$.          
\vspace{0.2cm}

{\bf Proposition (A.9)}. 
{\em Let $Y = \mathrm{Spec}(A)$ be an affine variety with a 
good $\mathbf{C}^*$-action and let $f: X \to Y$ 
be a birational projective morphism with $X$ normal.  
Assume that 
$Y$ has only rational singularities, 
and $X$ is {\bf Q}-factorial. Then 
$X^{an}$ is {\bf Q}-factorial.} 
\vspace{0.15cm}

{\em Proof}. Let $g: Z \to Y$ be a $\mathbf{C}^*$-equivariant 
projective resolution. Let $0 \in Y$ be the fixed origin of 
the $\mathbf{C}^*$-action and let $\hat{Y} := \mathrm{Spec}(\hat{A})$ 
where $\hat{A}$ is the completion of $A$ at $0$. 
We put $\hat{Z} := Z \times_Y \hat{Y}$ and denote by $\hat{g}: 
\hat{Z} \to \hat{Y}$ the induced morphism. 
Since $Y$ has only rational singularities, 
$\mathrm{Pic}(Z^{an}) \cong H^2(Z^{an}, \mathbf{Z})$, which 
is discrete. Hence every element $\mathcal{L} 
\in \mathrm{Pic}(Z^{an})$ is fixed by the 
$\mathbf{C}^*$-action. Take an arbitrary line 
bundle $\mathcal{L}$. We shall prove that, 
for some $m > 0$, $\mathcal{L}^{\otimes m}$ 
comes from an algebraic line bundle.    
As in the proof of Proposition 26, $\mathcal{L}$ 
defines a line bundle $\hat{L}$ on $\hat{X}$. 
By Lemma (A.8), $\hat{L}^{\otimes m}$ is 
$\mathbf{C}^*$-linearized for some $m$. 
By Proposition (A.6), $\hat{L}^{\otimes m}$ 
extends to a $\mathbf{C}^*$-linearized line bundle 
$M$ on $Z$. By the construction, there is an 
open neighborhood $U$ of $0 \in Y^{an}$ such that 
$M^{an}\vert_{(g^{an})^{-1}(U)} \cong 
\mathcal{L}^{\otimes m}\vert_{(g^{an})^{-1}(U)}$. 
Since $Y^{an}$ has a good $\mathbf{C}^*$-action, 
one can assume that $H^2(Z^{an}, \mathbf{Z}) 
\cong H^2((g^{an})^{-1}(U), \mathbf{Z})$; this 
implies that  
$\mathrm{Pic}(Z^{an}) \cong \mathrm{Pic}((g^{an})^{-1}(U))$.
Thus, $M^{an} \cong \mathcal{L}^{\otimes m}$. 
Let us take a common resolution of $Z$ and $X$: 
$h_1: W \to Z$ and $h_2: W \to X$. 
Let $D$ be an irreducible (analytic) Weil 
divisor of $X^{an}$. Take an irreducible 
component $D'$ of $(h_2^{an})^{-1}(D)$ such that 
$(h_2^{an})(D') = D$. We put $\bar{D} := (h_1^{an})(D')$. 
We first assume that $\bar{D}$ is a divisor 
of $Z^{an}$. Then the line bundle  
$\mathcal{O}_{Z^{an}}(r\bar{D})$ becomes 
algebraic for some $r > 0$. Hence 
$\mathcal{O}_{W^{an}}(rD')$ is  
algebraic. Finally, the direct 
image $(h_2^{an})_*\mathcal{O}_{W^{an}}(rD')$ 
is algebraic, and its double dual is 
also algebraic. Thus we conclude that 
$\mathcal{O}_{X^{an}}(rD)$ is an algebraic 
reflexive sheaf of rank 1.  
We next assume that $\bar{D}$ is not a divisor. 
Then $D'$ is an exceptional divisor of $h_1$. 
In this case, $\mathcal{O}_{W^{an}}(D')$ is 
algebraic, and the same argument as the first case 
shows that $\mathcal{O}_{X^{an}}(D)$ is algebraic.  
\vspace{0.2cm}

{\bf Corollary (A.10)}. {\em Let $Y$ be an affine 
symplectic variety with a good $\mathbf{C}^*$-action. 
Then the following hold.} 
\vspace{0.12cm}

(i) {\em If $f: X \to Y$ is a {\bf Q}-factorial 
terminalization, then $X^{an}$ is {\bf Q}-factorial 
as an analytic space.} 
\vspace{0.12cm}

(ii){\em If $\bar{f}: \mathcal{X} \to Y^{an}$ is a 
{\bf Q}-factorial terminalization as an analytic 
space, then there is a projective birational 
morphism $f: X \to Y$ such that $X^{an} 
= \mathcal{X}$ and $f^{an} = \bar{f}$.}

\quad \\
\quad\\

Yoshinori Namikawa \\
Departement of Mathematics, 
Graduate School of Science, Osaka University, JAPAN \\
namikawa@math.sci.osaka-u.ac.jp

\end{document}